\documentclass{gtart}

\def\ifplaintex{\expandafter\ifx\csname documentclass\endcsname\relax}

\def\gtp{{\mathsurround=0pt\it $\cal G\mskip-2mu$eometry \&\ 
$\cal T\!\!$opology $\cal P\!$ublications}}  

\def\recd{{\small Received:\qua\receiveddate\ifx\reviseddate\relax
\else\qquad Revised:\qua\reviseddate\fi\par}} 

\def\lognumber#1{\def\thelognumber{#1}}
\def\volumenumber#1{\def\thevolumenumber{#1}}
\def\volumeyear#1{\def\thevolumeyear{#1}}
\def\papernumber#1{\def\thepapernumber{#1}}
\def\pagenumbers#1#2{\def\startpage{#1}\def\finishpage{#2}}
\def\published#1{\def\publishdate{#1}}

\def\received#1{\def\receiveddate{#1}}

\def\accepted#1{\def\accepteddate{#1}}

\long\def\asciiabstract#1{\long\def\theasciiabstract{#1}}

\let\\\par\let\thelognumber\relax\let\thevolumenumber\relax
\let\thepapernumber\relax\let\thevolumeyear\relax\let\startpage\relax
\let\finishpage\relax\let\publishdate\relax\let\receiveddate\relax
\let\reviseddate\relax\let\accepteddate\relax\let\theasciititle\relax
\let\theasciiauthors\relax
\let\theasciiabstract\relax

\let\theasciiemail\relax

\ifplaintex
\font\logobig=cmssbx10 scaled 3836
\font\logomed=cmssbx10 scaled 2557
\else
\font\logobig=cmssbx10 scaled 4200
\font\logomed=cmssbx10 scaled 2800
\fi

\long\def\makeagttitle{   
\count0=\startpage
\agt\hfill      
\hbox to 45truept{\vbox to 0pt{\vglue -13truept{\logomed A\kern -.37em{\logobig 
T}\kern -.38em G}\vss}\hss}
\break
{\small Volume \thevolumenumber\ (\thevolumeyear)
\startpage--\finishpage\nl
Published: \publishdate}

\vglue .25truein

{\parskip=0pt\leftskip 0pt plus
1fil\def\\{\par\smallskip}{\Large\bf\thetitle}\par\medskip} \vglue
0.05truein

{\parskip=0pt\leftskip 0pt plus 1fil\def\\{\par}{\sc\theauthors}
\par\medskip}%
 
\vglue 0.03truein 

{\small\leftskip 25truept\rightskip 25truept{\bf Abstract}\stdspace\theabstract

{\bf AMS Classification}\stdspace\theprimaryclass
\ifx\thesecondaryclass\relax\else; \thesecondaryclass\fi\par
{\bf Keywords}\stdspace \thekeywords\par}\vglue 7truept

}   

\ifplaintex
\font\phead=cmsl9 scaled 950
\font\pnum=cmbx10 scaled 913
\font\pfoot=cmsl9 scaled 950
\headline{\vbox to 0pt{\vskip -4.5mm\line{\small\phead\ifnum
\count0=\startpage ISSN 1472-2739 (on-line) 1472-2747 (printed)
\hfill {\pnum\folio}\else\ifodd\count0\def\\{ }%
\ifx\theshorttitle\relax\thetitle\else\theshorttitle\fi\hfill{\pnum\folio}
\else\def\\{ and }{\pnum\folio}\hfill\ifx\theshortauthors\relax\theauthors
\else\theshortauthors\fi\fi\fi}\vss}}
\footline{\vbox to 0pt{\vglue 0mm\line{\small\pfoot\ifnum\count0=\startpage
\copyright\ \gtp\hfill\else
\agt, Volume \thevolumenumber\ (\thevolumeyear)\hfill\fi}\vss}}
\else
\font=cmsl9 scaled 1050
\font\lnum=cmbx10 
\font=cmsl9 scaled 1050
\makeatletter
\def\@oddhead{{\small\lhead\ifnum\count0=\startpage ISSN 1472-2739 
(on-line) 1472-2747 (printed)\hfill {\lnum\number\count0}\else\ifodd\count0
\def\\{ }\ifx\theshorttitle\relax \thetitle \else\theshorttitle\fi\hfill
{\lnum\number\count0}\else\def\\{ and }{\lnum\number\count0}
\hfill\ifx\theshortauthors\relax 
\theauthors\else\theshortauthors\fi\fi\fi}}\def\@evenhead{\@oddhead}
\def\@oddfoot{\small\lfoot\ifnum\count0=\startpage\copyright\ \gtp\hfill\else
\agt, Volume \thevolumenumber\ (\thevolumeyear)\hfill\fi}
\def\@evenfoot{\@oddfoot}
\makeatother
\fi
\let\maketitlepage\makeagttitle
\let\makeshorttitle\maketitlepage
\let\maketitle\maketitlepage

\newwrite\gtoutfile
\long\gdef\makeheadfile{  
{\def\\{, }\def\s{ }
\immediate\openout\gtoutfile head.xxx
\immediate\write\gtoutfile{To: math@arxiv.org}
\immediate\write\gtoutfile{Subject: put OR rep NNNNN:ppppp}
\immediate\write\gtoutfile{--text follows this line--}
\immediate\write\gtoutfile{Proxy-for: \ifx\theasciiauthors\relax
\theauthors\else\theasciiauthors\fi\s<\ifx\theasciiemail\relax\theemail\else\theasciiemail\fi>}
\immediate\write\gtoutfile{\noexpand\\}
\immediate\write\gtoutfile{Authors: \ifx\theasciiauthors\relax
\theauthors\else\theasciiauthors\fi}
{\def\\{ }\immediate\write\gtoutfile{Title: \ifx\theasciititle\relax
\thetitle\else\theasciititle\fi}}
\immediate\write\gtoutfile{Subj-class: GT or SG, GR etc}
\immediate\write\gtoutfile{MSC-class: \theprimaryclass\ifx\thesecondaryclass\relax\else, \thesecondaryclass\fi}
\immediate\write\gtoutfile{Journal-ref: Algebr. Geom. Topol. \thevolumenumber\s
(\thevolumeyear) \startpage-\finishpage}
\immediate\write\gtoutfile{Comments: Published by Algebraic and
Geometric Topology at}
\immediate\write\gtoutfile{\s\s\s  http://www.maths.warwick.ac.uk/agt/AGTVol\thevolumenumber/agt-\thevolumenumber-\thepapernumber.abs.html}
\immediate\write\gtoutfile{\noexpand\\}
\immediate\write\gtoutfile{}
\ifx\theasciiabstract\relax
\immediate\write\gtoutfile{\theabstract}\else
\immediate\write\gtoutfile{\theasciiabstract}\fi
\immediate\write\gtoutfile{}
\immediate\write\gtoutfile{\noexpand\\}
\immediate\write\gtoutfile{}
\immediate\closeout\gtoutfile}}  

\def\maketitlepage{\makeagttitle\makeheadfile}
\let\makeshorttitle\maketitlepage
\let\maketitle\maketitlepage

\def\ifplaintex{\expandafter\ifx\csname documentclass\endcsname\relax}

\def\gtp{{\mathsurround=0pt\it $\cal G\mskip-2mu$eometry \&\ 
$\cal T\!\!$opology $\cal P\!$ublications}}  

\def\recd{{\small Received:\qua\receiveddate\ifx\reviseddate\relax
\else\qquad Revised:\qua\reviseddate\fi\par}} 

\def\lognumber#1{\def\thelognumber{#1}}
\def\volumenumber#1{\def\thevolumenumber{#1}}
\def\volumeyear#1{\def\thevolumeyear{#1}}
\def\papernumber#1{\def\thepapernumber{#1}}
\def\pagenumbers#1#2{\def\startpage{#1}\def\finishpage{#2}}
\def\published#1{\def\publishdate{#1}}

\def\received#1{\def\receiveddate{#1}}

\def\accepted#1{\def\accepteddate{#1}}

\long\def\asciiabstract#1{\long\def\theasciiabstract{#1}}

\let\\\par\let\thelognumber\relax\let\thevolumenumber\relax
\let\thepapernumber\relax\let\thevolumeyear\relax\let\startpage\relax
\let\finishpage\relax\let\publishdate\relax\let\receiveddate\relax
\let\reviseddate\relax\let\accepteddate\relax\let\theasciititle\relax
\let\theasciiauthors\relax
\let\theasciiabstract\relax

\let\theasciiemail\relax

\ifplaintex
\font\logobig=cmssbx10 scaled 3836
\font\logomed=cmssbx10 scaled 2557
\else
\font\logobig=cmssbx10 scaled 4200
\font\logomed=cmssbx10 scaled 2800
\fi

\long\def\makeagttitle{   
\count0=\startpage
\agt\hfill      
\hbox to 45truept{\vbox to 0pt{\vglue -13truept{\logomed A\kern -.37em{\logobig 
T}\kern -.38em G}\vss}\hss}
\break
{\small Volume \thevolumenumber\ (\thevolumeyear)
\startpage--\finishpage\nl
Published: \publishdate}

\vglue .25truein

{\parskip=0pt\leftskip 0pt plus
1fil\def\\{\par\smallskip}{\Large\bf\thetitle}\par\medskip} \vglue
0.05truein

{\parskip=0pt\leftskip 0pt plus 1fil\def\\{\par}{\sc\theauthors}
\par\medskip}%
 
\vglue 0.03truein 

{\small\leftskip 25truept\rightskip 25truept{\bf Abstract}\stdspace\theabstract

{\bf AMS Classification}\stdspace\theprimaryclass
\ifx\thesecondaryclass\relax\else; \thesecondaryclass\fi\par
{\bf Keywords}\stdspace \thekeywords\par}\vglue 7truept

}   

\ifplaintex
\font\phead=cmsl9 scaled 950
\font\pnum=cmbx10 scaled 913
\font\pfoot=cmsl9 scaled 950
\headline{\vbox to 0pt{\vskip -4.5mm\line{\small\phead\ifnum
\count0=\startpage ISSN 1472-2739 (on-line) 1472-2747 (printed)
\hfill {\pnum\folio}\else\ifodd\count0\def\\{ }%
\ifx\theshorttitle\relax\thetitle\else\theshorttitle\fi\hfill{\pnum\folio}
\else\def\\{ and }{\pnum\folio}\hfill\ifx\theshortauthors\relax\theauthors
\else\theshortauthors\fi\fi\fi}\vss}}
\footline{\vbox to 0pt{\vglue 0mm\line{\small\pfoot\ifnum\count0=\startpage
\copyright\ \gtp\hfill\else
\agt, Volume \thevolumenumber\ (\thevolumeyear)\hfill\fi}\vss}}
\else
\font=cmsl9 scaled 1050
\font\lnum=cmbx10 
\font=cmsl9 scaled 1050
\makeatletter
\def\@oddhead{{\small\lhead\ifnum\count0=\startpage ISSN 1472-2739 
(on-line) 1472-2747 (printed)\hfill {\lnum\number\count0}\else\ifodd\count0
\def\\{ }\ifx\theshorttitle\relax \thetitle \else\theshorttitle\fi\hfill
{\lnum\number\count0}\else\def\\{ and }{\lnum\number\count0}
\hfill\ifx\theshortauthors\relax 
\theauthors\else\theshortauthors\fi\fi\fi}}\def\@evenhead{\@oddhead}
\def\@oddfoot{\small\lfoot\ifnum\count0=\startpage\copyright\ \gtp\hfill\else
\agt, Volume \thevolumenumber\ (\thevolumeyear)\hfill\fi}
\def\@evenfoot{\@oddfoot}
\makeatother
\fi
\let\maketitlepage\makeagttitle
\let\makeshorttitle\maketitlepage
\let\maketitle\maketitlepage

\newwrite\gtoutfile
\long\gdef\makeheadfile{  
{\def\\{, }\def\s{ }
\immediate\openout\gtoutfile head.xxx
\immediate\write\gtoutfile{To: math@arxiv.org}
\immediate\write\gtoutfile{Subject: put OR rep NNNNN:ppppp}
\immediate\write\gtoutfile{--text follows this line--}
\immediate\write\gtoutfile{Proxy-for: \ifx\theasciiauthors\relax
\theauthors\else\theasciiauthors\fi\s<\ifx\theasciiemail\relax\theemail\else\theasciiemail\fi>}
\immediate\write\gtoutfile{\noexpand\\}
\immediate\write\gtoutfile{Authors: \ifx\theasciiauthors\relax
\theauthors\else\theasciiauthors\fi}
{\def\\{ }\immediate\write\gtoutfile{Title: \ifx\theasciititle\relax
\thetitle\else\theasciititle\fi}}
\immediate\write\gtoutfile{Subj-class: GT or SG, GR etc}
\immediate\write\gtoutfile{MSC-class: \theprimaryclass\ifx\thesecondaryclass\relax\else, \thesecondaryclass\fi}
\immediate\write\gtoutfile{Journal-ref: Algebr. Geom. Topol. \thevolumenumber\s
(\thevolumeyear) \startpage-\finishpage}
\immediate\write\gtoutfile{Comments: Published by Algebraic and
Geometric Topology at}
\immediate\write\gtoutfile{\s\s\s  http://www.maths.warwick.ac.uk/agt/AGTVol\thevolumenumber/agt-\thevolumenumber-\thepapernumber.abs.html}
\immediate\write\gtoutfile{\noexpand\\}
\immediate\write\gtoutfile{}
\ifx\theasciiabstract\relax
\immediate\write\gtoutfile{\theabstract}\else
\immediate\write\gtoutfile{\theasciiabstract}\fi
\immediate\write\gtoutfile{}
\immediate\write\gtoutfile{\noexpand\\}
\immediate\write\gtoutfile{}
\immediate\closeout\gtoutfile}}  

\def\maketitlepage{\makeagttitle\makeheadfile}
\let\makeshorttitle\maketitlepage
\let\maketitle\maketitlepage

\lognumber{16}
\volumenumber{2}
\volumeyear{2002}
\papernumber{16}
\published{21 May 2002}
\pagenumbers{337}{370}
\received{13 January 2002}
\accepted{27 February 2002}

\usepackage{amsmath,amssymb}
\usepackage{picins,slashbox,dbnsymb,
  graphicx,verbatim,longtable,
  epic,eepic,stmaryrd,rotating}
\usepackage[all]{xy}

\theoremstyle{plain}

\newtheorem{theorem}{Theorem}
\newtheorem{proposition}{Proposition}[section]

\newtheorem{lemma}[proposition]{Lemma}

\newtheorem{claim}[proposition]{Claim}

\newtheorem{conjecture}{Conjecture}

\theoremstyle{definition}
\newtheorem{definition}[proposition]{Definition}

\theoremstyle{remark}

\newtheorem{exercise}[proposition]{Exercise}

\newtheorem{remark}[proposition]{Remark}

\newlength{\standardunitlength}
\catcode`\@=11
\long\def\@makecaption#1#2{%
    \vskip 10pt
    \setbox\@tempboxa\hbox{
      \small{#1: }\ignorespaces #2}%
    \ifdim \wd\@tempboxa >\captionwidth {%
        \rightskip=\@captionmargin\leftskip=\@captionmargin
        \unhbox\@tempboxa\par}%
      \else
        \hbox to\hsize{\hfil\box\@tempboxa\hfil}%
    \fi}

\newdimen\@captionmargin\@captionmargin=2\parindent
\newdimen\captionwidth\captionwidth=\hsize
\catcode`\@=12

\newcommand{\Span}{\operatorname{span}}

\newlength{\globalparindent}
\def\la{\langle}
\def\ra{\rangle}

\def\bbF{{\mathbb F}}
\def\bbQ{{\mathbb Q}}

\def\bbZ{{\mathbb Z}}
\def\hatJ{{\hat J}}
\def\calC{{\mathcal C}}
\def\calF{{\mathcal F}}
\def\calH{{\mathcal H}}
\def\calO{{\mathcal O}}
\def\calX{{\mathcal X}}

\newcommand{\Kh}{{\text{\it Kh}}}
\newcommand{\qdim}{\operatorname{{\it q}dim}}

\def\KB#1#2{{\left\la
  \setlength{\unitlength}{#2\standardunitlength}
  \begin{array}{c}
    {\input figs/#1.tex } 
  \end{array}
\right\ra}}
\def\KhB#1#2{{
  \left\llbracket\hspace{-1mm}
  \setlength{\unitlength}{#2\standardunitlength}
  \begin{array}{c}
    {\input figs/#1.tex } 
  \end{array}
\hspace{-1mm}\right\rrbracket
}}

\begin{document}

\title{On Khovanov's categorification of the\\Jones polynomial}

\author{Dror Bar-Natan}
\address{Institute of Mathematics, The Hebrew University\\Giv'at-Ram, 
Jerusalem 91904, Israel}
\email{drorbn@math.huji.ac.il}
\url{http://www.ma.huji.ac.il/\char'176drorbn}


\begin{abstract}
  The working mathematician fears complicated words but loves pictures and
  diagrams. We thus give a no-fancy-anything picture rich glimpse into
  Khovanov's novel construction of ``the categorification of the Jones
  polynomial''. For the same low cost we also provide some
  computations, including one that shows that Khovanov's invariant is
  strictly stronger than the Jones polynomial and including a table of the
  values of Khovanov's invariant for all prime knots with up to $11$
  crossings.
\end{abstract}
\asciiabstract{
  The working mathematician fears complicated words but loves pictures and
  diagrams. We thus give a no-fancy-anything picture rich glimpse into
  Khovanov's novel construction of `the categorification of the Jones
  polynomial'. For the same low cost we also provide some
  computations, including one that shows that Khovanov's invariant is
  strictly stronger than the Jones polynomial and including a table of the
  values of Khovanov's invariant for all prime knots with up to 11
  crossings.}

\maketitle

\newcounter{In}

\def\In{\stepcounter{In}\vskip 6pt \par\noindent\footnotesize
  In[\theIn]$:=$\tt
  \catcode`\#=12 \catcode`\&=12 \catcode`\~=12 \catcode`\%=12
  \catcode`\$=12 \catcode`\^=12 \catcode`\_=12
  \obeylines\obeyspaces}
\def\Out{\vskip 6pt \par\noindent\footnotesize\parindent=12pt
  Out[\theIn]$=$\tt
  \catcode`\#=12 \catcode`\&=12 \catcode`\~=12 \catcode`\%=12
  \catcode`\$=12 \catcode`\^=12 \catcode`\_=12 \catcode`\{=12
  \catcode`\}=12 \obeyspaces
}
\def\Print{\vskip 6pt \par\noindent\footnotesize\tt
  \catcode`\#=12 \catcode`\&=12 \catcode`\~=12 \catcode`\%=12
  \catcode`\$=12 \catcode`\^=12 \catcode`\_=12 \catcode`\{=12
  \catcode`\}=12 \obeyspaces
}

\section{Introduction} In the summer of 2001 the author of this
note spent a week at Harvard University visiting David Kazhdan
and Dylan Thurston.  Our hope for the week was to understand and
improve Khovanov's seminal work on the categorification of the Jones
polynomial~\cite{Khovanov:Categorification, Khovanov:Functor}. We've
hardly achieved the first goal and certainly not the second; but we
did convince ourselves that there is something very new and novel in
Khovanov's work both on the deep conceptual level (not discussed here)
and on the shallower surface level. For on the surface level Khovanov
presents invariants of links which contain and generalize the Jones
polynomial but whose construction is like nothing ever seen in knot theory
before. Not being able to really digest it we decided to just chew some,
and then provide our output as a note containing a description of his
construction, complete and consistent and accompanied by computer code
and examples but stripped of all philosophy and of all the linguistic
gymnastics that is necessary for the philosophy but isn't necessary
for the mere purpose of having a working construction. Such a note may
be more accessible than the original papers. It may lead more people to
read Khovanov at the source, and maybe somebody reading such a note will
figure out what the Khovanov invariants really are. Congratulations! You
are reading this note right now.

\rk{1.1 Executive summary} \label{subsec:ExecutiveSummary} In
very brief words, Khovanov's idea is to replace the Kauffman bracket
$\la L\ra$ of a link projection $L$ by what we call ``the
Khovanov bracket'' $\llbracket L\rrbracket$, which is a chain complex
of graded vector spaces whose graded Euler characteristic is $\la
L\ra$. The Kauffman bracket is defined by the axioms
\[
  \la\emptyset\ra=1; \quad
  \la\bigcirc L\ra = (q+q^{-1})\la L\ra; \quad
  \la\backoverslash\ra = \la\hsmoothing\ra
    - q\la\smoothing\ra.
\]
Likewise, the definition of the Khovanov bracket can be summarized by the
axioms
\[
  \llbracket\emptyset\rrbracket=0\to\bbZ\to 0;
  \quad \llbracket\bigcirc L\rrbracket = V\otimes\llbracket L\rrbracket;
  \quad \llbracket\backoverslash\rrbracket =
  \calF\left(
    0\to\llbracket\hsmoothing\rrbracket
    \overset{d}{\to}\llbracket\smoothing\rrbracket\{1\}
    \to 0
  \right).
\]
Here $V$ is a vector space of graded dimension $q+q^{-1}$, the operator
$\{1\}$ is the ``degree shift by 1'' operation, which is the appropriate
replacement of ``multiplication by $q$'', $\calF$ is the ``flatten''
operation which takes a double complex to a single complex by taking
direct sums along diagonals, and a key ingredient, the differential $d$,
is yet to be defined.

The (unnormalized) Jones polynomial is a minor renormalization of the
Kauffman bracket, $\hatJ(L) = (-1)^{n_-}q^{n_+-2n_-}\la L\ra$.  The
Khovanov invariant $\calH(L)$ is the homology of a similar renormalization
$\llbracket L\rrbracket[-n_-]\{n_+-2n_-\}$ of the Khovanov bracket. The
``main theorem'' states that the Khovanov invariant is indeed a link
invariant and that its graded Euler characteristic is $\hatJ(L)$.
Anything in $\calH(L)$ beyond its Euler characteristic appears to be new,
and direct computations show that there really is more in $\calH(L)$
than in its Euler characteristic.

{\bf 1.2 Acknowledgements}\qua I wish to thank David Kazhdan and Dylan
Thurston for the week at Harvard that led to writing of this note and for
their help since then. I also wish to thank G.~Bergman, S.~Garoufalidis,
J.~Hoste, V.~Jones, M.~Khovanov, A.~Kricker, G.~Kuperberg, A.~Stoimenow
and M.~Thistlethwaite for further assistance, comments and suggestions.

\vspace{-4mm}

\section{The Jones polynomial}

\vspace{-4mm}

\parpic[r]{\begin{tabular}{c}
  \includegraphics[width=0.8in]{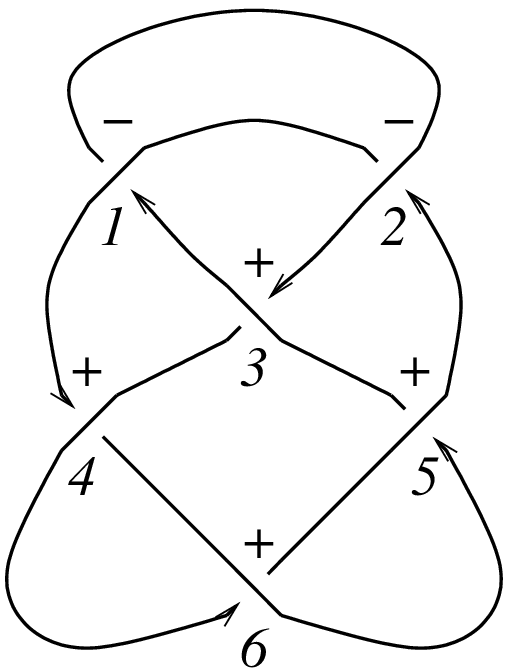} \\
  \small $n_+=4$; $n_-=2$
\end{tabular}}
All of our links are oriented links in an oriented Euclidean space. We
will present links using their projections to the plane as shown in the
example on the right. Let $L$ be a link projection, let $\calX$ be the
set of crossings of $L$, let $n=|\calX|$, let us number the elements of
$\calX$ from $1$ to $n$ in some arbitrary way and let us write
$n=n_++n_-$ where $n_+$ ($n_-$) is the number of right-handed
(left-handed) crossings in $\calX$. (again, look to the right).

Recall that the Kauffman bracket~\cite{Kauffman:OnKnots} of $L$ is
defined by the formulas\footnote{Our slightly unorthodox conventions
follow~\cite{Khovanov:Categorification}. At some minor regrading and
renaming cost, we could have used more standard conventions just as well.}
$\la\emptyset\ra=1$, $\la\bigcirc L\ra = (q+q^{-1})\la
L\ra$ and $\la\backoverslash\ra = \la\hsmoothing\ra
- q\la\smoothing\ra$, that the unnormalized Jones polynomial is
defined by $\hatJ(L) = (-1)^{n_-}q^{n_+-2n_-}\la L\ra$, and that
the Jones polynomial of $L$ is simply $J(L):=\hatJ(L)/(q+q^{-1})$.
We name $\hsmoothing$ and $\smoothing$ the 0- and 1-smoothing of
$\backoverslash$, respectively. With this naming convention each vertex
$\alpha\in\{0,1\}^\calX$ of the $n$-dimensional cube $\{0,1\}^\calX$
corresponds in a natural way to a ``complete smoothing'' $S_\alpha$ of $L$
where all the crossings are smoothed and the result is just a union of
planar cycles. To compute the unnormalized Jones polynomial, we replace
each such union $S_\alpha$ of (say) $k$ cycles with a term of the form
$(-1)^rq^r(q+q^{-1})^k$, where $r$ is the ``height'' of a smoothing, the
number of 1-smoothings used in it.  We then sum all these terms over all
$\alpha\in\{0,1\}^\calX$ and multiply by the final normalization term,
$(-1)^{n_-}q^{n_+-2n_-}$.  Thus the whole procedure (in the case of the
trefoil knot) can be depicted as in the diagram below. Notice that in this
diagram we have split the summation over the vertices of $\{0,1\}^\calX$
to a summation over vertices of a given height followed by a summation
over the possible heights. This allows us to factor out the $(-1)^r$
factor and turn the final summation into an alternating summation:


{\centering\resizebox*{\textwidth}{!}{\parbox{6.2in}{\begin{gather}
  \label{eq:JCube}
  \def\S#1#2{\fbox{
    \!\!\parbox[c]{10mm}{\includegraphics[width=10mm]{figs/#1.eps}}
    \hspace{-3mm}\raisebox{3.5mm}{${\scriptstyle #2}$}
    \hspace{-5mm}\raisebox{-3.5mm}{$\scriptscriptstyle #1$} \hspace{-3mm}
  }}
  \begin{array}{c} \xymatrix@R=1.414cm@C=0.866cm{
    \includegraphics[width=0.7in]{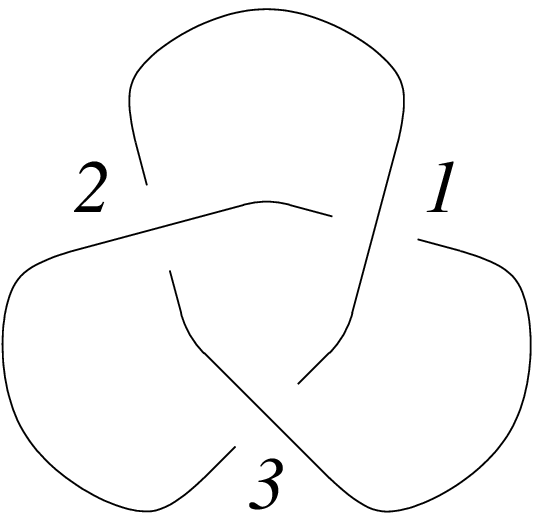}
      & \S{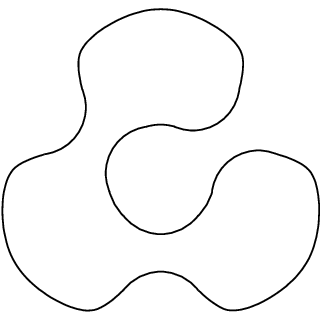}{q(q+q^{-1})} \ar@{-}[r] \ar@{-}[rd]|\hole \ar@{.}[d]_+
      & \S{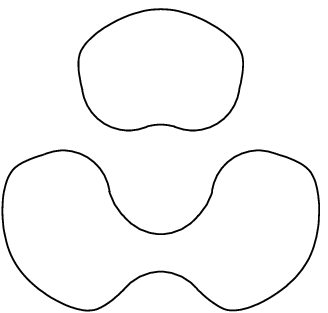}{q^2(q+q^{-1})^2} \ar@{-}[rd] \ar@{.}[d]^+
      & \\
    \S{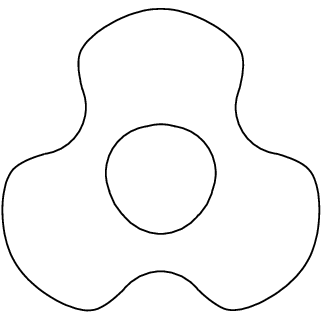}{(q+q^{-1})^2} \ar@{-}[ru] \ar@{-}[r] \ar@{-}[rd] \ar@{.>}[dd]
      & \S{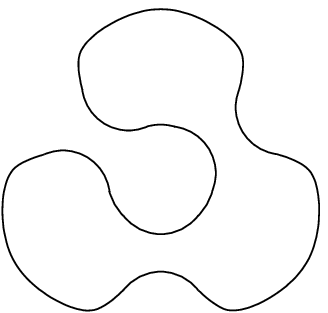}{q(q+q^{-1})} \ar@{-}[ru] \ar@{-}[rd] \ar@{.}[d]_+
      & \S{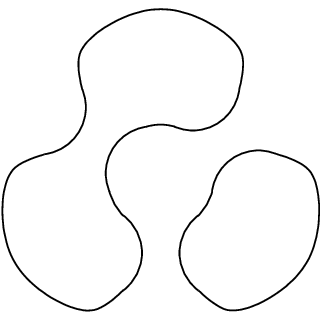}{q^2(q+q^{-1})^2} \ar@{-}[r] \ar@{.}[d]^+
      & \S{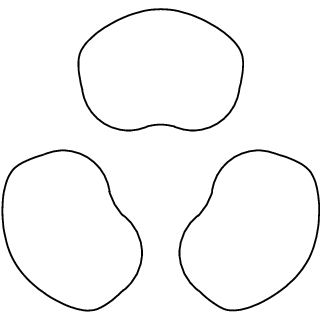}{q^3(q+q^{-1})^3} \ar@{.>}[dd] \\
    & \S{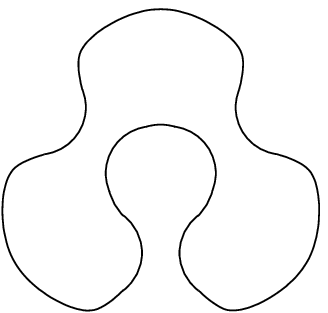}{q(q+q^{-1})} \ar@{-}[ru]|\hole \ar@{-}[r] \ar@{.>}[d]
      & \S{011}{q^2(q+q^{-1})^2} \ar@{-}[ru] \ar@{.>}[d] & \\
    (q+q^{-1})^2 \ar@{}[r]|{-}
      & 3q(q+q^{-1}) \ar@{}[r]|{+}
      & 3q^2(q+q^{-1})^2 \ar@{}[r]|{-}
      & q^3(q+q^{-1})^3 
  } \end{array}
\\
  = \ q^{-2} + 1 + q^2 - q^6 \ 
  \xrightarrow[\text{(with $(n_+,n_-)=(3,0)$)}]{\cdot(-1)^{n_-}q^{n_+-2n_-}}
  \ q + q^3 + q^5 - q^9
  \ \xrightarrow{\cdot(q+q^{-1})^{-1}} \ 
  J(\righttrefoil) = q^2+q^6-q^8. \notag
\end{gather}}}}

\section{Categorification}

\subsection{Spaces}

Khovanov's ``categorification'' idea is to replace polynomials
by graded vector spaces\footnote{Everything that we do works just
as well (with some linguistic differences) over $\bbZ$. In fact,
in~\cite{Khovanov:Categorification} Khovanov works over the even more
general ground ring $\bbZ[c]$ where $\deg c=2$.} of the appropriate
``graded dimension'', so as to turn the Jones polynomial into a
homological object. With the diagram~\eqref{eq:JCube} as the starting
point the process is straight forward and essentially unique. Let us
start with a brief on some necessary generalities:

\begin{definition} Let $W=\bigoplus_m W_m$ be a graded vector space
with homogeneous components $\{W_m\}$.  The graded dimension of $W$ is
the power series $\qdim W:=\sum_mq^m\dim W_m$.
\end{definition}

\begin{definition} Let $\cdot\{l\}$ be the ``degree shift'' operation
on graded vector spaces. That is, if $W=\bigoplus_m W_m$ is a graded
vector space, we set $W\{l\}_m:=W_{m-l}$, so that $\qdim W\{l\} =
q^l\qdim W$.
\end{definition}

\begin{definition} Likewise, let $\cdot[s]$ be the ``height shift''
operation on chain complexes. That is, if $\bar\calC$ is a chain
complex $\ldots\to\bar\calC^r\overset{d^r}{\to}\bar\calC^{r+1}\ldots$
of (possibly graded) vector spaces (we call $r$ the ``height'' of a
piece $\bar\calC^r$ of that complex), and if $\calC=\bar\calC[s]$, then
$\calC^r=\bar\calC^{r-s}$ (with all differentials shifted
accordingly).
\end{definition}

Armed with these three notions, we can proceed with ease. Let $L$,
$\calX$, $n$ and $n_\pm$ be as in the previous section. Let $V$ be the
graded vector space with two basis elements $v_\pm$ whose degrees are
$\pm 1$ respectively, so that $\qdim V=q+q^{-1}$. With every vertex
$\alpha\in\{0,1\}^\calX$ of the cube $\{0,1\}^\calX$ we associate
the graded vector space $V_\alpha(L):=V^{\otimes k}\{r\}$, where $k$
is the number of cycles in the smoothing of $L$ corresponding to
$\alpha$ and $r$ is the height $|\alpha|=\sum_i\alpha_i$ of $\alpha$
(so that $\qdim V_\alpha(L)$ is the polynomial that appears at the
vertex $\alpha$ in the cube at~\eqref{eq:JCube}). We then set the
$r$th chain group $\llbracket L\rrbracket^r$ (for $0\leq r\leq n$) to
be the direct sum of all the vector spaces at height $r$: $\llbracket
L\rrbracket^r:=\bigoplus_{\alpha:r=|\alpha|}V_\alpha(L)$. Finally
(for this long paragraph), we gracefully ignore the fact that
$\llbracket L\rrbracket$ is not yet a complex, for we have not yet
endowed it with a differential, and we set $\calC(L):=\llbracket
L\rrbracket[-n_-]\{n_+-2n_-\}$.  Thus the diagram~\eqref{eq:JCube}
(in the case of the trefoil knot) becomes:

\newpage

\begin{multline} \label{eq:VCube}
  \def\S#1#2{\fbox{
    \!\!\parbox[c]{10mm}{\includegraphics[width=10mm]{figs/#1.eps}}
    \hspace{-3mm}\raisebox{3.5mm}{${\scriptstyle #2}$}
    \hspace{-5mm}\raisebox{-3.5mm}{$\scriptscriptstyle #1$} \hspace{-3mm}
  }}
  \begin{array}{c} \xymatrix@R=1.414cm@C=0.866cm{
    \includegraphics[width=0.7in]{figs/Trefoil.eps}
      & \S{100}{V\{1\}} \ar@{-}[r] \ar@{-}[rd]|\hole \ar@{.}[d]_\oplus
      & \S{110}{V^{\otimes 2}\{2\}} \ar@{-}[rd] \ar@{.}[d]^\oplus
      & \\
    \S{000}{V^{\otimes 2}} \ar@{-}[ru] \ar@{-}[r] \ar@{-}[rd] \ar@{.>}[dd]
      & \S{010}{V\{1\}} \ar@{-}[ru] \ar@{-}[rd] \ar@{.}[d]_\oplus
      & \S{101}{V^{\otimes 2}\{2\}} \ar@{-}[r] \ar@{.}[d]^\oplus
      & \S{111}{V^{\otimes 3}\{3\}} \ar@{.>}[dd] \\
    & \S{001}{V\{1\}} \ar@{-}[ru]|\hole \ar@{-}[r] \ar@{.>}[d]
      & \S{011}{V^{\otimes 2}\{2\}} \ar@{-}[ru] \ar@{.>}[d] & \\
    \Big(\llbracket\righttrefoil\rrbracket^0 \ar@{}|;[r]
      & \llbracket\righttrefoil\rrbracket^1 \ar@{}|;[r]
      & \llbracket\righttrefoil\rrbracket^2 \ar@{}|;[r]
      & \llbracket\righttrefoil\rrbracket^3 \Big)
  } \end{array}
\\
  \qquad\qquad = \quad \llbracket\righttrefoil\rrbracket \quad
  \xrightarrow[\text{(with $(n_+,n_-)=(3,0)$)}]{\ \cdot[-n_-]\{n_+-2n_-\}\ }
  \quad \calC(\righttrefoil).
\end{multline}

The graded Euler characteristic $\chi_q(\calC)$ of a chain complex
$\calC$ is defined to be the alternating sum of the graded dimensions
of its homology groups, and, if the degree of the differential $d$ is
$0$ and all chain groups are finite dimensional, it is also equal to
the alternating sum of the graded dimensions of the chain groups. A few
paragraphs down we will endow $\calC(L)$ with a degree $0$ differential.
This granted and given that the chains of $\calC(L)$ are already
defined, we can state and prove the following theorem:

\begin{theorem} \label{thm:Euler} The graded Euler characteristic of
$\calC(L)$ is the unnormalized Jones polynomial of $L$:
$$\chi_q(\calC(L)) = \hatJ(L). $$
\end{theorem}

\begin{proof} The theorem is trivial by design; just compare
diagrams~\eqref{eq:JCube} and~\eqref{eq:VCube} and all the relevant
definitions. Thus rather than a proof we comment on the statement and the
construction preceding it: If one wishes our theorem to hold,
everything in the construction of diagram~\eqref{eq:VCube} is forced,
except the height shift $[-n_-]$. The parity of this shift is determined by
the $(-1)^{n_-}$ factor in the definition of $\hatJ(L)$. The given choice of
magnitude is dictated within the proof of Theorem~\ref{thm:Main}.
\end{proof}

\subsection{Maps} \label{subsec:Maps}

Next, we wish to turn the sequence of spaces $\calC(L)$ into a chain
complex. Let us flash the answer upfront, and only then go through the
traditional ceremony of formal declarations:

{\centering\resizebox*{\textwidth}{!}{\parbox{5.85in}{\begin{gather}
  \def\S#1#2{\fbox{
    \!\!\parbox[c]{10mm}{\includegraphics[width=10mm]{figs/#1.eps}}
    \hspace{-3mm}\raisebox{3.5mm}{${\scriptstyle #2}\,$}
    \hspace{-6mm}\raisebox{-3.5mm}{$\scriptscriptstyle #1$} \hspace{-3mm}
  }}
  \def\C#1{\hspace{-2mm}\includegraphics[width=6mm]{figs/C#1.eps}}
  \def\neg{\hspace{-2mm}}
  \begin{array}{c} \xymatrix@R=1.697cm@C=2cm{
    \includegraphics[width=0.7in]{figs/Trefoil.eps}
      & \S{100}{V\{1\}}
        \ar@{o->}[r]^{\C{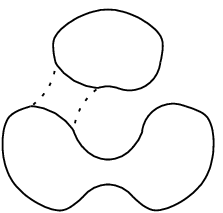}}_{d_{1\star 0}}
        \ar@{o->}[rd]|\hole^(0.36){\C{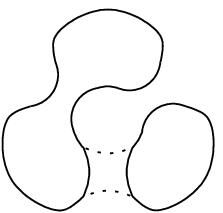}}_(0.6){d_{10\star}\neg}
        \ar@{.}[d]_\oplus
      & \S{110}{V^{\otimes 2}\{2\}}
        \ar[rd]^{\C{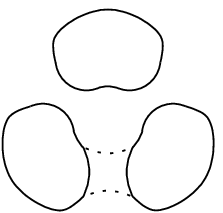}}_{d_{11\star}\neg}
        \ar@{.}[d]^\oplus
      & \\
    \S{000}{V^{\otimes 2}}
        \ar[ru]^{\C{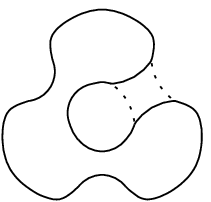}\neg}_{d_{\star 00}}
        \ar[r]^{\C{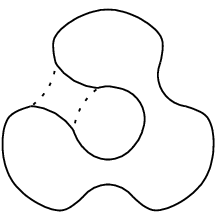}}_{d_{0\star 0}}
        \ar[rd]^{\C{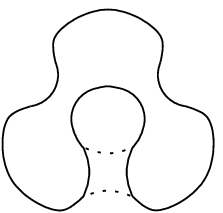}}_{d_{00\star}}
        \ar@{.>}[dd]
      & \S{010}{V\{1\}}
        \ar[ru]^(0.3){\C{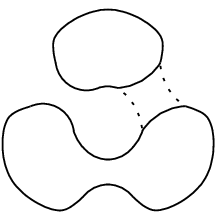}\neg}_(0.66){d_{\star 10}}
        \ar@{o->}[rd]^(0.36){\C{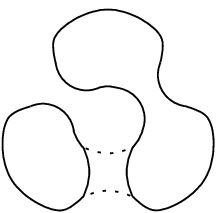}}_(0.6){d_{01\star}\neg}
        \ar@{.}[d]_\oplus
      & \S{101}{V^{\otimes 2}\{2\}}
        \ar@{o->}[r]^{\C{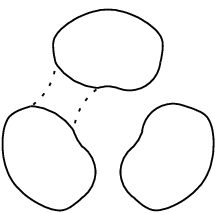}}_{d_{1\star 1}}
        \ar@{.}[d]^\oplus
      & \S{111}{V^{\otimes 3}\{3\}}
        \ar@{.>}[dd] \\
    \ar@{}[r]_\ ="1"
      & \S{001}{V\{1\}}
        \ar[ru]|\hole^(0.3){\C{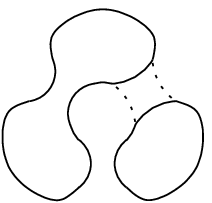}\neg}_(0.66){d_{\star 01}}
        \ar[r]^{\C{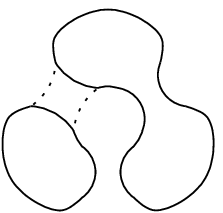}}_{d_{0\star 1}}="3"
        \ar@{.>}[d]
      & \S{011}{V^{\otimes 2}\{2\}}
        \ar[ru]^{\C{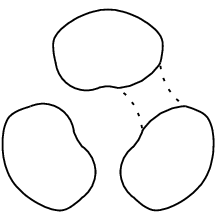}\neg}_{d_{\star 11}}
        \ar@{.>}[d]
        \ar@{}[r]_\ ="5" & \\
    \llbracket\righttrefoil\rrbracket^0
        \ar[r]^{d^0}="2"
      & \llbracket\righttrefoil\rrbracket^1
        \ar[r]^{d^1}="4"
      & \llbracket\righttrefoil\rrbracket^2
        \ar[r]^{d^2}="6"
      & \llbracket\righttrefoil\rrbracket^3
    \ar@{.>}"1";"2"|{\underset{|\xi|=0}{\sum}(-1)^\xi d_\xi}
    \ar@{.>}"3";"4"|{\underset{|\xi|=1}{\sum}(-1)^\xi d_\xi}
    \ar@{.>}"5";"6"|{\underset{|\xi|=2}{\sum}(-1)^\xi d_\xi}
  } \end{array}
\notag \\ \label{eq:KCube}
  \qquad\qquad = \quad \llbracket\righttrefoil\rrbracket \quad
  \xrightarrow[\text{(with $(n_+,n_-)=(3,0)$)}]{\ \cdot[-n_-]\{n_+-2n_-\}\ }
  \quad \calC(\righttrefoil).
\end{gather}}}}

This diagram certainly looks threatening, but in fact, it's quite
harmless.  Just hold on tight for about a page!  The chain groups
$\llbracket L\rrbracket^r$ are, as we have already seen, direct sums
of the vector spaces that appear in the vertices of the cube along the
columns above each one of the $\llbracket L\rrbracket^r$ spaces. We do
the same for the arrows $d^r$ --- we turn each edge $\xi$ of the cube to
map between the vector spaces at its ends, and then we add up these maps
along columns as shown above.  The edges of the cube $\{0,1\}^\calX$ can
be labeled by sequences in $\{0,1,\star\}^\calX$ with just one $\star$
(so the tail of such an edge is found by setting $\star\to 0$ and the
head by setting $\star\to 1$). The height $|\xi|$ of an edge $\xi$
is defined to be the height of its tail, and hence if the maps on the
edges are called $d_\xi$ (as in the diagram), then the vertical collapse
of the cube to a complex becomes $d^r:=\sum_{|\xi|=r}(-1)^\xi d_\xi$.

It remains to explain the signs $(-1)^\xi$ and to define the per-edge
maps $d_\xi$. The former is easy. To get the differential $d$ to
satisfy $d\circ d=0$, it is enough that all square faces of the cube
would anti-commute. But it is easier to arrange the $d_\xi$'s so that
these faces would (positively) commute; so we do that and then sprinkle
signs to make the faces anti-commutative. One may verify that this can
be done by multiplying $d_\xi$ by $(-1)^\xi:=(-1)^{\sum_{i<j}\xi_i}$,
where $j$ is the location of the $\star$ in $\xi$. In
diagram~\eqref{eq:KCube} we've indicated the edges $\xi$ for which
$(-1)^\xi=-1$ with little circles at their tails. The reader is welcome to
verify that there is an odd number of such circles around each face of the
cube shown.

It remains to find maps $d_\xi$ that make the cube commutative (when taken
with no signs) and that
are of degree $0$ so as not to undermine Theorem~\ref{thm:Euler}. The
space $V_\alpha$ on each vertex $\alpha$ has as many tensor factors as
there are cycles in the smoothing $S_\alpha$. Thus we put these tensor
factors in $V_\alpha$ and cycles in $S_\alpha$ in bijective
correspondence once and for all. Now for any edge $\xi$, the smoothing
at the tail of $\xi$ differs from the smoothing at the head of $\xi$ by
just a little: either two of the cycles merge into one (see say
$\xi=0\!\star\!0$ above) or one of the cycles splits in two (see say
$\xi=1\!\star\!1$ above). So for any $\xi$, we set $d_\xi$ to be the
identity on the tensor factors corresponding to the cycles that don't
participate, and then we complete the definition of $\xi$ using two
linear maps $m:V\otimes V\to V$ and $\Delta:V\to V\otimes V$ as
follows:
{\def\neg{\hspace{-1mm}}
  \begin{align} \label{eq:mdef}
    \big(\neg\begin{array}{c}
      \includegraphics[width=1in]{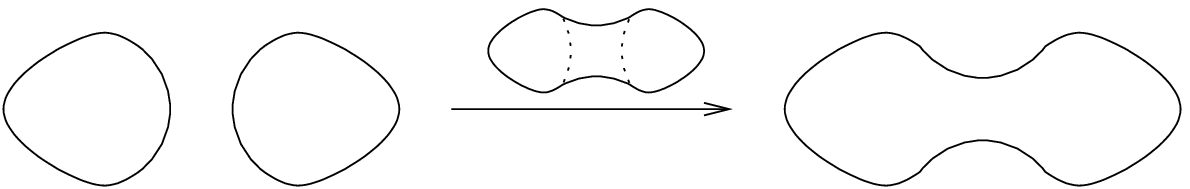}
    \end{array}\neg\big)
    & \longrightarrow \big(V\otimes V\overset{m}{\rightarrow}V\big)
    &\quad& m:\begin{cases}
      v_+\otimes v_-\mapsto v_- &
      v_+\otimes v_+\mapsto v_+ \\
      v_-\otimes v_+\mapsto v_- &
      v_-\otimes v_-\mapsto 0
    \end{cases}
  \\ \label{eq:DeltaDef}
    \big(\neg\begin{array}{c}
      \includegraphics[width=1in]{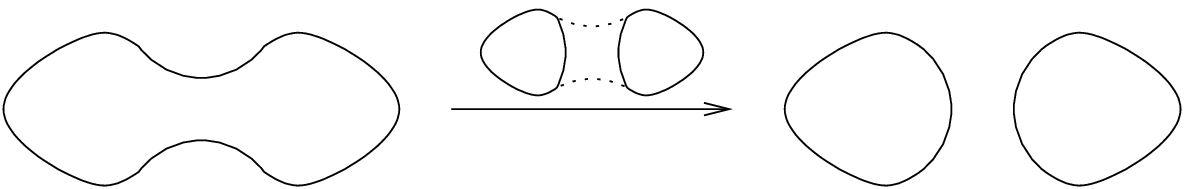}
    \end{array}\neg\big)
    & \longrightarrow \big(V\overset{\Delta}{\rightarrow}V\otimes V\big)
    &\quad& \Delta:\begin{cases}
      v_+ \mapsto v_+\otimes v_- + v_-\otimes v_+ &\\
      v_- \mapsto v_-\otimes v_- &
    \end{cases}
  \end{align}
}

We note that because of the degree shifts in the definition of the
$V_\alpha$'s and because we want the $d_\xi$'s to be of degree $0$, the
maps $m$ and $\Delta$ must be of degree $-1$. Also, as there is no
canonical order on the cycles in $S_\alpha$ (and hence on the tensor
factors of $V_\alpha$), $m$ and $\Delta$ must be commutative and
co-commutative respectively. These requirements force the equality
$m(v_+\otimes v_-)=m(v_-\otimes v_+)$ and force the values of $m$ and
$\Delta$ to be as shown above up to scalars.

\begin{remark} It is worthwhile to note, though not strictly necessary to
the understanding of this note, that the cube in diagram~\eqref{eq:KCube}
is related to a certain $(1+1)$-dimensional topological quantum field
theory (TQFT). Indeed, given any $(1+1)$-dimensional TQFT one may assign
vector spaces to the vertices of $\{0,1\}^\calX$ and maps to the edges ---
on each vertex we have a union of cycles which is a $1$-manifold that gets
mapped to a vector space via the TQFT, and on each edge we can place the
obvious $2$-dimensional saddle-like cobordism between the $1$-manifolds
on its ends, and then get a map between vector spaces using the TQFT. The
cube in diagram~\eqref{eq:KCube} comes from this construction if one
starts from the TQFT corresponding to the Frobenius algebra defined by
$V$, $m$, $\Delta$, the unit $v_+$ and the co-unit $\epsilon\in V^\star$
defined by $\epsilon(v_+)=0$, $\epsilon(v_-)=1$. See more
in~\cite{Khovanov:Categorification}.
\end{remark}

\begin{exercise} Verify that the definitions given in this section agree
with the ``executive summary'' (Section~\ref{subsec:ExecutiveSummary}).
\end{exercise}

\subsection{A notational digression} \label{subsec:Notation}

For notational and computational reasons\footnote{You may skip this
section if the previous section was clear enough and you don't intend
to read the computational Section~\ref{sec:ComputerTalk}.} it is
convenient to also label the edges of $L$. Our convention is to reserve
separate interval of integers for each component, and then to label the
edges within this component in an ascending order (except for one jump
down) --- see Figure~\ref{fig:LinkNotation} in
Section~\ref{sec:ComputerTalk}. Given $\alpha\in\{0,1\}^\calX$, we
label every cycle in the smoothing $S_\alpha$ by the minimal edge that
appears in it, and then we label the tensor factor in $V_\alpha$
accordingly. So for example (with $L=\righttrefoil$ labeled as in
Figure~\ref{fig:LinkNotation}), the big and small components of
$S_{011}=\parbox[c]{10mm}{\includegraphics[width=10mm]{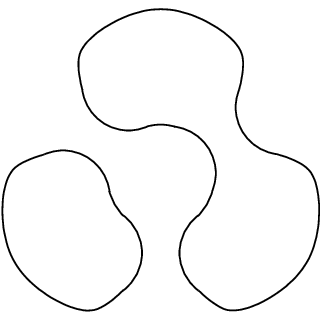}}$
would be labeled $1$ and $3$ respectively, and thus $V_{011}$ would be
$V_1\otimes V_3\{2\}$. The indices in the latter space have only a
notational meaning that allows us easier access to its tensor factors.
Thus $V_1\otimes V_3\cong V\otimes V$, yet the standard basis elements
of $V_1\otimes V_3$ have nice standard names:  $\{v_+^1v_+^3,
v_+^1v_-^3, v_-^1v_+^3, v_-^1v_-^3\}$.

With this notation, we can make the cube of Equation~\eqref{eq:KCube} a
little more explicit. We denote by $\Delta^{ij}$ the map which acts on
a tensor product of labeled copies of $V$ as the identity on all
factors except the one labeled $V_{\min(i,j)}$ which gets mapped by
$\Delta$ of Equation~\eqref{eq:DeltaDef} to $V_i\otimes V_j$. Likewise
$m_{ij}$ denotes the natural extension by identity maps of
$m:V_i\otimes V_j\to V_{\min(i,j)}$. All said, the cube in
diagram~\eqref{eq:KCube} becomes:

{\centering\resizebox*{\textwidth}{!}{\parbox{6.2in}{\begin{equation}
  \label{eq:LCube}
  \def\S#1#2{\fbox{
    \!\!\parbox[c]{10mm}{\includegraphics[width=10mm]{figs/#1.eps}}
    \hspace{-3mm}\raisebox{3.5mm}{${\scriptstyle #2}\,$}
    \hspace{-6mm}\raisebox{-3.5mm}{$\scriptscriptstyle #1$} \hspace{-3mm}
  }}
  \def\neg{\hspace{-2mm}}
  \begin{array}{c} \xymatrix@R=1.697cm@C=1.6cm{
    \includegraphics[width=0.7in]{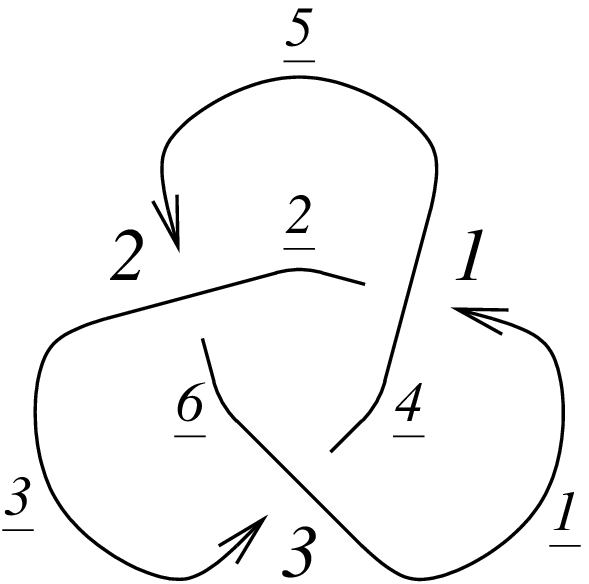}
      & \S{100}{V_1\{1\}}
        \ar@{o->}[r]^{\Delta^{12}}_{d_{1\star 0}}
        \ar@{o->}[rd]|\hole^(0.36){\Delta^{12}}_(0.6){d_{10\star}\neg}
      & \S{110}{V_1\otimes V_2\{2\}}
        \ar[rd]^{\Delta^{13}}_{d_{11\star}\neg}
      & \\
    \S{000}{V_1\otimes V_2}
        \ar[ru]^{m_{12}}_{d_{\star 00}}
        \ar[r]^{m_{12}}_{d_{0\star 0}}
        \ar[rd]^{m_{12}}_{d_{00\star}}
      & \S{010}{V_1\{1\}}
        \ar[ru]^(0.3){\Delta^{12}\neg}_(0.66){d_{\star 10}}
        \ar@{o->}[rd]^(0.36){\Delta^{13}}_(0.6){d_{01\star}\neg}
      & \S{101}{V_1\otimes V_2\{2\}}
        \ar@{o->}[r]^{\Delta^{23}}_{d_{1\star 1}}
      & \S{111}{V_1\otimes V_2\otimes V_3\{3\}} \\
    \ar@{}[r]
      & \S{001}{V_1\{1\}}
        \ar[ru]|\hole^(0.3){\Delta^{12}\neg}_(0.66){d_{\star 01}}
        \ar[r]^{\Delta^{13}}_{d_{0\star 1}}
      & \S{011}{V_1\otimes V_3\{2\}}
        \ar[ru]^{\Delta^{12}\neg}_{d_{\star 11}}
        \ar@{}[r] &
  } \end{array}
\end{equation}}}}

\subsection{The main theorem}

\begin{claim} The $n$-dimensional cube as in Equation~\eqref{eq:KCube}
(just as well,~\eqref{eq:LCube}) is commutative (for any $L$, and provided
all maps are taken with no signs) and hence
the sequences $\llbracket L\rrbracket$ and $\calC(L)$ are chain complexes.
\end{claim}

\begin{proof} A routine verification.
\end{proof}

Let $\calH^r(L)$ denote the $r$th cohomology of the complex $\calC(L)$. It
is a graded vector space depending on the link projection $L$. Let $\Kh(L)$
denote the graded Poincar\'e polynomial of the complex $\calC(L)$ in the
variable $t$; i.e., let
\[ \Kh(L):=\sum_r t^r\qdim\calH^r(L). \]
(When we wish to emphasize the ground field $\bbF$, we write
$\Kh_\bbF(L)$.)

\begin{theorem}[Khovanov~\cite{Khovanov:Categorification}]  
\label{thm:Main}
The graded dimensions of the
homology groups $\calH^r(L)$ are link invariants, and hence $\Kh(L)$,
a polynomial in the variables $t$ and $q$, is a link invariant that
specializes to the unnormalized Jones polynomial at $t=-1$.
\end{theorem}

\subsection{Proof of the main theorem}

To prove Theorem~\ref{thm:Main}, we need to study the behavior of
$\llbracket L\rrbracket$ under the three Reidemeister moves\footnote{We
leave it to the reader to confirm that no further variants of
these moves need to be considered. For example, we check only the
``right twist'' version of (R1). The left twist
version follows from it and from (R2).} (R1):
$
  \setlength{\unitlength}{0.2\standardunitlength}
  \begin{array}{c}
    {\begingroup\makeatletter\ifx\SetFigFont\undefined%
\gdef\SetFigFont#1#2#3#4#5{%
  \reset@font\fontsize{#1}{#2pt}%
  \fontfamily{#3}\fontseries{#4}\fontshape{#5}%
  \selectfont}%
\fi\endgroup%
{\renewcommand{\dashlinestretch}{30}
\begin{picture}(1244,859)(0,-10)
\thicklines
\path(1222,22)(1219,22)(1212,23)
	(1200,24)(1183,25)(1161,28)
	(1137,31)(1110,35)(1083,39)
	(1056,45)(1030,51)(1006,57)
	(982,65)(960,74)(938,84)
	(916,95)(895,108)(872,122)
	(855,134)(837,147)(818,161)
	(799,176)(779,192)(759,209)
	(739,227)(718,246)(697,266)
	(677,286)(656,307)(636,329)
	(616,351)(597,373)(579,395)
	(562,416)(546,438)(531,459)
	(518,479)(506,499)(495,518)
	(486,537)(478,555)(472,572)
	(467,589)(463,606)(461,622)
	(460,638)(460,654)(461,669)
	(464,684)(467,699)(472,713)
	(478,726)(485,739)(493,751)
	(501,763)(511,773)(521,782)
	(531,791)(542,798)(553,805)
	(565,810)(576,814)(588,818)
	(599,820)(611,822)(622,822)
	(637,821)(653,819)(668,815)
	(684,809)(700,802)(716,794)
	(731,784)(747,773)(762,762)
	(776,749)(789,736)(801,723)
	(812,710)(821,696)(830,683)
	(837,671)(843,659)(847,647)
	(851,630)(852,614)(850,598)
	(845,580)(836,562)(824,542)
	(811,521)(797,502)(784,487)
	(776,477)(773,473)(772,472)
\path(472,172)(469,169)(462,163)
	(451,153)(437,140)(420,126)
	(402,111)(385,97)(369,85)
	(354,75)(339,66)(325,58)
	(311,52)(297,47)(284,43)
	(270,39)(255,36)(238,34)
	(219,32)(198,30)(174,28)
	(148,26)(122,25)(95,24)
	(70,23)(50,23)(35,22)
	(26,22)(23,22)(22,22)
\end{picture}
} } 
  \end{array}
\leftrightarrow
  \setlength{\unitlength}{0.2\standardunitlength}
  \begin{array}{c}
    {\begingroup\makeatletter\ifx\SetFigFont\undefined%
\gdef\SetFigFont#1#2#3#4#5{%
  \reset@font\fontsize{#1}{#2pt}%
  \fontfamily{#3}\fontseries{#4}\fontshape{#5}%
  \selectfont}%
\fi\endgroup%
{\renewcommand{\dashlinestretch}{30}
\begin{picture}(1244,259)(0,-10)
\thicklines
\path(22,22)(23,22)(26,22)
	(35,23)(51,23)(72,24)
	(97,26)(125,28)(153,31)
	(180,34)(205,37)(228,41)
	(249,46)(269,51)(287,57)
	(305,64)(322,72)(337,80)
	(353,89)(369,99)(385,109)
	(402,120)(419,131)(436,142)
	(454,154)(472,165)(490,176)
	(508,186)(525,195)(542,203)
	(559,210)(575,215)(591,219)
	(607,221)(622,222)(637,221)
	(653,219)(669,215)(685,210)
	(702,203)(719,195)(736,186)
	(754,176)(772,165)(790,154)
	(808,142)(825,131)(842,120)
	(859,109)(875,99)(891,89)
	(907,80)(922,72)(939,64)
	(957,57)(975,51)(995,46)
	(1016,41)(1039,37)(1064,34)
	(1091,31)(1119,28)(1147,26)
	(1172,24)(1193,23)(1209,23)
	(1218,22)(1221,22)(1222,22)
\end{picture}
} } 
  \end{array}
$, (R2):
$
  \setlength{\unitlength}{0.2\standardunitlength}
  \begin{array}{c}
    {\begingroup\makeatletter\ifx\SetFigFont\undefined%
\gdef\SetFigFont#1#2#3#4#5{%
  \reset@font\fontsize{#1}{#2pt}%
  \fontfamily{#3}\fontseries{#4}\fontshape{#5}%
  \selectfont}%
\fi\endgroup%
{\renewcommand{\dashlinestretch}{30}
\begin{picture}(2144,659)(0,-10)
\thicklines
\path(22,22)(23,22)(26,22)
	(35,23)(51,23)(72,24)
	(97,26)(125,28)(153,31)
	(180,34)(205,37)(228,41)
	(249,46)(269,51)(287,57)
	(305,64)(322,72)(337,80)
	(353,89)(369,99)(385,110)
	(402,122)(419,135)(436,149)
	(454,164)(472,179)(490,195)
	(508,211)(525,227)(542,243)
	(559,260)(575,276)(591,291)
	(607,307)(622,322)(637,337)
	(653,353)(669,368)(685,384)
	(702,401)(719,417)(736,433)
	(754,449)(772,465)(790,480)
	(808,495)(825,509)(842,522)
	(859,534)(875,545)(891,555)
	(907,564)(922,572)(937,579)
	(953,586)(969,591)(985,596)
	(1002,601)(1019,604)(1036,606)
	(1054,608)(1072,608)(1090,608)
	(1108,606)(1125,604)(1142,601)
	(1159,596)(1175,591)(1191,586)
	(1207,579)(1222,572)(1237,564)
	(1253,555)(1269,545)(1285,534)
	(1302,522)(1319,509)(1336,495)
	(1354,480)(1372,465)(1390,449)
	(1408,433)(1425,417)(1442,401)
	(1459,384)(1475,368)(1491,353)
	(1507,337)(1522,322)(1537,307)
	(1553,291)(1569,276)(1585,260)
	(1602,243)(1619,227)(1636,211)
	(1654,195)(1672,179)(1690,164)
	(1708,149)(1725,135)(1742,122)
	(1759,110)(1775,99)(1791,89)
	(1807,80)(1822,72)(1839,64)
	(1857,57)(1875,51)(1895,46)
	(1916,41)(1939,37)(1964,34)
	(1991,31)(2019,28)(2047,26)
	(2072,24)(2093,23)(2109,23)
	(2118,22)(2121,22)(2122,22)
\path(1672,472)(1675,475)(1682,481)
	(1693,491)(1707,504)(1724,518)
	(1742,533)(1759,547)(1775,559)
	(1790,569)(1805,578)(1819,586)
	(1833,592)(1847,597)(1860,601)
	(1874,605)(1889,608)(1906,610)
	(1925,612)(1946,614)(1970,616)
	(1996,618)(2022,619)(2049,620)
	(2074,621)(2094,621)(2109,622)
	(2118,622)(2121,622)(2122,622)
\path(772,172)(775,169)(782,163)
	(793,153)(807,140)(824,126)
	(842,111)(859,97)(875,85)
	(890,75)(905,66)(919,58)
	(933,52)(947,47)(962,42)
	(978,38)(995,35)(1013,32)
	(1032,31)(1052,29)(1072,29)
	(1092,29)(1112,31)(1131,32)
	(1149,35)(1166,38)(1182,42)
	(1197,47)(1211,52)(1225,58)
	(1239,66)(1254,75)(1269,85)
	(1285,97)(1302,111)(1320,126)
	(1337,140)(1351,153)(1362,163)
	(1369,169)(1372,172)
\path(22,622)(23,622)(26,622)
	(35,622)(50,621)(70,621)
	(95,620)(122,619)(148,618)
	(174,616)(198,614)(219,612)
	(238,610)(255,608)(270,605)
	(284,601)(297,597)(311,592)
	(325,586)(339,578)(354,569)
	(369,559)(385,547)(402,533)
	(420,518)(437,504)(451,491)
	(462,481)(469,475)(472,472)
\end{picture}
} } 
  \end{array}
\leftrightarrow
  \setlength{\unitlength}{0.2\standardunitlength}
  \begin{array}{c}
    {\begingroup\makeatletter\ifx\SetFigFont\undefined%
\gdef\SetFigFont#1#2#3#4#5{%
  \reset@font\fontsize{#1}{#2pt}%
  \fontfamily{#3}\fontseries{#4}\fontshape{#5}%
  \selectfont}%
\fi\endgroup%
{\renewcommand{\dashlinestretch}{30}
\begin{picture}(2144,659)(0,-10)
\thicklines
\path(22,22)(25,22)(32,22)
	(44,23)(62,24)(84,25)
	(110,26)(137,28)(166,31)
	(194,33)(221,36)(247,40)
	(272,43)(297,48)(321,53)
	(346,58)(371,65)(397,72)
	(416,77)(435,83)(456,90)
	(477,96)(499,103)(522,110)
	(546,118)(571,126)(597,134)
	(624,142)(651,149)(679,157)
	(707,165)(736,173)(765,180)
	(794,187)(823,193)(852,199)
	(881,204)(909,209)(937,213)
	(965,216)(992,219)(1019,220)
	(1046,222)(1072,222)(1098,222)
	(1125,220)(1152,219)(1179,216)
	(1207,213)(1235,209)(1263,204)
	(1292,199)(1321,193)(1350,187)
	(1379,180)(1408,173)(1437,165)
	(1465,157)(1493,149)(1520,142)
	(1547,134)(1573,126)(1598,118)
	(1622,110)(1645,103)(1667,96)
	(1688,90)(1709,83)(1728,77)
	(1747,72)(1773,65)(1798,58)
	(1823,53)(1847,48)(1872,43)
	(1897,40)(1923,36)(1950,33)
	(1978,31)(2007,28)(2034,26)
	(2060,25)(2082,24)(2100,23)
	(2112,22)(2119,22)(2122,22)
\path(2122,622)(2119,622)(2112,622)
	(2100,621)(2082,620)(2060,619)
	(2034,618)(2007,616)(1978,613)
	(1950,611)(1923,608)(1897,604)
	(1872,601)(1847,596)(1823,591)
	(1798,586)(1773,579)(1747,572)
	(1728,567)(1709,561)(1688,554)
	(1667,548)(1645,541)(1622,534)
	(1598,526)(1573,518)(1547,510)
	(1520,502)(1493,495)(1465,487)
	(1437,479)(1408,471)(1379,464)
	(1350,457)(1321,451)(1292,445)
	(1263,440)(1235,435)(1207,431)
	(1179,428)(1152,425)(1125,424)
	(1098,422)(1072,422)(1046,422)
	(1019,424)(992,425)(965,428)
	(937,431)(909,435)(881,440)
	(852,445)(823,451)(794,457)
	(765,464)(736,471)(707,479)
	(679,487)(651,495)(624,502)
	(597,510)(571,518)(546,526)
	(522,534)(499,541)(477,548)
	(456,554)(435,561)(416,567)
	(397,572)(371,579)(346,586)
	(321,591)(297,596)(272,601)
	(247,604)(221,608)(194,611)
	(166,613)(137,616)(110,618)
	(84,619)(62,620)(44,621)
	(32,622)(25,622)(22,622)
\end{picture}
} } 
  \end{array}
$ and (R3):
$
  \setlength{\unitlength}{0.15\standardunitlength}
  \begin{array}{c}
    {\begingroup\makeatletter\ifx\SetFigFont\undefined%
\gdef\SetFigFont#1#2#3#4#5{%
  \reset@font\fontsize{#1}{#2pt}%
  \fontfamily{#3}\fontseries{#4}\fontshape{#5}%
  \selectfont}%
\fi\endgroup%
{\renewcommand{\dashlinestretch}{30}
\begin{picture}(3044,3059)(0,-10)
\thicklines
\path(1297,2047)(1295,2049)(1290,2054)
	(1282,2063)(1270,2076)(1255,2093)
	(1235,2114)(1214,2137)(1191,2163)
	(1167,2190)(1144,2217)(1122,2243)
	(1101,2269)(1082,2293)(1065,2316)
	(1049,2338)(1035,2359)(1023,2380)
	(1012,2400)(1002,2420)(993,2439)
	(985,2459)(977,2480)(970,2501)
	(964,2523)(959,2547)(954,2572)
	(949,2599)(945,2628)(942,2659)
	(938,2693)(935,2729)(933,2767)
	(930,2806)(928,2846)(927,2885)
	(925,2920)(924,2952)(923,2978)
	(923,2998)(922,3011)(922,3019)(922,3022)
\path(2272,1147)(2275,1148)(2282,1150)
	(2293,1153)(2310,1158)(2331,1164)
	(2355,1171)(2381,1178)(2407,1186)
	(2433,1194)(2457,1202)(2480,1210)
	(2501,1218)(2521,1225)(2541,1233)
	(2559,1242)(2578,1250)(2597,1259)
	(2613,1268)(2629,1276)(2646,1285)
	(2664,1295)(2683,1306)(2703,1318)
	(2725,1332)(2749,1346)(2774,1362)
	(2801,1379)(2830,1398)(2860,1417)
	(2890,1436)(2919,1455)(2946,1472)
	(2969,1488)(2989,1500)(3004,1510)
	(3014,1517)(3019,1520)(3022,1522)
\path(1297,997)(1747,997)
\path(1747,1597)(1749,1595)(1752,1591)
	(1759,1584)(1768,1573)(1781,1558)
	(1797,1539)(1815,1517)(1836,1492)
	(1857,1466)(1878,1438)(1900,1410)
	(1920,1382)(1939,1354)(1957,1327)
	(1973,1300)(1988,1274)(2002,1248)
	(2014,1222)(2025,1196)(2035,1169)
	(2044,1142)(2052,1114)(2060,1084)
	(2065,1062)(2070,1038)(2074,1013)
	(2078,987)(2082,959)(2086,930)
	(2089,900)(2093,867)(2095,832)
	(2098,795)(2101,755)(2103,713)
	(2105,668)(2107,621)(2109,572)
	(2111,521)(2113,469)(2114,416)
	(2116,363)(2117,311)(2118,261)
	(2119,214)(2120,171)(2120,133)
	(2121,100)(2121,74)(2122,53)
	(2122,39)(2122,29)(2122,24)(2122,22)
\path(22,1522)(25,1520)(31,1517)
	(43,1510)(60,1500)(83,1488)
	(110,1472)(141,1455)(174,1436)
	(208,1417)(242,1398)(275,1379)
	(306,1362)(334,1346)(361,1332)
	(385,1318)(407,1306)(428,1295)
	(447,1285)(464,1276)(481,1268)
	(497,1259)(518,1249)(538,1239)
	(558,1230)(578,1222)(598,1213)
	(620,1204)(642,1195)(665,1186)
	(689,1177)(712,1169)(732,1161)
	(749,1155)(762,1151)(769,1148)(772,1147)
\path(922,22)(922,24)(922,28)
	(923,36)(923,49)(924,67)
	(925,90)(926,119)(928,154)
	(930,193)(933,237)(935,284)
	(939,335)(943,386)(947,439)
	(951,492)(956,545)(962,597)
	(968,647)(974,696)(981,743)
	(988,788)(995,831)(1003,873)
	(1012,912)(1021,951)(1031,987)
	(1041,1023)(1053,1058)(1065,1091)
	(1078,1124)(1092,1157)(1106,1190)
	(1122,1222)(1138,1253)(1154,1283)
	(1172,1314)(1190,1345)(1209,1377)
	(1229,1408)(1250,1440)(1272,1473)
	(1294,1505)(1317,1538)(1341,1572)
	(1366,1605)(1391,1639)(1417,1673)
	(1443,1707)(1469,1741)(1495,1774)
	(1522,1808)(1549,1842)(1575,1875)
	(1601,1907)(1627,1940)(1653,1971)
	(1678,2002)(1703,2033)(1727,2063)
	(1750,2092)(1772,2120)(1794,2148)
	(1815,2175)(1835,2201)(1854,2226)
	(1872,2251)(1890,2275)(1906,2299)
	(1922,2322)(1942,2353)(1961,2384)
	(1978,2414)(1994,2444)(2008,2474)
	(2022,2505)(2034,2536)(2045,2567)
	(2055,2600)(2064,2634)(2073,2670)
	(2080,2706)(2088,2743)(2094,2781)
	(2100,2819)(2105,2855)(2109,2890)
	(2113,2922)(2116,2950)(2118,2974)
	(2120,2992)(2121,3006)(2122,3015)
	(2122,3020)(2122,3022)
\end{picture}
} } 
  \end{array}
\leftrightarrow
  \setlength{\unitlength}{0.15\standardunitlength}
  \begin{array}{c}
    {\begingroup\makeatletter\ifx\SetFigFont\undefined%
\gdef\SetFigFont#1#2#3#4#5{%
  \reset@font\fontsize{#1}{#2pt}%
  \fontfamily{#3}\fontseries{#4}\fontshape{#5}%
  \selectfont}%
\fi\endgroup%
{\renewcommand{\dashlinestretch}{30}
\begin{picture}(3044,3059)(0,-10)
\thicklines
\path(1747,997)(1749,995)(1754,990)
	(1762,981)(1774,968)(1789,951)
	(1809,930)(1830,907)(1853,881)
	(1877,854)(1900,827)(1922,801)
	(1943,775)(1962,751)(1979,728)
	(1995,706)(2009,685)(2021,664)
	(2032,644)(2042,624)(2051,605)
	(2060,584)(2067,564)(2074,543)
	(2080,521)(2085,497)(2090,472)
	(2095,445)(2099,416)(2102,385)
	(2106,351)(2109,315)(2111,277)
	(2114,238)(2116,198)(2117,159)
	(2119,124)(2120,92)(2121,66)
	(2121,46)(2122,33)(2122,25)(2122,22)
\path(772,1897)(769,1896)(762,1894)
	(751,1891)(734,1886)(713,1880)
	(689,1873)(663,1866)(637,1858)
	(611,1850)(587,1842)(564,1834)
	(543,1826)(523,1819)(503,1811)
	(485,1802)(466,1794)(447,1784)
	(431,1776)(415,1768)(398,1759)
	(380,1749)(361,1738)(341,1726)
	(319,1712)(295,1698)(270,1682)
	(243,1665)(214,1646)(184,1627)
	(154,1608)(125,1589)(98,1572)
	(75,1556)(55,1544)(40,1534)
	(30,1527)(25,1524)(22,1522)
\path(1747,2047)(1297,2047)
\path(1297,1447)(1295,1449)(1292,1453)
	(1285,1460)(1276,1471)(1263,1486)
	(1247,1505)(1229,1527)(1208,1552)
	(1187,1578)(1166,1606)(1144,1634)
	(1124,1662)(1105,1690)(1087,1717)
	(1071,1744)(1056,1770)(1042,1796)
	(1030,1822)(1019,1848)(1009,1875)
	(1000,1902)(992,1930)(985,1959)
	(979,1982)(974,2006)(970,2031)
	(966,2057)(962,2085)(958,2114)
	(955,2144)(951,2177)(949,2212)
	(946,2249)(943,2289)(941,2331)
	(939,2376)(937,2423)(935,2472)
	(933,2523)(931,2575)(930,2628)
	(928,2681)(927,2733)(926,2783)
	(925,2830)(924,2873)(924,2911)
	(923,2944)(923,2970)(922,2991)
	(922,3005)(922,3015)(922,3020)(922,3022)
\path(3022,1522)(3019,1524)(3013,1527)
	(3001,1534)(2984,1544)(2961,1556)
	(2934,1572)(2903,1589)(2870,1608)
	(2836,1627)(2802,1646)(2769,1665)
	(2738,1682)(2710,1698)(2683,1712)
	(2659,1726)(2637,1738)(2616,1749)
	(2597,1759)(2580,1768)(2563,1776)
	(2547,1784)(2526,1795)(2506,1805)
	(2486,1814)(2466,1822)(2446,1831)
	(2424,1840)(2402,1849)(2379,1858)
	(2355,1867)(2332,1875)(2312,1883)
	(2295,1889)(2282,1893)(2275,1896)(2272,1897)
\path(2122,3022)(2122,3020)(2122,3016)
	(2121,3008)(2121,2995)(2120,2977)
	(2119,2954)(2118,2925)(2116,2890)
	(2114,2851)(2111,2807)(2109,2760)
	(2105,2709)(2101,2658)(2097,2605)
	(2093,2552)(2088,2499)(2082,2447)
	(2076,2397)(2070,2348)(2063,2301)
	(2056,2256)(2049,2213)(2041,2171)
	(2032,2132)(2023,2093)(2013,2057)
	(2003,2021)(1991,1986)(1979,1953)
	(1966,1920)(1952,1887)(1938,1854)
	(1922,1822)(1906,1791)(1890,1761)
	(1872,1730)(1854,1699)(1835,1667)
	(1815,1636)(1794,1604)(1772,1571)
	(1750,1539)(1727,1506)(1703,1472)
	(1678,1439)(1653,1405)(1627,1371)
	(1601,1337)(1575,1303)(1549,1270)
	(1522,1236)(1495,1202)(1469,1169)
	(1443,1137)(1417,1104)(1391,1073)
	(1366,1042)(1341,1011)(1317,981)
	(1294,952)(1272,924)(1250,896)
	(1229,869)(1209,843)(1190,818)
	(1172,793)(1154,769)(1138,745)
	(1122,722)(1102,691)(1083,660)
	(1066,630)(1050,600)(1036,570)
	(1022,539)(1010,508)(999,477)
	(989,444)(980,410)(971,374)
	(964,338)(956,301)(950,263)
	(944,225)(939,189)(935,154)
	(931,122)(928,94)(926,70)
	(924,52)(923,38)(922,29)
	(922,24)(922,22)
\end{picture}
} } 
  \end{array}
$. In the case of the
Kauffman bracket/Jones polynomial, this is done by reducing the Kauffman
bracket of the ``complicated side'' of each of these moves using the
rules $\la\backoverslash\ra=\la\hsmoothing\ra-q\la\smoothing\ra$
and $\la\bigcirc L\ra = (q+q^{-1})\la L\ra$ and then by
canceling terms until the ``easy side'' is reached. (Example:
{$\def\KB#1{{\left\la
  \setlength{\unitlength}{0.2\standardunitlength}
  \begin{array}{c}
    {\input figs/#1.tex } 
  \end{array}
\right\ra}}\KB{RT}=
\KB{RT0} - q \KB{RT1}= (q+q^{-1})\KB{RTResolved} - q \KB{RTResolved}=
q^{-1}\KB{RTResolved}$}). We do nearly the same in the case of the Khovanov
bracket. We first need to introduce a ``cancellation principle'' for chain
complexes:

\begin{lemma} \label{lem:Cancellation} Let $\calC$ be a chain complex
and let $\calC'\subset\calC$ be a sub chain complex.
\begin{itemize}
\item If $\calC'$ is acyclic (has no homology), then it can be
  ``canceled''. That is, in that case the homology $H(\calC)$
  of $\calC$ is equal to the homology $H(\calC/\calC')$ of
  $\calC/\calC'$.
\item Likewise, if $\calC/\calC'$ is acyclic then
  $H(\calC)=H(\calC')$.
\end{itemize}
\end{lemma}

\begin{proof} Both assertions follow easily from the long exact sequence
\[
  \ldots\longrightarrow H^r(\calC')\longrightarrow
  H^r(\calC)\longrightarrow H^r(\calC/\calC')
  \longrightarrow H^{r+1}(\calC')\longrightarrow\ldots
\]
associated with the short exact sequence
$0\longrightarrow\calC'\longrightarrow\calC
\longrightarrow\calC/\calC'\longrightarrow 0$.
\end{proof}

\subsubsection{Invariance under (R1).} In computing
$\calH(
  \setlength{\unitlength}{0.2\standardunitlength}
  \begin{array}{c}
    {\begingroup\makeatletter\ifx\SetFigFont\undefined%
\gdef\SetFigFont#1#2#3#4#5{%
  \reset@font\fontsize{#1}{#2pt}%
  \fontfamily{#3}\fontseries{#4}\fontshape{#5}%
  \selectfont}%
\fi\endgroup%
{\renewcommand{\dashlinestretch}{30}
\begin{picture}(1244,859)(0,-10)
\thicklines
\path(1222,22)(1219,22)(1212,23)
	(1200,24)(1183,25)(1161,28)
	(1137,31)(1110,35)(1083,39)
	(1056,45)(1030,51)(1006,57)
	(982,65)(960,74)(938,84)
	(916,95)(895,108)(872,122)
	(855,134)(837,147)(818,161)
	(799,176)(779,192)(759,209)
	(739,227)(718,246)(697,266)
	(677,286)(656,307)(636,329)
	(616,351)(597,373)(579,395)
	(562,416)(546,438)(531,459)
	(518,479)(506,499)(495,518)
	(486,537)(478,555)(472,572)
	(467,589)(463,606)(461,622)
	(460,638)(460,654)(461,669)
	(464,684)(467,699)(472,713)
	(478,726)(485,739)(493,751)
	(501,763)(511,773)(521,782)
	(531,791)(542,798)(553,805)
	(565,810)(576,814)(588,818)
	(599,820)(611,822)(622,822)
	(637,821)(653,819)(668,815)
	(684,809)(700,802)(716,794)
	(731,784)(747,773)(762,762)
	(776,749)(789,736)(801,723)
	(812,710)(821,696)(830,683)
	(837,671)(843,659)(847,647)
	(851,630)(852,614)(850,598)
	(845,580)(836,562)(824,542)
	(811,521)(797,502)(784,487)
	(776,477)(773,473)(772,472)
\path(472,172)(469,169)(462,163)
	(451,153)(437,140)(420,126)
	(402,111)(385,97)(369,85)
	(354,75)(339,66)(325,58)
	(311,52)(297,47)(284,43)
	(270,39)(255,36)(238,34)
	(219,32)(198,30)(174,28)
	(148,26)(122,25)(95,24)
	(70,23)(50,23)(35,22)
	(26,22)(23,22)(22,22)
\end{picture}
} } 
  \end{array}
)$ we encounter the complex
\begin{equation} \label{eq:RTComplex}
  \calC = \KhB{RT}{0.2} = \left(
    \KhB{RT0}{0.2}\overset{m}{\longrightarrow}\KhB{RT1}{0.2}\{1\}
  \right).
\end{equation}
(Each of the terms in this complex is itself a complex, coming from a
whole cube of spaces and maps. We implicitly ``flatten'' such complexes
of complexes to single complexes as in Section~\ref{subsec:Maps} without
further comment). The complex in Equation~\eqref{eq:RTComplex} has a natural
subcomplex
\[ \calC' = \left(
    \KhB{RT0}{0.2}_{v_+}\overset{m}{\longrightarrow}\KhB{RT1}{0.2}\{1\}
  \right)
\]
We need to pause to explain the notation. Recall that $\llbracket
L\rrbracket$ is a direct sum over the smoothings of $L$ of tensor
powers of $V$, with one tensor factor corresponding to  each cycle
in any given smoothing. Such tensor powers can be viewed as spaces
of linear combinations of marked smoothings of $L$, where each
cycle in any smoothing of $L$ is marked by an element of $V$. For
$L=
  \setlength{\unitlength}{0.2\standardunitlength}
  \begin{array}{c}
    {\begingroup\makeatletter\ifx\SetFigFont\undefined%
\gdef\SetFigFont#1#2#3#4#5{%
  \reset@font\fontsize{#1}{#2pt}%
  \fontfamily{#3}\fontseries{#4}\fontshape{#5}%
  \selectfont}%
\fi\endgroup%
{\renewcommand{\dashlinestretch}{30}
\begin{picture}(1244,859)(0,-10)
\thicklines
\path(22,22)(23,22)(26,22)
	(35,23)(51,23)(72,24)
	(97,26)(125,28)(153,31)
	(180,34)(205,37)(228,41)
	(249,46)(269,51)(287,57)
	(305,64)(322,72)(337,80)
	(353,89)(369,99)(385,109)
	(402,120)(419,131)(436,142)
	(454,154)(472,165)(490,176)
	(508,186)(525,195)(542,203)
	(559,210)(575,215)(591,219)
	(607,221)(622,222)(637,221)
	(653,219)(669,215)(685,210)
	(702,203)(719,195)(736,186)
	(754,176)(772,165)(790,154)
	(808,142)(825,131)(842,120)
	(859,109)(875,99)(891,89)
	(907,80)(922,72)(939,64)
	(957,57)(975,51)(995,46)
	(1016,41)(1039,37)(1064,34)
	(1091,31)(1119,28)(1147,26)
	(1172,24)(1193,23)(1209,23)
	(1218,22)(1221,22)(1222,22)
\path(422,622)(423,608)(425,594)
	(428,580)(433,566)(438,552)
	(446,538)(454,525)(464,511)
	(474,498)(486,486)(498,474)
	(511,464)(525,454)(538,446)
	(552,438)(566,433)(580,428)
	(594,425)(608,423)(622,422)
	(636,423)(650,425)(664,428)
	(678,433)(692,438)(706,446)
	(719,454)(733,464)(746,474)
	(758,486)(770,498)(780,511)
	(790,525)(798,538)(806,552)
	(811,566)(816,580)(819,594)
	(821,608)(822,622)(821,636)
	(819,650)(816,664)(811,678)
	(806,692)(798,706)(790,719)
	(780,733)(770,746)(758,758)
	(746,770)(733,780)(719,790)
	(706,798)(692,806)(678,811)
	(664,816)(650,819)(636,821)
	(622,822)(608,821)(594,819)
	(580,816)(566,811)(552,806)
	(538,798)(525,790)(511,780)
	(498,770)(486,758)(474,746)
	(464,733)(454,719)(446,706)
	(438,692)(433,678)(428,664)
	(425,650)(423,636)(422,622)
\end{picture}
} } 
  \end{array}
$ all smoothings have one special cycle, the
one appearing within the icon $
  \setlength{\unitlength}{0.2\standardunitlength}
  \begin{array}{c}
    {\begingroup\makeatletter\ifx\SetFigFont\undefined%
\gdef\SetFigFont#1#2#3#4#5{%
  \reset@font\fontsize{#1}{#2pt}%
  \fontfamily{#3}\fontseries{#4}\fontshape{#5}%
  \selectfont}%
\fi\endgroup%
{\renewcommand{\dashlinestretch}{30}
\begin{picture}(1244,859)(0,-10)
\thicklines
\path(22,22)(23,22)(26,22)
	(35,23)(51,23)(72,24)
	(97,26)(125,28)(153,31)
	(180,34)(205,37)(228,41)
	(249,46)(269,51)(287,57)
	(305,64)(322,72)(337,80)
	(353,89)(369,99)(385,109)
	(402,120)(419,131)(436,142)
	(454,154)(472,165)(490,176)
	(508,186)(525,195)(542,203)
	(559,210)(575,215)(591,219)
	(607,221)(622,222)(637,221)
	(653,219)(669,215)(685,210)
	(702,203)(719,195)(736,186)
	(754,176)(772,165)(790,154)
	(808,142)(825,131)(842,120)
	(859,109)(875,99)(891,89)
	(907,80)(922,72)(939,64)
	(957,57)(975,51)(995,46)
	(1016,41)(1039,37)(1064,34)
	(1091,31)(1119,28)(1147,26)
	(1172,24)(1193,23)(1209,23)
	(1218,22)(1221,22)(1222,22)
\path(422,622)(423,608)(425,594)
	(428,580)(433,566)(438,552)
	(446,538)(454,525)(464,511)
	(474,498)(486,486)(498,474)
	(511,464)(525,454)(538,446)
	(552,438)(566,433)(580,428)
	(594,425)(608,423)(622,422)
	(636,423)(650,425)(664,428)
	(678,433)(692,438)(706,446)
	(719,454)(733,464)(746,474)
	(758,486)(770,498)(780,511)
	(790,525)(798,538)(806,552)
	(811,566)(816,580)(819,594)
	(821,608)(822,622)(821,636)
	(819,650)(816,664)(811,678)
	(806,692)(798,706)(790,719)
	(780,733)(770,746)(758,758)
	(746,770)(733,780)(719,790)
	(706,798)(692,806)(678,811)
	(664,816)(650,819)(636,821)
	(622,822)(608,821)(594,819)
	(580,816)(566,811)(552,806)
	(538,798)(525,790)(511,780)
	(498,770)(486,758)(474,746)
	(464,733)(454,719)(446,706)
	(438,692)(433,678)(428,664)
	(425,650)(423,636)(422,622)
\end{picture}
} } 
  \end{array}
$. The subscript
$v_+$ in $\KhB{RT0}{0.2}_{v_+}$ means ``the subspace of $\KhB{RT0}{0.2}$
in which the special cycle is always marked $v_+$''.

It is easy to check that $\calC'$ is indeed a subcomplex of $\calC$, and
as $v_+$ is a unit for the product $m$ (see~\eqref{eq:mdef}), $\calC'$ is
acyclic. Thus by Lemma~\ref{lem:Cancellation} we are reduced to studying
the quotient complex
\[ \calC/\calC' = \left(\KhB{RT0}{0.2}_{/v_+=0}\rightarrow 0\right) \]
where the subscript ``$/v_+=0$'' means ``mod out (within the
tensor factor corresponding to the special cycle) by $v_+=0$''. But
$V/(v_+=0)$ is one dimensional and generated by $v_-$, and hence apart
from a shift in degrees, $\KhB{RT0}{0.2}_{/v_+=0}$ is isomorphic to
$\KhB{RTResolved}{0.2}$. The reader may verify that this shift precisely
gets canceled by the shifts $[-n_-]\{n_+-2n_-\}$ in the definition of
$\calC(L)$ from $\llbracket L\rrbracket$. \qed

\subsubsection{Invariance under (R2), first proof.} In computing
$\calH(
  \setlength{\unitlength}{0.2\standardunitlength}
  \begin{array}{c}
    {\begingroup\makeatletter\ifx\SetFigFont\undefined%
\gdef\SetFigFont#1#2#3#4#5{%
  \reset@font\fontsize{#1}{#2pt}%
  \fontfamily{#3}\fontseries{#4}\fontshape{#5}%
  \selectfont}%
\fi\endgroup%
{\renewcommand{\dashlinestretch}{30}
\begin{picture}(2144,659)(0,-10)
\thicklines
\path(22,22)(23,22)(26,22)
	(35,23)(51,23)(72,24)
	(97,26)(125,28)(153,31)
	(180,34)(205,37)(228,41)
	(249,46)(269,51)(287,57)
	(305,64)(322,72)(337,80)
	(353,89)(369,99)(385,110)
	(402,122)(419,135)(436,149)
	(454,164)(472,179)(490,195)
	(508,211)(525,227)(542,243)
	(559,260)(575,276)(591,291)
	(607,307)(622,322)(637,337)
	(653,353)(669,368)(685,384)
	(702,401)(719,417)(736,433)
	(754,449)(772,465)(790,480)
	(808,495)(825,509)(842,522)
	(859,534)(875,545)(891,555)
	(907,564)(922,572)(937,579)
	(953,586)(969,591)(985,596)
	(1002,601)(1019,604)(1036,606)
	(1054,608)(1072,608)(1090,608)
	(1108,606)(1125,604)(1142,601)
	(1159,596)(1175,591)(1191,586)
	(1207,579)(1222,572)(1237,564)
	(1253,555)(1269,545)(1285,534)
	(1302,522)(1319,509)(1336,495)
	(1354,480)(1372,465)(1390,449)
	(1408,433)(1425,417)(1442,401)
	(1459,384)(1475,368)(1491,353)
	(1507,337)(1522,322)(1537,307)
	(1553,291)(1569,276)(1585,260)
	(1602,243)(1619,227)(1636,211)
	(1654,195)(1672,179)(1690,164)
	(1708,149)(1725,135)(1742,122)
	(1759,110)(1775,99)(1791,89)
	(1807,80)(1822,72)(1839,64)
	(1857,57)(1875,51)(1895,46)
	(1916,41)(1939,37)(1964,34)
	(1991,31)(2019,28)(2047,26)
	(2072,24)(2093,23)(2109,23)
	(2118,22)(2121,22)(2122,22)
\path(1672,472)(1675,475)(1682,481)
	(1693,491)(1707,504)(1724,518)
	(1742,533)(1759,547)(1775,559)
	(1790,569)(1805,578)(1819,586)
	(1833,592)(1847,597)(1860,601)
	(1874,605)(1889,608)(1906,610)
	(1925,612)(1946,614)(1970,616)
	(1996,618)(2022,619)(2049,620)
	(2074,621)(2094,621)(2109,622)
	(2118,622)(2121,622)(2122,622)
\path(772,172)(775,169)(782,163)
	(793,153)(807,140)(824,126)
	(842,111)(859,97)(875,85)
	(890,75)(905,66)(919,58)
	(933,52)(947,47)(962,42)
	(978,38)(995,35)(1013,32)
	(1032,31)(1052,29)(1072,29)
	(1092,29)(1112,31)(1131,32)
	(1149,35)(1166,38)(1182,42)
	(1197,47)(1211,52)(1225,58)
	(1239,66)(1254,75)(1269,85)
	(1285,97)(1302,111)(1320,126)
	(1337,140)(1351,153)(1362,163)
	(1369,169)(1372,172)
\path(22,622)(23,622)(26,622)
	(35,622)(50,621)(70,621)
	(95,620)(122,619)(148,618)
	(174,616)(198,614)(219,612)
	(238,610)(255,608)(270,605)
	(284,601)(297,597)(311,592)
	(325,586)(339,578)(354,569)
	(369,559)(385,547)(402,533)
	(420,518)(437,504)(451,491)
	(462,481)(469,475)(472,472)
\end{picture}
} } 
  \end{array}
)$ we encounter the complex $\calC$
of Figure~1. This complex has a subcomplex $\calC'$ (see
Figure~1), which is clearly acyclic.  The quotient complex
$\calC/\calC'$ (see Figure~1) has a subcomplex $\calC''$
(see Figure~1), and the quotient $(\calC/\calC')/\calC''$ (see
Figure~1) is acyclic because modulo $v_+=0$, the map $\Delta$
is an isomorphism. Hence using both parts of Lemma~\ref{lem:Cancellation}
we find that $H(\calC)=H(\calC/\calC')=H(\calC'')$. But up to shifts in
degree and height, $\calC''$ is just $\KhB{TMResolved}{0.2}$. Again,
these shifts get canceled by the shifts $[-n_-]\{n_+-2n_-\}$ in the
definition of $\calC(L)$ from $\llbracket L\rrbracket$. \qed

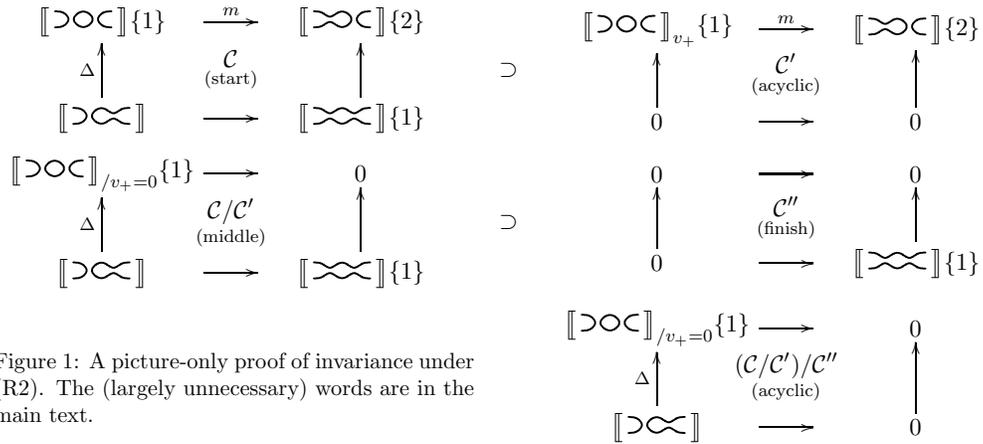
\begin{figure}
{\centering\resizebox*{\textwidth}{!}{\parbox{6.2in}{
\def\entry#1{{\parbox{1.2in}{\centering $#1$}}}
\begin{eqnarray*}
  \begin{array}{c}\xymatrix{
    \entry{\KhB{TM01}{0.2}\{1\}} \ar[r]^m
      \ar@{}[rd]|{\displaystyle \underset{\text{(start)}}{\calC}}
    & \entry{\KhB{TM11}{0.2}\{2\}} \\
    \entry{\KhB{TM00}{0.2}} \ar[u]^\Delta \ar[r]
    & \entry{\KhB{TM10}{0.2}\{1\}} \ar[u]
  }\end{array}
  & \supset & \begin{array}{c}\xymatrix{
    \entry{\KhB{TM01}{0.2}_{v_+}\{1\}} \ar[r]^m
      \ar@{}[rd]|{\displaystyle \underset{\text{(acyclic)}}{\calC'}}
    & \entry{\KhB{TM11}{0.2}\{2\}} \\
    \entry{0} \ar[u] \ar[r]
    & \entry{0} \ar[u]
  }\end{array}
  \\ \begin{array}{c}\xymatrix{
    \entry{\KhB{TM01}{0.2}_{/v_+=0}\{1\}} \ar[r]
      \ar@{}[rd]|{\displaystyle \underset{\text{(middle)}}{\calC/\calC'}}
    & \entry{0} \\
    \entry{\KhB{TM00}{0.2}} \ar[u]^\Delta \ar[r]
    & \entry{\KhB{TM10}{0.2}\{1\}} \ar[u]
  }\end{array}
  & \supset & \begin{array}{c}\xymatrix{
    \entry{0} \ar[r]
      \ar@{}[rd]|{\displaystyle \underset{\text{(finish)}}{\calC''}}
    & \entry{0} \\
    \entry{0} \ar[u] \ar[r]
    & \entry{\KhB{TM10}{0.2}\{1\}} \ar[u]
  }\end{array}
  \\
  \parbox{3in}{\caption{
    A picture-only proof of invariance under {\rm(R2)}. The (largely
    unnecessary) words are in the main text.
  }}
  & \quad & \begin{array}{c}\xymatrix{
    \entry{\KhB{TM01}{0.2}_{/v_+=0}\{1\}} \ar[r]
      \ar@{}[rd]|{\displaystyle
        \underset{\text{(acyclic)}}{(\calC/\calC')/\calC''}
      }
    & \entry{0} \\
    \entry{\KhB{TM00}{0.2}} \ar[u]^\Delta \ar[r]
    & \entry{0} \ar[u]
  }\end{array}
\end{eqnarray*}
}}}
\end{figure}

\subsubsection{Invariance under (R3), unsuccessful attempt.} For the
Kauffman bracket, invariance under (R3) follows from invariance under
(R2). Indeed, assuming relations of the form $\la\bigcirc L\ra = d\la
L\ra$ and $\la\backoverslash\ra = A \la\hsmoothing\ra + B\la\smoothing\ra$
the move (R3) follows from (R2) without imposing any constraints on $A$,
$B$ and $d$ (beyond those that are necessary for (R2) to hold):
\[ 
  \def\sKB#1{{\left\la
    \hspace{-1mm}
  \setlength{\unitlength}{0.12\standardunitlength}
  \begin{array}{c}
    {\input figs/#1.tex } 
  \end{array}
\hspace{-1mm}
  \right\ra}}
  \sKB{Du} = A\sKB{D0} + B\sKB{D1} \overset{\text{(R2)}}{=}
    A \sKB{U0} + B \sKB{U1} = \sKB{Uu}
\]
The case of the Khovanov bracket is unfortunately not as lucky. Invariance
under (R2) does play a key role, but more is needed. Let us see how it
works.

If we fully smooth the two sides of (R3), we get the following two cubes of
complexes (to save space we suppress the Khovanov bracket notation
$\llbracket\cdot\rrbracket$ and the degree shifts $\{\cdot\}$):

\begin{equation} \label{eq:R3Cubes}
  \def\neg{\hspace{-2mm}}
  \def\incDU#1{
    \parbox{10mm}{\includegraphics[width=10mm]{figs/#1.eps}}
  }
  \def\D#1{\neg
  \setlength{\unitlength}{0.125\standardunitlength}
  \begin{array}{c}
    {\input figs/D#1.tex } 
  \end{array}
\neg}
  \def\U#1{\neg
  \setlength{\unitlength}{0.125\standardunitlength}
  \begin{array}{c}
    {\input figs/U#1.tex } 
  \end{array}
\neg}
  \begin{array}{c} \xymatrix@R=4mm@C=4mm{
    \incDU{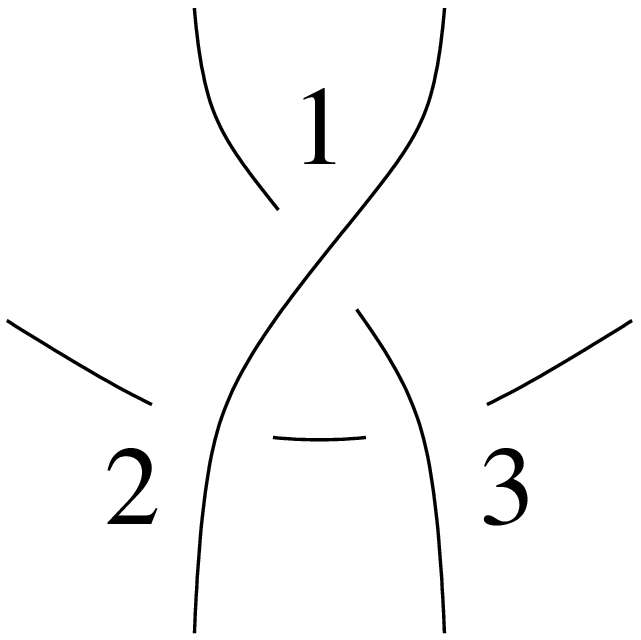} & \D{101} \ar@{->}[rr]^m \ar@{<-}'[d][dd]
      & & \D{111} \ar@{<-}[dd] \\
    \D{100} \ar@{->}[ur]^\Delta \ar@{->}[rr] \ar@{<-}[dd] &
      & \D{110} \ar@{->}[ur] \ar@{<-}[dd] \\
    & \D{001} \ar@{->}'[r][rr] & & \D{011} \\
    \D{000} \ar@{->}[rr] \ar@{->}[ur] & & \D{010} \ar@{->}[ur]
  } \end{array}
  \qquad
  \begin{array}{c} \xymatrix@R=4mm@C=4mm{
    \incDU{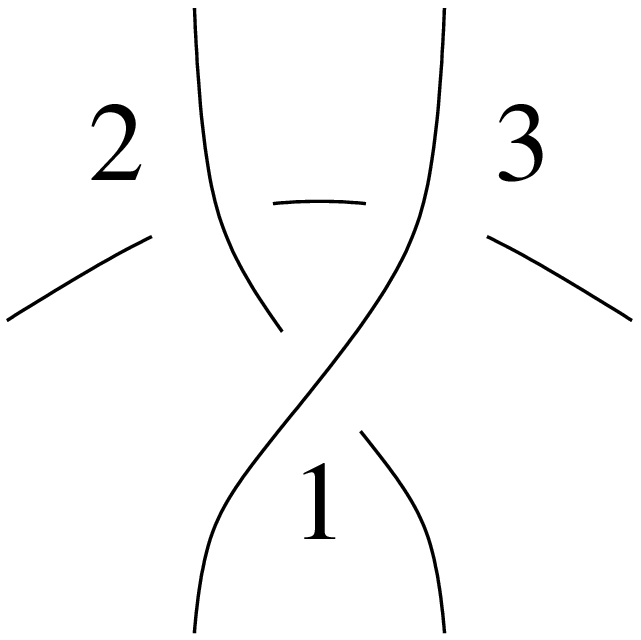} & \U{101} \ar@{->}[rr] \ar@{<-}'[d][dd]
      & & \U{111} \ar@{<-}[dd] \\
    \U{100} \ar@{->}[ur] \ar@{->}[rr]^(0.7)\Delta \ar@{<-}[dd] &
      & \U{110} \ar@{->}[ur]^m \ar@{<-}[dd] \\
    & \U{001} \ar@{->}'[r][rr] & & \U{011} \\
    \U{000} \ar@{->}[rr] \ar@{->}[ur] & & \U{010} \ar@{->}[ur]
  } \end{array}
\end{equation}

The bottom layers of these two cubes correspond to the partial smoothings
$
  \setlength{\unitlength}{0.12\standardunitlength}
  \begin{array}{c}
    {\begingroup\makeatletter\ifx\SetFigFont\undefined%
\gdef\SetFigFont#1#2#3#4#5{%
  \reset@font\fontsize{#1}{#2pt}%
  \fontfamily{#3}\fontseries{#4}\fontshape{#5}%
  \selectfont}%
\fi\endgroup%
{\renewcommand{\dashlinestretch}{30}
\begin{picture}(3044,3059)(0,-10)
\thicklines
\path(922,22)(922,25)(922,33)
	(922,46)(922,66)(922,92)
	(923,124)(923,159)(923,198)
	(924,238)(924,277)(925,315)
	(925,351)(926,385)(927,416)
	(927,445)(928,472)(929,497)
	(930,521)(932,543)(933,564)
	(935,584)(936,607)(938,629)
	(941,650)(943,672)(946,694)
	(949,716)(952,738)(956,760)
	(960,782)(965,804)(969,826)
	(974,847)(980,868)(985,888)
	(991,907)(997,926)(1003,945)
	(1009,962)(1015,980)(1022,997)
	(1029,1014)(1037,1032)(1044,1049)
	(1053,1068)(1062,1086)(1072,1105)
	(1082,1125)(1092,1145)(1103,1166)
	(1115,1186)(1126,1207)(1138,1228)
	(1150,1248)(1163,1269)(1175,1289)
	(1187,1308)(1199,1328)(1211,1347)
	(1223,1366)(1235,1384)(1247,1404)
	(1259,1423)(1271,1443)(1284,1463)
	(1297,1485)(1309,1506)(1322,1529)
	(1334,1552)(1346,1575)(1356,1599)
	(1366,1622)(1375,1646)(1383,1670)
	(1389,1693)(1394,1715)(1398,1738)
	(1400,1759)(1401,1780)(1400,1801)
	(1397,1822)(1393,1841)(1388,1860)
	(1381,1879)(1374,1899)(1364,1919)
	(1354,1940)(1342,1962)(1330,1984)
	(1316,2007)(1302,2030)(1287,2054)
	(1271,2078)(1255,2102)(1239,2127)
	(1223,2151)(1207,2176)(1192,2200)
	(1177,2224)(1162,2248)(1148,2273)
	(1135,2297)(1122,2322)(1112,2343)
	(1102,2365)(1092,2388)(1082,2412)
	(1073,2437)(1064,2463)(1055,2491)
	(1045,2521)(1036,2554)(1026,2588)
	(1016,2625)(1006,2663)(996,2704)
	(986,2745)(976,2787)(966,2828)
	(957,2868)(948,2905)(941,2937)
	(935,2965)(930,2987)(926,3003)
	(924,3014)(923,3019)(922,3022)
\path(2122,3022)(2121,3019)(2120,3014)
	(2118,3003)(2114,2987)(2109,2965)
	(2103,2937)(2096,2905)(2087,2868)
	(2078,2828)(2068,2787)(2058,2745)
	(2048,2704)(2038,2663)(2028,2625)
	(2018,2588)(2008,2554)(1999,2521)
	(1989,2491)(1980,2463)(1971,2437)
	(1962,2412)(1952,2388)(1942,2365)
	(1932,2343)(1922,2322)(1909,2297)
	(1896,2273)(1882,2248)(1867,2224)
	(1852,2200)(1837,2176)(1821,2151)
	(1805,2127)(1789,2102)(1773,2078)
	(1757,2054)(1742,2030)(1728,2007)
	(1714,1984)(1702,1962)(1690,1940)
	(1680,1919)(1670,1899)(1663,1879)
	(1656,1860)(1651,1841)(1647,1822)
	(1644,1801)(1643,1780)(1644,1759)
	(1646,1738)(1650,1715)(1655,1693)
	(1661,1670)(1669,1646)(1678,1622)
	(1688,1599)(1698,1575)(1710,1552)
	(1722,1529)(1735,1506)(1747,1485)
	(1760,1463)(1773,1443)(1785,1423)
	(1797,1404)(1810,1384)(1821,1366)
	(1833,1347)(1845,1328)(1857,1308)
	(1869,1289)(1881,1269)(1894,1248)
	(1906,1228)(1918,1207)(1929,1186)
	(1941,1166)(1952,1145)(1962,1125)
	(1972,1105)(1982,1086)(1991,1068)
	(2000,1049)(2007,1032)(2015,1014)
	(2022,997)(2029,980)(2035,962)
	(2041,945)(2047,926)(2053,907)
	(2059,888)(2064,868)(2070,847)
	(2075,826)(2079,804)(2084,782)
	(2088,760)(2092,738)(2095,716)
	(2098,694)(2101,672)(2103,650)
	(2106,629)(2108,607)(2110,584)
	(2111,564)(2112,543)(2114,521)
	(2115,497)(2116,472)(2117,445)
	(2117,416)(2118,385)(2119,351)
	(2119,315)(2120,277)(2120,238)
	(2121,198)(2121,159)(2121,124)
	(2122,92)(2122,66)(2122,46)
	(2122,33)(2122,25)(2122,22)
\path(22,1522)(25,1520)(31,1517)
	(43,1510)(60,1500)(83,1488)
	(110,1472)(141,1455)(174,1436)
	(208,1417)(242,1398)(275,1379)
	(306,1362)(334,1346)(361,1332)
	(385,1318)(407,1306)(428,1295)
	(447,1285)(464,1276)(481,1268)
	(497,1259)(518,1249)(538,1239)
	(558,1230)(578,1222)(598,1213)
	(620,1204)(642,1195)(665,1186)
	(689,1177)(712,1169)(732,1161)
	(749,1155)(762,1151)(769,1148)(772,1147)
\path(1297,997)(1747,997)
\path(2272,1147)(2275,1148)(2282,1150)
	(2293,1153)(2310,1158)(2331,1164)
	(2355,1171)(2381,1178)(2407,1186)
	(2433,1194)(2457,1202)(2480,1210)
	(2501,1218)(2521,1225)(2541,1233)
	(2559,1242)(2578,1250)(2597,1259)
	(2613,1268)(2629,1276)(2646,1285)
	(2664,1295)(2683,1306)(2703,1318)
	(2725,1332)(2749,1346)(2774,1362)
	(2801,1379)(2830,1398)(2860,1417)
	(2890,1436)(2919,1455)(2946,1472)
	(2969,1488)(2989,1500)(3004,1510)
	(3014,1517)(3019,1520)(3022,1522)
\end{picture}
} } 
  \end{array}
$ and $
  \setlength{\unitlength}{0.12\standardunitlength}
  \begin{array}{c}
    {\begingroup\makeatletter\ifx\SetFigFont\undefined%
\gdef\SetFigFont#1#2#3#4#5{%
  \reset@font\fontsize{#1}{#2pt}%
  \fontfamily{#3}\fontseries{#4}\fontshape{#5}%
  \selectfont}%
\fi\endgroup%
{\renewcommand{\dashlinestretch}{30}
\begin{picture}(3044,3059)(0,-10)
\thicklines
\path(922,3022)(922,3019)(922,3011)
	(922,2998)(922,2978)(922,2952)
	(923,2920)(923,2885)(923,2846)
	(924,2806)(924,2767)(925,2729)
	(925,2693)(926,2659)(927,2628)
	(927,2599)(928,2572)(929,2547)
	(930,2523)(932,2501)(933,2480)
	(935,2459)(936,2437)(938,2415)
	(941,2394)(943,2372)(946,2350)
	(949,2328)(952,2306)(956,2284)
	(960,2262)(965,2240)(969,2218)
	(974,2197)(980,2176)(985,2156)
	(991,2137)(997,2118)(1003,2099)
	(1009,2082)(1015,2064)(1022,2047)
	(1029,2030)(1037,2012)(1044,1995)
	(1053,1976)(1062,1958)(1072,1939)
	(1082,1919)(1092,1899)(1103,1878)
	(1115,1858)(1126,1837)(1138,1816)
	(1150,1796)(1163,1775)(1175,1755)
	(1187,1736)(1199,1716)(1211,1697)
	(1223,1678)(1235,1659)(1247,1640)
	(1259,1621)(1271,1601)(1284,1581)
	(1297,1559)(1309,1538)(1322,1515)
	(1334,1492)(1346,1469)(1356,1445)
	(1366,1422)(1375,1398)(1383,1374)
	(1389,1351)(1394,1329)(1398,1306)
	(1400,1285)(1401,1264)(1400,1243)
	(1397,1222)(1393,1203)(1388,1184)
	(1381,1165)(1374,1145)(1364,1125)
	(1354,1104)(1342,1082)(1330,1060)
	(1316,1037)(1302,1014)(1287,990)
	(1271,966)(1255,942)(1239,917)
	(1223,893)(1207,868)(1192,844)
	(1177,820)(1162,796)(1148,771)
	(1135,747)(1122,722)(1112,701)
	(1102,679)(1092,656)(1082,632)
	(1073,607)(1064,581)(1055,553)
	(1045,523)(1036,490)(1026,456)
	(1016,419)(1006,381)(996,340)
	(986,299)(976,257)(966,216)
	(957,176)(948,139)(941,107)
	(935,79)(930,57)(926,41)
	(924,30)(923,25)(922,22)
\path(2122,22)(2121,25)(2120,30)
	(2118,41)(2114,57)(2109,79)
	(2103,107)(2096,139)(2087,176)
	(2078,216)(2068,257)(2058,299)
	(2048,340)(2038,381)(2028,419)
	(2018,456)(2008,490)(1999,523)
	(1989,553)(1980,581)(1971,607)
	(1962,632)(1952,656)(1942,679)
	(1932,701)(1922,722)(1909,747)
	(1896,771)(1882,796)(1867,820)
	(1852,844)(1837,868)(1821,893)
	(1805,917)(1789,942)(1773,966)
	(1757,990)(1742,1014)(1728,1037)
	(1714,1060)(1702,1082)(1690,1104)
	(1680,1125)(1670,1145)(1663,1165)
	(1656,1184)(1651,1203)(1647,1222)
	(1644,1243)(1643,1264)(1644,1285)
	(1646,1306)(1650,1329)(1655,1351)
	(1661,1374)(1669,1398)(1678,1422)
	(1688,1445)(1698,1469)(1710,1492)
	(1722,1515)(1735,1538)(1747,1559)
	(1760,1581)(1773,1601)(1785,1621)
	(1797,1640)(1810,1659)(1821,1678)
	(1833,1697)(1845,1716)(1857,1736)
	(1869,1755)(1881,1775)(1894,1796)
	(1906,1816)(1918,1837)(1929,1858)
	(1941,1878)(1952,1899)(1962,1919)
	(1972,1939)(1982,1958)(1991,1976)
	(2000,1995)(2007,2012)(2015,2030)
	(2022,2047)(2029,2064)(2035,2082)
	(2041,2099)(2047,2118)(2053,2137)
	(2059,2156)(2064,2176)(2070,2197)
	(2075,2218)(2079,2240)(2084,2262)
	(2088,2284)(2092,2306)(2095,2328)
	(2098,2350)(2101,2372)(2103,2394)
	(2106,2415)(2108,2437)(2110,2459)
	(2111,2480)(2112,2501)(2114,2523)
	(2115,2547)(2116,2572)(2117,2599)
	(2117,2628)(2118,2659)(2119,2693)
	(2119,2729)(2120,2767)(2120,2806)
	(2121,2846)(2121,2885)(2121,2920)
	(2122,2952)(2122,2978)(2122,2998)
	(2122,3011)(2122,3019)(2122,3022)
\path(22,1522)(25,1524)(31,1527)
	(43,1534)(60,1544)(83,1556)
	(110,1572)(141,1589)(174,1608)
	(208,1627)(242,1646)(275,1665)
	(306,1682)(334,1698)(361,1712)
	(385,1726)(407,1738)(428,1749)
	(447,1759)(464,1768)(481,1776)
	(497,1784)(518,1795)(538,1805)
	(558,1814)(578,1822)(598,1831)
	(620,1840)(642,1849)(665,1858)
	(689,1867)(712,1875)(732,1883)
	(749,1889)(762,1893)(769,1896)(772,1897)
\path(1297,2047)(1747,2047)
\path(2272,1897)(2275,1896)(2282,1894)
	(2293,1891)(2310,1886)(2331,1880)
	(2355,1873)(2381,1866)(2407,1858)
	(2433,1850)(2457,1842)(2480,1834)
	(2501,1826)(2521,1819)(2541,1811)
	(2559,1802)(2578,1794)(2597,1784)
	(2613,1776)(2629,1768)(2646,1759)
	(2664,1749)(2683,1738)(2703,1726)
	(2725,1712)(2749,1698)(2774,1682)
	(2801,1665)(2830,1646)(2860,1627)
	(2890,1608)(2919,1589)(2946,1572)
	(2969,1556)(2989,1544)(3004,1534)
	(3014,1527)(3019,1524)(3022,1522)
\end{picture}
} } 
  \end{array}
$ and are therefore
isomorphic.  The top layers correspond to $
  \setlength{\unitlength}{0.12\standardunitlength}
  \begin{array}{c}
    {\begingroup\makeatletter\ifx\SetFigFont\undefined%
\gdef\SetFigFont#1#2#3#4#5{%
  \reset@font\fontsize{#1}{#2pt}%
  \fontfamily{#3}\fontseries{#4}\fontshape{#5}%
  \selectfont}%
\fi\endgroup%
{\renewcommand{\dashlinestretch}{30}
\begin{picture}(3044,3059)(0,-10)
\thicklines
\path(922,22)(922,25)(922,33)
	(922,46)(922,66)(922,92)
	(923,124)(923,159)(923,198)
	(924,238)(924,277)(925,315)
	(925,351)(926,385)(927,416)
	(927,445)(928,472)(929,497)
	(930,521)(932,543)(933,564)
	(935,584)(936,607)(938,629)
	(941,650)(943,672)(946,694)
	(949,716)(952,738)(956,760)
	(960,782)(965,804)(969,826)
	(974,847)(980,868)(985,888)
	(991,907)(997,926)(1003,945)
	(1009,962)(1015,980)(1022,997)
	(1029,1014)(1037,1032)(1044,1049)
	(1053,1068)(1062,1086)(1072,1105)
	(1082,1125)(1092,1145)(1103,1166)
	(1115,1186)(1126,1207)(1138,1228)
	(1150,1248)(1163,1269)(1175,1289)
	(1187,1308)(1199,1328)(1211,1347)
	(1223,1366)(1235,1384)(1247,1404)
	(1259,1423)(1272,1442)(1285,1462)
	(1299,1482)(1313,1502)(1328,1522)
	(1343,1541)(1358,1560)(1374,1578)
	(1389,1595)(1405,1610)(1421,1624)
	(1436,1637)(1451,1647)(1466,1656)
	(1480,1663)(1494,1668)(1508,1671)
	(1522,1672)(1536,1671)(1550,1668)
	(1564,1663)(1578,1656)(1593,1647)
	(1608,1637)(1623,1624)(1639,1610)
	(1655,1595)(1670,1578)(1686,1560)
	(1701,1541)(1716,1522)(1731,1502)
	(1745,1482)(1759,1462)(1772,1442)
	(1785,1423)(1797,1404)(1810,1384)
	(1821,1366)(1833,1347)(1845,1328)
	(1857,1308)(1869,1289)(1881,1269)
	(1894,1248)(1906,1228)(1918,1207)
	(1929,1186)(1941,1166)(1952,1145)
	(1962,1125)(1972,1105)(1982,1086)
	(1991,1068)(2000,1049)(2007,1032)
	(2015,1014)(2022,997)(2029,980)
	(2035,962)(2041,945)(2047,926)
	(2053,907)(2059,888)(2064,868)
	(2070,847)(2075,826)(2079,804)
	(2084,782)(2088,760)(2092,738)
	(2095,716)(2098,694)(2101,672)
	(2103,650)(2106,629)(2108,607)
	(2110,584)(2111,564)(2112,543)
	(2114,521)(2115,497)(2116,472)
	(2117,445)(2117,416)(2118,385)
	(2119,351)(2119,315)(2120,277)
	(2120,238)(2121,198)(2121,159)
	(2121,124)(2122,92)(2122,66)
	(2122,46)(2122,33)(2122,25)(2122,22)
\path(922,3022)(923,3019)(924,3014)
	(926,3003)(930,2987)(935,2965)
	(941,2937)(948,2905)(957,2868)
	(966,2828)(976,2787)(986,2745)
	(996,2704)(1006,2663)(1016,2625)
	(1026,2588)(1036,2554)(1045,2521)
	(1055,2491)(1064,2463)(1073,2437)
	(1082,2412)(1092,2388)(1102,2365)
	(1112,2343)(1122,2322)(1134,2299)
	(1146,2277)(1159,2255)(1172,2233)
	(1187,2211)(1201,2190)(1217,2169)
	(1233,2149)(1250,2129)(1268,2110)
	(1286,2092)(1304,2075)(1322,2059)
	(1341,2044)(1360,2030)(1379,2018)
	(1397,2007)(1416,1998)(1434,1990)
	(1452,1983)(1470,1978)(1487,1975)
	(1505,1973)(1522,1972)(1539,1973)
	(1557,1975)(1574,1978)(1592,1983)
	(1610,1990)(1628,1998)(1647,2007)
	(1665,2018)(1684,2030)(1703,2044)
	(1722,2059)(1740,2075)(1758,2092)
	(1776,2110)(1794,2129)(1811,2149)
	(1827,2169)(1843,2190)(1857,2211)
	(1872,2233)(1885,2255)(1898,2277)
	(1910,2299)(1922,2322)(1932,2343)
	(1942,2365)(1952,2388)(1962,2412)
	(1971,2437)(1980,2463)(1989,2491)
	(1999,2521)(2008,2554)(2018,2588)
	(2028,2625)(2038,2663)(2048,2704)
	(2058,2745)(2068,2787)(2078,2828)
	(2087,2868)(2096,2905)(2103,2937)
	(2109,2965)(2114,2987)(2118,3003)
	(2120,3014)(2121,3019)(2122,3022)
\path(22,1522)(25,1520)(31,1517)
	(43,1510)(60,1500)(83,1488)
	(110,1472)(141,1455)(174,1436)
	(208,1417)(242,1398)(275,1379)
	(306,1362)(334,1346)(361,1332)
	(385,1318)(407,1306)(428,1295)
	(447,1285)(464,1276)(481,1268)
	(497,1259)(518,1249)(538,1239)
	(558,1230)(578,1222)(598,1213)
	(620,1204)(642,1195)(665,1186)
	(689,1177)(712,1169)(732,1161)
	(749,1155)(762,1151)(769,1148)(772,1147)
\path(1297,997)(1747,997)
\path(2272,1147)(2275,1148)(2282,1150)
	(2293,1153)(2310,1158)(2331,1164)
	(2355,1171)(2381,1178)(2407,1186)
	(2433,1194)(2457,1202)(2480,1210)
	(2501,1218)(2521,1225)(2541,1233)
	(2559,1242)(2578,1250)(2597,1259)
	(2613,1268)(2629,1276)(2646,1285)
	(2664,1295)(2683,1306)(2703,1318)
	(2725,1332)(2749,1346)(2774,1362)
	(2801,1379)(2830,1398)(2860,1417)
	(2890,1436)(2919,1455)(2946,1472)
	(2969,1488)(2989,1500)(3004,1510)
	(3014,1517)(3019,1520)(3022,1522)
\end{picture}
} } 
  \end{array}
$
and $
  \setlength{\unitlength}{0.12\standardunitlength}
  \begin{array}{c}
    {\begingroup\makeatletter\ifx\SetFigFont\undefined%
\gdef\SetFigFont#1#2#3#4#5{%
  \reset@font\fontsize{#1}{#2pt}%
  \fontfamily{#3}\fontseries{#4}\fontshape{#5}%
  \selectfont}%
\fi\endgroup%
{\renewcommand{\dashlinestretch}{30}
\begin{picture}(3044,3059)(0,-10)
\thicklines
\path(922,3022)(922,3019)(922,3011)
	(922,2998)(922,2978)(922,2952)
	(923,2920)(923,2885)(923,2846)
	(924,2806)(924,2767)(925,2729)
	(925,2693)(926,2659)(927,2628)
	(927,2599)(928,2572)(929,2547)
	(930,2523)(932,2501)(933,2480)
	(935,2459)(936,2437)(938,2415)
	(941,2394)(943,2372)(946,2350)
	(949,2328)(952,2306)(956,2284)
	(960,2262)(965,2240)(969,2218)
	(974,2197)(980,2176)(985,2156)
	(991,2137)(997,2118)(1003,2099)
	(1009,2082)(1015,2064)(1022,2047)
	(1029,2030)(1037,2012)(1044,1995)
	(1053,1976)(1062,1958)(1072,1939)
	(1082,1919)(1092,1899)(1103,1878)
	(1115,1858)(1126,1837)(1138,1816)
	(1150,1796)(1163,1775)(1175,1755)
	(1187,1736)(1199,1716)(1211,1697)
	(1223,1678)(1235,1659)(1247,1640)
	(1259,1621)(1272,1602)(1285,1582)
	(1299,1562)(1313,1542)(1328,1522)
	(1343,1503)(1358,1484)(1374,1466)
	(1389,1449)(1405,1434)(1421,1420)
	(1436,1407)(1451,1397)(1466,1388)
	(1480,1381)(1494,1376)(1508,1373)
	(1522,1372)(1536,1373)(1550,1376)
	(1564,1381)(1578,1388)(1593,1397)
	(1608,1407)(1623,1420)(1639,1434)
	(1655,1449)(1670,1466)(1686,1484)
	(1701,1503)(1716,1522)(1731,1542)
	(1745,1562)(1759,1582)(1772,1602)
	(1785,1621)(1797,1640)(1810,1659)
	(1821,1678)(1833,1697)(1845,1716)
	(1857,1736)(1869,1755)(1881,1775)
	(1894,1796)(1906,1816)(1918,1837)
	(1929,1858)(1941,1878)(1952,1899)
	(1962,1919)(1972,1939)(1982,1958)
	(1991,1976)(2000,1995)(2007,2012)
	(2015,2030)(2022,2047)(2029,2064)
	(2035,2082)(2041,2099)(2047,2118)
	(2053,2137)(2059,2156)(2064,2176)
	(2070,2197)(2075,2218)(2079,2240)
	(2084,2262)(2088,2284)(2092,2306)
	(2095,2328)(2098,2350)(2101,2372)
	(2103,2394)(2106,2415)(2108,2437)
	(2110,2459)(2111,2480)(2112,2501)
	(2114,2523)(2115,2547)(2116,2572)
	(2117,2599)(2117,2628)(2118,2659)
	(2119,2693)(2119,2729)(2120,2767)
	(2120,2806)(2121,2846)(2121,2885)
	(2121,2920)(2122,2952)(2122,2978)
	(2122,2998)(2122,3011)(2122,3019)(2122,3022)
\path(922,22)(923,25)(924,30)
	(926,41)(930,57)(935,79)
	(941,107)(948,139)(957,176)
	(966,216)(976,257)(986,299)
	(996,340)(1006,381)(1016,419)
	(1026,456)(1036,490)(1045,523)
	(1055,553)(1064,581)(1073,607)
	(1082,632)(1092,656)(1102,679)
	(1112,701)(1122,722)(1134,745)
	(1146,767)(1159,789)(1172,811)
	(1187,833)(1201,854)(1217,875)
	(1233,895)(1250,915)(1268,934)
	(1286,952)(1304,969)(1322,985)
	(1341,1000)(1360,1014)(1379,1026)
	(1397,1037)(1416,1046)(1434,1054)
	(1452,1061)(1470,1066)(1487,1069)
	(1505,1071)(1522,1072)(1539,1071)
	(1557,1069)(1574,1066)(1592,1061)
	(1610,1054)(1628,1046)(1647,1037)
	(1665,1026)(1684,1014)(1703,1000)
	(1722,985)(1740,969)(1758,952)
	(1776,934)(1794,915)(1811,895)
	(1827,875)(1843,854)(1857,833)
	(1872,811)(1885,789)(1898,767)
	(1910,745)(1922,722)(1932,701)
	(1942,679)(1952,656)(1962,632)
	(1971,607)(1980,581)(1989,553)
	(1999,523)(2008,490)(2018,456)
	(2028,419)(2038,381)(2048,340)
	(2058,299)(2068,257)(2078,216)
	(2087,176)(2096,139)(2103,107)
	(2109,79)(2114,57)(2118,41)
	(2120,30)(2121,25)(2122,22)
\path(22,1522)(25,1524)(31,1527)
	(43,1534)(60,1544)(83,1556)
	(110,1572)(141,1589)(174,1608)
	(208,1627)(242,1646)(275,1665)
	(306,1682)(334,1698)(361,1712)
	(385,1726)(407,1738)(428,1749)
	(447,1759)(464,1768)(481,1776)
	(497,1784)(518,1795)(538,1805)
	(558,1814)(578,1822)(598,1831)
	(620,1840)(642,1849)(665,1858)
	(689,1867)(712,1875)(732,1883)
	(749,1889)(762,1893)(769,1896)(772,1897)
\path(1297,2047)(1747,2047)
\path(2272,1897)(2275,1896)(2282,1894)
	(2293,1891)(2310,1886)(2331,1880)
	(2355,1873)(2381,1866)(2407,1858)
	(2433,1850)(2457,1842)(2480,1834)
	(2501,1826)(2521,1819)(2541,1811)
	(2559,1802)(2578,1794)(2597,1784)
	(2613,1776)(2629,1768)(2646,1759)
	(2664,1749)(2683,1738)(2703,1726)
	(2725,1712)(2749,1698)(2774,1682)
	(2801,1665)(2830,1646)(2860,1627)
	(2890,1608)(2919,1589)(2946,1572)
	(2969,1556)(2989,1544)(3004,1534)
	(3014,1527)(3019,1524)(3022,1522)
\end{picture}
} } 
  \end{array}
$ and it is tempting to use (R2) on both to
reduce to

\[
  \def\neg{\hspace{-2mm}}
  \def\D#1{\neg
  \setlength{\unitlength}{0.125\standardunitlength}
  \begin{array}{c}
    {\input figs/D#1.tex } 
  \end{array}
\neg}
  \def\U#1{\neg
  \setlength{\unitlength}{0.125\standardunitlength}
  \begin{array}{c}
    {\input figs/U#1.tex } 
  \end{array}
\neg}
  \begin{array}{c} \xymatrix@R=4mm@C=4mm{
    & 0 \ar@{->}[rr] \ar@{<-}'[d][dd]
      & & 0 \ar@{<-}[dd] \\
    0 \ar@{->}[ur] \ar@{->}[rr] \ar@{<-}[dd] &
      & \D{110} \ar@{->}[ur] \ar@{<-}[dd] \\
    & \D{001} \ar@{->}'[r][rr] & & \D{011} \\
    \D{000} \ar@{->}[rr] \ar@{->}[ur] & & \D{010} \ar@{->}[ur]
  } \end{array}
  \qquad
  \begin{array}{c} \xymatrix@R=4mm@C=4mm{
    & \U{101} \ar@{->}[rr] \ar@{<-}'[d][dd]
      & & 0 \ar@{<-}[dd] \\
    0 \ar@{->}[ur] \ar@{->}[rr] \ar@{<-}[dd] &
      & 0 \ar@{->}[ur] \ar@{<-}[dd] \\
    & \U{001} \ar@{->}'[r][rr] & & \U{011} \\
    \U{000} \ar@{->}[rr] \ar@{->}[ur] & & \U{010} \ar@{->}[ur]
  } \end{array}
\]

But this fails for two reasons. These cubes aren't isomorphic (their
bottom layers are isomorphic and their top layers are isomorphic, but the
maps between them are different), and the (R2)-style reduction used to get
here is invalid, for in the presence of the bottom layers what would be the
analog of $\calC''$ simply isn't a subcomplex. Fortunately, there is a
somewhat more complicated proof of invariance under (R2) that does lead to
a correct argument for invariance under (R3).

\subsubsection{Invariance under (R2), second proof.} We start
in the same way as in the first proof and reduce to the complex
$\calC/\calC'$ which is displayed once again in Figure~2 (except
this time we suppress the $\llbracket\cdot\rrbracket$ brackets
and the degree shift $\{\cdot\}$ symbols). In $\calC/\calC'$
the vertical arrow $\Delta$ is a bijection so we can invert it
and compose with the horizontal arrow $d_{\star 0}$ to get a map
$\tau:
  \setlength{\unitlength}{0.2\standardunitlength}
  \begin{array}{c}
    {\begingroup\makeatletter\ifx\SetFigFont\undefined%
\gdef\SetFigFont#1#2#3#4#5{%
  \reset@font\fontsize{#1}{#2pt}%
  \fontfamily{#3}\fontseries{#4}\fontshape{#5}%
  \selectfont}%
\fi\endgroup%
{\renewcommand{\dashlinestretch}{30}
\begin{picture}(2144,659)(0,-10)
\thicklines
\path(22,22)(23,22)(26,22)
	(35,23)(51,23)(72,24)
	(97,26)(125,28)(153,31)
	(180,34)(205,37)(228,41)
	(249,46)(269,51)(287,57)
	(305,64)(322,72)(337,80)
	(353,89)(368,99)(384,110)
	(399,122)(415,135)(430,149)
	(444,164)(458,179)(471,195)
	(483,211)(493,227)(502,243)
	(509,260)(515,276)(519,291)
	(521,307)(522,322)(521,337)
	(519,353)(515,368)(509,384)
	(502,401)(493,417)(483,433)
	(471,449)(458,465)(444,480)
	(430,495)(415,509)(399,522)
	(384,534)(368,545)(353,555)
	(337,564)(322,572)(305,580)
	(287,587)(269,593)(249,598)
	(228,603)(205,607)(180,610)
	(153,613)(125,616)(97,618)
	(72,620)(51,621)(35,621)
	(26,622)(23,622)(22,622)
\path(722,322)(723,307)(725,291)
	(729,276)(735,260)(742,243)
	(751,227)(761,211)(773,195)
	(786,179)(800,164)(814,149)
	(829,135)(845,122)(860,110)
	(876,99)(891,89)(907,80)
	(922,72)(937,65)(953,58)
	(969,53)(985,48)(1002,43)
	(1019,40)(1036,38)(1054,36)
	(1072,36)(1090,36)(1108,38)
	(1125,40)(1142,43)(1159,48)
	(1175,53)(1191,58)(1207,65)
	(1222,72)(1237,80)(1253,89)
	(1268,99)(1284,110)(1299,122)
	(1315,135)(1330,149)(1344,164)
	(1358,179)(1371,195)(1383,211)
	(1393,227)(1402,243)(1409,260)
	(1415,276)(1419,291)(1421,307)
	(1422,322)(1421,337)(1419,353)
	(1415,368)(1409,384)(1402,401)
	(1393,417)(1383,433)(1371,449)
	(1358,465)(1344,480)(1330,495)
	(1315,509)(1299,522)(1284,534)
	(1268,545)(1253,555)(1237,564)
	(1222,572)(1207,579)(1191,586)
	(1175,591)(1159,596)(1142,601)
	(1125,604)(1108,606)(1090,608)
	(1072,608)(1054,608)(1036,606)
	(1019,604)(1002,601)(985,596)
	(969,591)(953,586)(937,579)
	(922,572)(907,564)(891,555)
	(876,545)(860,534)(845,522)
	(829,509)(814,495)(800,480)
	(786,465)(773,449)(761,433)
	(751,417)(742,401)(735,384)
	(729,368)(725,353)(723,337)(722,322)
\path(2122,622)(2121,622)(2118,622)
	(2109,621)(2093,621)(2072,620)
	(2047,618)(2019,616)(1991,613)
	(1964,610)(1939,607)(1916,603)
	(1895,598)(1875,593)(1857,587)
	(1839,580)(1822,572)(1807,564)
	(1791,555)(1776,545)(1760,534)
	(1745,522)(1729,509)(1714,495)
	(1700,480)(1686,465)(1673,449)
	(1661,433)(1651,417)(1642,401)
	(1635,384)(1629,368)(1625,353)
	(1623,337)(1622,322)(1623,307)
	(1625,291)(1629,276)(1635,260)
	(1642,243)(1651,227)(1661,211)
	(1673,195)(1686,179)(1700,164)
	(1714,149)(1729,135)(1745,122)
	(1760,110)(1776,99)(1791,89)
	(1807,80)(1822,72)(1839,64)
	(1857,57)(1875,51)(1895,46)
	(1916,41)(1939,37)(1964,34)
	(1991,31)(2019,28)(2047,26)
	(2072,24)(2093,23)(2109,23)
	(2118,22)(2121,22)(2122,22)
\end{picture}
} } 
  \end{array}
\!\!_{/v_+=0}\to
  \setlength{\unitlength}{0.2\standardunitlength}
  \begin{array}{c}
    {\begingroup\makeatletter\ifx\SetFigFont\undefined%
\gdef\SetFigFont#1#2#3#4#5{%
  \reset@font\fontsize{#1}{#2pt}%
  \fontfamily{#3}\fontseries{#4}\fontshape{#5}%
  \selectfont}%
\fi\endgroup%
{\renewcommand{\dashlinestretch}{30}
\begin{picture}(2144,659)(0,-10)
\thicklines
\path(22,22)(23,22)(26,22)
	(35,23)(51,23)(72,24)
	(97,26)(125,28)(153,31)
	(180,34)(205,37)(228,41)
	(249,46)(269,51)(287,57)
	(305,64)(322,72)(337,80)
	(353,89)(369,99)(385,109)
	(402,120)(419,131)(436,142)
	(454,154)(472,165)(490,176)
	(508,186)(525,195)(542,203)
	(559,210)(575,215)(591,219)
	(607,221)(622,222)(637,221)
	(653,219)(669,215)(685,210)
	(702,203)(719,195)(736,186)
	(754,176)(772,165)(790,154)
	(808,142)(825,131)(842,120)
	(859,109)(875,99)(891,89)
	(907,80)(922,72)(937,65)
	(953,58)(969,53)(985,48)
	(1002,43)(1019,40)(1036,38)
	(1054,36)(1072,36)(1090,36)
	(1108,38)(1125,40)(1142,43)
	(1159,48)(1175,53)(1191,58)
	(1207,65)(1222,72)(1237,80)
	(1253,89)(1269,99)(1285,109)
	(1302,120)(1319,131)(1336,142)
	(1354,154)(1372,165)(1390,176)
	(1408,186)(1425,195)(1442,203)
	(1459,210)(1475,215)(1491,219)
	(1507,221)(1522,222)(1537,221)
	(1553,219)(1569,215)(1585,210)
	(1602,203)(1619,195)(1636,186)
	(1654,176)(1672,165)(1690,154)
	(1708,142)(1725,131)(1742,120)
	(1759,109)(1775,99)(1791,89)
	(1807,80)(1822,72)(1839,64)
	(1857,57)(1875,51)(1895,46)
	(1916,41)(1939,37)(1964,34)
	(1991,31)(2019,28)(2047,26)
	(2072,24)(2093,23)(2109,23)
	(2118,22)(2121,22)(2122,22)
\path(2122,622)(2121,622)(2118,622)
	(2109,621)(2093,621)(2072,620)
	(2047,618)(2019,616)(1991,613)
	(1964,610)(1939,607)(1916,603)
	(1895,598)(1875,593)(1857,587)
	(1839,580)(1822,572)(1807,564)
	(1791,555)(1775,545)(1759,535)
	(1742,524)(1725,513)(1708,502)
	(1690,490)(1672,479)(1654,468)
	(1636,458)(1619,449)(1602,441)
	(1585,434)(1569,429)(1553,425)
	(1537,423)(1522,422)(1507,423)
	(1491,425)(1475,429)(1459,434)
	(1442,441)(1425,449)(1408,458)
	(1390,468)(1372,479)(1354,490)
	(1336,502)(1319,513)(1302,524)
	(1285,535)(1269,545)(1253,555)
	(1237,564)(1222,572)(1207,579)
	(1191,586)(1175,591)(1159,596)
	(1142,601)(1125,604)(1108,606)
	(1090,608)(1072,608)(1054,608)
	(1036,606)(1019,604)(1002,601)
	(985,596)(969,591)(953,586)
	(937,579)(922,572)(907,564)
	(891,555)(875,545)(859,535)
	(842,524)(825,513)(808,502)
	(790,490)(772,479)(754,468)
	(736,458)(719,449)(702,441)
	(685,434)(669,429)(653,425)
	(637,423)(622,422)(607,423)
	(591,425)(575,429)(559,434)
	(542,441)(525,449)(508,458)
	(490,468)(472,479)(454,490)
	(436,502)(419,513)(402,524)
	(385,535)(369,545)(353,555)
	(337,564)(322,572)(305,580)
	(287,587)(269,593)(249,598)
	(228,603)(205,607)(180,610)
	(153,613)(125,616)(97,618)
	(72,620)(51,621)(35,621)
	(26,622)(23,622)(22,622)
\end{picture}
} } 
  \end{array}
$. We
now let $\calC'''$ be the subcomplex of $\calC/\calC'$ containing
all $\alpha\in
  \setlength{\unitlength}{0.2\standardunitlength}
  \begin{array}{c}
    {\begingroup\makeatletter\ifx\SetFigFont\undefined%
\gdef\SetFigFont#1#2#3#4#5{%
  \reset@font\fontsize{#1}{#2pt}%
  \fontfamily{#3}\fontseries{#4}\fontshape{#5}%
  \selectfont}%
\fi\endgroup%
{\renewcommand{\dashlinestretch}{30}
\begin{picture}(2144,659)(0,-10)
\thicklines
\path(22,22)(23,22)(26,22)
	(35,23)(51,23)(72,24)
	(97,26)(125,28)(153,31)
	(180,34)(205,37)(228,41)
	(249,46)(269,51)(287,57)
	(305,64)(322,72)(337,80)
	(353,89)(368,99)(384,110)
	(399,122)(415,135)(430,149)
	(444,164)(458,179)(471,195)
	(483,211)(493,227)(502,243)
	(509,260)(515,276)(519,291)
	(521,307)(522,322)(521,337)
	(519,353)(515,368)(509,384)
	(502,401)(493,417)(483,433)
	(471,449)(458,465)(444,480)
	(430,495)(415,509)(399,522)
	(384,534)(368,545)(353,555)
	(337,564)(322,572)(305,580)
	(287,587)(269,593)(249,598)
	(228,603)(205,607)(180,610)
	(153,613)(125,616)(97,618)
	(72,620)(51,621)(35,621)
	(26,622)(23,622)(22,622)
\path(2122,622)(2121,622)(2118,622)
	(2109,621)(2093,621)(2072,620)
	(2047,618)(2019,616)(1991,613)
	(1964,610)(1939,607)(1916,603)
	(1895,598)(1875,593)(1857,587)
	(1839,580)(1822,572)(1807,564)
	(1791,555)(1775,545)(1759,535)
	(1742,524)(1725,513)(1708,502)
	(1690,490)(1672,479)(1654,468)
	(1636,458)(1619,449)(1602,441)
	(1585,434)(1569,429)(1553,425)
	(1537,423)(1522,422)(1507,423)
	(1491,425)(1475,429)(1459,434)
	(1442,441)(1425,449)(1408,458)
	(1390,468)(1372,479)(1354,490)
	(1336,502)(1319,513)(1302,524)
	(1285,535)(1269,545)(1253,555)
	(1237,564)(1222,572)(1207,579)
	(1191,586)(1175,591)(1159,596)
	(1142,601)(1125,604)(1108,606)
	(1090,608)(1072,608)(1054,608)
	(1036,606)(1019,604)(1002,601)
	(985,596)(969,591)(953,586)
	(937,579)(922,572)(907,564)
	(891,555)(876,545)(860,534)
	(845,522)(829,509)(814,495)
	(800,480)(786,465)(773,449)
	(761,433)(751,417)(742,401)
	(735,384)(729,368)(725,353)
	(723,337)(722,322)(723,307)
	(725,291)(729,276)(735,260)
	(742,243)(751,227)(761,211)
	(773,195)(786,179)(800,164)
	(814,149)(829,135)(845,122)
	(860,110)(876,99)(891,89)
	(907,80)(922,72)(937,65)
	(953,58)(969,53)(985,48)
	(1002,43)(1019,40)(1036,38)
	(1054,36)(1072,36)(1090,36)
	(1108,38)(1125,40)(1142,43)
	(1159,48)(1175,53)(1191,58)
	(1207,65)(1222,72)(1237,80)
	(1253,89)(1269,99)(1285,109)
	(1302,120)(1319,131)(1336,142)
	(1354,154)(1372,165)(1390,176)
	(1408,186)(1425,195)(1442,203)
	(1459,210)(1475,215)(1491,219)
	(1507,221)(1522,222)(1537,221)
	(1553,219)(1569,215)(1585,210)
	(1602,203)(1619,195)(1636,186)
	(1654,176)(1672,165)(1690,154)
	(1708,142)(1725,131)(1742,120)
	(1759,109)(1775,99)(1791,89)
	(1807,80)(1822,72)(1839,64)
	(1857,57)(1875,51)(1895,46)
	(1916,41)(1939,37)(1964,34)
	(1991,31)(2019,28)(2047,26)
	(2072,24)(2093,23)(2109,23)
	(2118,22)(2121,22)(2122,22)
\end{picture}
} } 
  \end{array}
$ and all pairs of the form
$(\beta,\tau\beta) \in 
  \setlength{\unitlength}{0.2\standardunitlength}
  \begin{array}{c}
    {\begingroup\makeatletter\ifx\SetFigFont\undefined%
\gdef\SetFigFont#1#2#3#4#5{%
  \reset@font\fontsize{#1}{#2pt}%
  \fontfamily{#3}\fontseries{#4}\fontshape{#5}%
  \selectfont}%
\fi\endgroup%
{\renewcommand{\dashlinestretch}{30}
\begin{picture}(2144,659)(0,-10)
\thicklines
\path(22,22)(23,22)(26,22)
	(35,23)(51,23)(72,24)
	(97,26)(125,28)(153,31)
	(180,34)(205,37)(228,41)
	(249,46)(269,51)(287,57)
	(305,64)(322,72)(337,80)
	(353,89)(368,99)(384,110)
	(399,122)(415,135)(430,149)
	(444,164)(458,179)(471,195)
	(483,211)(493,227)(502,243)
	(509,260)(515,276)(519,291)
	(521,307)(522,322)(521,337)
	(519,353)(515,368)(509,384)
	(502,401)(493,417)(483,433)
	(471,449)(458,465)(444,480)
	(430,495)(415,509)(399,522)
	(384,534)(368,545)(353,555)
	(337,564)(322,572)(305,580)
	(287,587)(269,593)(249,598)
	(228,603)(205,607)(180,610)
	(153,613)(125,616)(97,618)
	(72,620)(51,621)(35,621)
	(26,622)(23,622)(22,622)
\path(722,322)(723,307)(725,291)
	(729,276)(735,260)(742,243)
	(751,227)(761,211)(773,195)
	(786,179)(800,164)(814,149)
	(829,135)(845,122)(860,110)
	(876,99)(891,89)(907,80)
	(922,72)(937,65)(953,58)
	(969,53)(985,48)(1002,43)
	(1019,40)(1036,38)(1054,36)
	(1072,36)(1090,36)(1108,38)
	(1125,40)(1142,43)(1159,48)
	(1175,53)(1191,58)(1207,65)
	(1222,72)(1237,80)(1253,89)
	(1268,99)(1284,110)(1299,122)
	(1315,135)(1330,149)(1344,164)
	(1358,179)(1371,195)(1383,211)
	(1393,227)(1402,243)(1409,260)
	(1415,276)(1419,291)(1421,307)
	(1422,322)(1421,337)(1419,353)
	(1415,368)(1409,384)(1402,401)
	(1393,417)(1383,433)(1371,449)
	(1358,465)(1344,480)(1330,495)
	(1315,509)(1299,522)(1284,534)
	(1268,545)(1253,555)(1237,564)
	(1222,572)(1207,579)(1191,586)
	(1175,591)(1159,596)(1142,601)
	(1125,604)(1108,606)(1090,608)
	(1072,608)(1054,608)(1036,606)
	(1019,604)(1002,601)(985,596)
	(969,591)(953,586)(937,579)
	(922,572)(907,564)(891,555)
	(876,545)(860,534)(845,522)
	(829,509)(814,495)(800,480)
	(786,465)(773,449)(761,433)
	(751,417)(742,401)(735,384)
	(729,368)(725,353)(723,337)(722,322)
\path(2122,622)(2121,622)(2118,622)
	(2109,621)(2093,621)(2072,620)
	(2047,618)(2019,616)(1991,613)
	(1964,610)(1939,607)(1916,603)
	(1895,598)(1875,593)(1857,587)
	(1839,580)(1822,572)(1807,564)
	(1791,555)(1776,545)(1760,534)
	(1745,522)(1729,509)(1714,495)
	(1700,480)(1686,465)(1673,449)
	(1661,433)(1651,417)(1642,401)
	(1635,384)(1629,368)(1625,353)
	(1623,337)(1622,322)(1623,307)
	(1625,291)(1629,276)(1635,260)
	(1642,243)(1651,227)(1661,211)
	(1673,195)(1686,179)(1700,164)
	(1714,149)(1729,135)(1745,122)
	(1760,110)(1776,99)(1791,89)
	(1807,80)(1822,72)(1839,64)
	(1857,57)(1875,51)(1895,46)
	(1916,41)(1939,37)(1964,34)
	(1991,31)(2019,28)(2047,26)
	(2072,24)(2093,23)(2109,23)
	(2118,22)(2121,22)(2122,22)
\end{picture}
} } 
  \end{array}
\!\!_{/v_+=0} \oplus

  \setlength{\unitlength}{0.2\standardunitlength}
  \begin{array}{c}
    {\begingroup\makeatletter\ifx\SetFigFont\undefined%
\gdef\SetFigFont#1#2#3#4#5{%
  \reset@font\fontsize{#1}{#2pt}%
  \fontfamily{#3}\fontseries{#4}\fontshape{#5}%
  \selectfont}%
\fi\endgroup%
{\renewcommand{\dashlinestretch}{30}
\begin{picture}(2144,659)(0,-10)
\thicklines
\path(22,22)(23,22)(26,22)
	(35,23)(51,23)(72,24)
	(97,26)(125,28)(153,31)
	(180,34)(205,37)(228,41)
	(249,46)(269,51)(287,57)
	(305,64)(322,72)(337,80)
	(353,89)(369,99)(385,109)
	(402,120)(419,131)(436,142)
	(454,154)(472,165)(490,176)
	(508,186)(525,195)(542,203)
	(559,210)(575,215)(591,219)
	(607,221)(622,222)(637,221)
	(653,219)(669,215)(685,210)
	(702,203)(719,195)(736,186)
	(754,176)(772,165)(790,154)
	(808,142)(825,131)(842,120)
	(859,109)(875,99)(891,89)
	(907,80)(922,72)(937,65)
	(953,58)(969,53)(985,48)
	(1002,43)(1019,40)(1036,38)
	(1054,36)(1072,36)(1090,36)
	(1108,38)(1125,40)(1142,43)
	(1159,48)(1175,53)(1191,58)
	(1207,65)(1222,72)(1237,80)
	(1253,89)(1269,99)(1285,109)
	(1302,120)(1319,131)(1336,142)
	(1354,154)(1372,165)(1390,176)
	(1408,186)(1425,195)(1442,203)
	(1459,210)(1475,215)(1491,219)
	(1507,221)(1522,222)(1537,221)
	(1553,219)(1569,215)(1585,210)
	(1602,203)(1619,195)(1636,186)
	(1654,176)(1672,165)(1690,154)
	(1708,142)(1725,131)(1742,120)
	(1759,109)(1775,99)(1791,89)
	(1807,80)(1822,72)(1839,64)
	(1857,57)(1875,51)(1895,46)
	(1916,41)(1939,37)(1964,34)
	(1991,31)(2019,28)(2047,26)
	(2072,24)(2093,23)(2109,23)
	(2118,22)(2121,22)(2122,22)
\path(2122,622)(2121,622)(2118,622)
	(2109,621)(2093,621)(2072,620)
	(2047,618)(2019,616)(1991,613)
	(1964,610)(1939,607)(1916,603)
	(1895,598)(1875,593)(1857,587)
	(1839,580)(1822,572)(1807,564)
	(1791,555)(1775,545)(1759,535)
	(1742,524)(1725,513)(1708,502)
	(1690,490)(1672,479)(1654,468)
	(1636,458)(1619,449)(1602,441)
	(1585,434)(1569,429)(1553,425)
	(1537,423)(1522,422)(1507,423)
	(1491,425)(1475,429)(1459,434)
	(1442,441)(1425,449)(1408,458)
	(1390,468)(1372,479)(1354,490)
	(1336,502)(1319,513)(1302,524)
	(1285,535)(1269,545)(1253,555)
	(1237,564)(1222,572)(1207,579)
	(1191,586)(1175,591)(1159,596)
	(1142,601)(1125,604)(1108,606)
	(1090,608)(1072,608)(1054,608)
	(1036,606)(1019,604)(1002,601)
	(985,596)(969,591)(953,586)
	(937,579)(922,572)(907,564)
	(891,555)(875,545)(859,535)
	(842,524)(825,513)(808,502)
	(790,490)(772,479)(754,468)
	(736,458)(719,449)(702,441)
	(685,434)(669,429)(653,425)
	(637,423)(622,422)(607,423)
	(591,425)(575,429)(559,434)
	(542,441)(525,449)(508,458)
	(490,468)(472,479)(454,490)
	(436,502)(419,513)(402,524)
	(385,535)(369,545)(353,555)
	(337,564)(322,572)(305,580)
	(287,587)(269,593)(249,598)
	(228,603)(205,607)(180,610)
	(153,613)(125,616)(97,618)
	(72,620)(51,621)(35,621)
	(26,622)(23,622)(22,622)
\end{picture}
} } 
  \end{array}
$ (see Figure~2). The map $\Delta$ is bijective
in $\calC'''$ and hence $\calC'''$ is acyclic and thus it is enough to
study $(\calC/\calC')/\calC'''$.

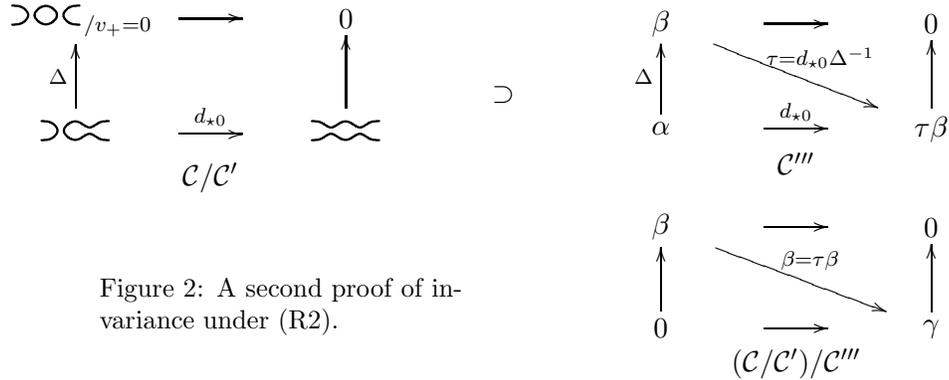
\begin{figure}
\def\entry#1{{\parbox{1in}{\centering $#1$}}}
\begin{eqnarray*}
  \begin{array}{c}
    \xymatrix{
      \entry{
  \setlength{\unitlength}{0.2\standardunitlength}
  \begin{array}{c}
    {\begingroup\makeatletter\ifx\SetFigFont\undefined%
\gdef\SetFigFont#1#2#3#4#5{%
  \reset@font\fontsize{#1}{#2pt}%
  \fontfamily{#3}\fontseries{#4}\fontshape{#5}%
  \selectfont}%
\fi\endgroup%
{\renewcommand{\dashlinestretch}{30}
\begin{picture}(2144,659)(0,-10)
\thicklines
\path(22,22)(23,22)(26,22)
	(35,23)(51,23)(72,24)
	(97,26)(125,28)(153,31)
	(180,34)(205,37)(228,41)
	(249,46)(269,51)(287,57)
	(305,64)(322,72)(337,80)
	(353,89)(368,99)(384,110)
	(399,122)(415,135)(430,149)
	(444,164)(458,179)(471,195)
	(483,211)(493,227)(502,243)
	(509,260)(515,276)(519,291)
	(521,307)(522,322)(521,337)
	(519,353)(515,368)(509,384)
	(502,401)(493,417)(483,433)
	(471,449)(458,465)(444,480)
	(430,495)(415,509)(399,522)
	(384,534)(368,545)(353,555)
	(337,564)(322,572)(305,580)
	(287,587)(269,593)(249,598)
	(228,603)(205,607)(180,610)
	(153,613)(125,616)(97,618)
	(72,620)(51,621)(35,621)
	(26,622)(23,622)(22,622)
\path(722,322)(723,307)(725,291)
	(729,276)(735,260)(742,243)
	(751,227)(761,211)(773,195)
	(786,179)(800,164)(814,149)
	(829,135)(845,122)(860,110)
	(876,99)(891,89)(907,80)
	(922,72)(937,65)(953,58)
	(969,53)(985,48)(1002,43)
	(1019,40)(1036,38)(1054,36)
	(1072,36)(1090,36)(1108,38)
	(1125,40)(1142,43)(1159,48)
	(1175,53)(1191,58)(1207,65)
	(1222,72)(1237,80)(1253,89)
	(1268,99)(1284,110)(1299,122)
	(1315,135)(1330,149)(1344,164)
	(1358,179)(1371,195)(1383,211)
	(1393,227)(1402,243)(1409,260)
	(1415,276)(1419,291)(1421,307)
	(1422,322)(1421,337)(1419,353)
	(1415,368)(1409,384)(1402,401)
	(1393,417)(1383,433)(1371,449)
	(1358,465)(1344,480)(1330,495)
	(1315,509)(1299,522)(1284,534)
	(1268,545)(1253,555)(1237,564)
	(1222,572)(1207,579)(1191,586)
	(1175,591)(1159,596)(1142,601)
	(1125,604)(1108,606)(1090,608)
	(1072,608)(1054,608)(1036,606)
	(1019,604)(1002,601)(985,596)
	(969,591)(953,586)(937,579)
	(922,572)(907,564)(891,555)
	(876,545)(860,534)(845,522)
	(829,509)(814,495)(800,480)
	(786,465)(773,449)(761,433)
	(751,417)(742,401)(735,384)
	(729,368)(725,353)(723,337)(722,322)
\path(2122,622)(2121,622)(2118,622)
	(2109,621)(2093,621)(2072,620)
	(2047,618)(2019,616)(1991,613)
	(1964,610)(1939,607)(1916,603)
	(1895,598)(1875,593)(1857,587)
	(1839,580)(1822,572)(1807,564)
	(1791,555)(1776,545)(1760,534)
	(1745,522)(1729,509)(1714,495)
	(1700,480)(1686,465)(1673,449)
	(1661,433)(1651,417)(1642,401)
	(1635,384)(1629,368)(1625,353)
	(1623,337)(1622,322)(1623,307)
	(1625,291)(1629,276)(1635,260)
	(1642,243)(1651,227)(1661,211)
	(1673,195)(1686,179)(1700,164)
	(1714,149)(1729,135)(1745,122)
	(1760,110)(1776,99)(1791,89)
	(1807,80)(1822,72)(1839,64)
	(1857,57)(1875,51)(1895,46)
	(1916,41)(1939,37)(1964,34)
	(1991,31)(2019,28)(2047,26)
	(2072,24)(2093,23)(2109,23)
	(2118,22)(2121,22)(2122,22)
\end{picture}
} } 
  \end{array}
\!\!_{/v_+=0}} \ar[r]
      & \entry{0} \\
      \entry{
  \setlength{\unitlength}{0.2\standardunitlength}
  \begin{array}{c}
    {\begingroup\makeatletter\ifx\SetFigFont\undefined%
\gdef\SetFigFont#1#2#3#4#5{%
  \reset@font\fontsize{#1}{#2pt}%
  \fontfamily{#3}\fontseries{#4}\fontshape{#5}%
  \selectfont}%
\fi\endgroup%
{\renewcommand{\dashlinestretch}{30}
\begin{picture}(2144,659)(0,-10)
\thicklines
\path(22,22)(23,22)(26,22)
	(35,23)(51,23)(72,24)
	(97,26)(125,28)(153,31)
	(180,34)(205,37)(228,41)
	(249,46)(269,51)(287,57)
	(305,64)(322,72)(337,80)
	(353,89)(368,99)(384,110)
	(399,122)(415,135)(430,149)
	(444,164)(458,179)(471,195)
	(483,211)(493,227)(502,243)
	(509,260)(515,276)(519,291)
	(521,307)(522,322)(521,337)
	(519,353)(515,368)(509,384)
	(502,401)(493,417)(483,433)
	(471,449)(458,465)(444,480)
	(430,495)(415,509)(399,522)
	(384,534)(368,545)(353,555)
	(337,564)(322,572)(305,580)
	(287,587)(269,593)(249,598)
	(228,603)(205,607)(180,610)
	(153,613)(125,616)(97,618)
	(72,620)(51,621)(35,621)
	(26,622)(23,622)(22,622)
\path(2122,622)(2121,622)(2118,622)
	(2109,621)(2093,621)(2072,620)
	(2047,618)(2019,616)(1991,613)
	(1964,610)(1939,607)(1916,603)
	(1895,598)(1875,593)(1857,587)
	(1839,580)(1822,572)(1807,564)
	(1791,555)(1775,545)(1759,535)
	(1742,524)(1725,513)(1708,502)
	(1690,490)(1672,479)(1654,468)
	(1636,458)(1619,449)(1602,441)
	(1585,434)(1569,429)(1553,425)
	(1537,423)(1522,422)(1507,423)
	(1491,425)(1475,429)(1459,434)
	(1442,441)(1425,449)(1408,458)
	(1390,468)(1372,479)(1354,490)
	(1336,502)(1319,513)(1302,524)
	(1285,535)(1269,545)(1253,555)
	(1237,564)(1222,572)(1207,579)
	(1191,586)(1175,591)(1159,596)
	(1142,601)(1125,604)(1108,606)
	(1090,608)(1072,608)(1054,608)
	(1036,606)(1019,604)(1002,601)
	(985,596)(969,591)(953,586)
	(937,579)(922,572)(907,564)
	(891,555)(876,545)(860,534)
	(845,522)(829,509)(814,495)
	(800,480)(786,465)(773,449)
	(761,433)(751,417)(742,401)
	(735,384)(729,368)(725,353)
	(723,337)(722,322)(723,307)
	(725,291)(729,276)(735,260)
	(742,243)(751,227)(761,211)
	(773,195)(786,179)(800,164)
	(814,149)(829,135)(845,122)
	(860,110)(876,99)(891,89)
	(907,80)(922,72)(937,65)
	(953,58)(969,53)(985,48)
	(1002,43)(1019,40)(1036,38)
	(1054,36)(1072,36)(1090,36)
	(1108,38)(1125,40)(1142,43)
	(1159,48)(1175,53)(1191,58)
	(1207,65)(1222,72)(1237,80)
	(1253,89)(1269,99)(1285,109)
	(1302,120)(1319,131)(1336,142)
	(1354,154)(1372,165)(1390,176)
	(1408,186)(1425,195)(1442,203)
	(1459,210)(1475,215)(1491,219)
	(1507,221)(1522,222)(1537,221)
	(1553,219)(1569,215)(1585,210)
	(1602,203)(1619,195)(1636,186)
	(1654,176)(1672,165)(1690,154)
	(1708,142)(1725,131)(1742,120)
	(1759,109)(1775,99)(1791,89)
	(1807,80)(1822,72)(1839,64)
	(1857,57)(1875,51)(1895,46)
	(1916,41)(1939,37)(1964,34)
	(1991,31)(2019,28)(2047,26)
	(2072,24)(2093,23)(2109,23)
	(2118,22)(2121,22)(2122,22)
\end{picture}
} } 
  \end{array}
} \ar[u]^\Delta \ar[r]^{d_{\star 0}}
      & \entry{
  \setlength{\unitlength}{0.2\standardunitlength}
  \begin{array}{c}
    {\begingroup\makeatletter\ifx\SetFigFont\undefined%
\gdef\SetFigFont#1#2#3#4#5{%
  \reset@font\fontsize{#1}{#2pt}%
  \fontfamily{#3}\fontseries{#4}\fontshape{#5}%
  \selectfont}%
\fi\endgroup%
{\renewcommand{\dashlinestretch}{30}
\begin{picture}(2144,659)(0,-10)
\thicklines
\path(22,22)(23,22)(26,22)
	(35,23)(51,23)(72,24)
	(97,26)(125,28)(153,31)
	(180,34)(205,37)(228,41)
	(249,46)(269,51)(287,57)
	(305,64)(322,72)(337,80)
	(353,89)(369,99)(385,109)
	(402,120)(419,131)(436,142)
	(454,154)(472,165)(490,176)
	(508,186)(525,195)(542,203)
	(559,210)(575,215)(591,219)
	(607,221)(622,222)(637,221)
	(653,219)(669,215)(685,210)
	(702,203)(719,195)(736,186)
	(754,176)(772,165)(790,154)
	(808,142)(825,131)(842,120)
	(859,109)(875,99)(891,89)
	(907,80)(922,72)(937,65)
	(953,58)(969,53)(985,48)
	(1002,43)(1019,40)(1036,38)
	(1054,36)(1072,36)(1090,36)
	(1108,38)(1125,40)(1142,43)
	(1159,48)(1175,53)(1191,58)
	(1207,65)(1222,72)(1237,80)
	(1253,89)(1269,99)(1285,109)
	(1302,120)(1319,131)(1336,142)
	(1354,154)(1372,165)(1390,176)
	(1408,186)(1425,195)(1442,203)
	(1459,210)(1475,215)(1491,219)
	(1507,221)(1522,222)(1537,221)
	(1553,219)(1569,215)(1585,210)
	(1602,203)(1619,195)(1636,186)
	(1654,176)(1672,165)(1690,154)
	(1708,142)(1725,131)(1742,120)
	(1759,109)(1775,99)(1791,89)
	(1807,80)(1822,72)(1839,64)
	(1857,57)(1875,51)(1895,46)
	(1916,41)(1939,37)(1964,34)
	(1991,31)(2019,28)(2047,26)
	(2072,24)(2093,23)(2109,23)
	(2118,22)(2121,22)(2122,22)
\path(2122,622)(2121,622)(2118,622)
	(2109,621)(2093,621)(2072,620)
	(2047,618)(2019,616)(1991,613)
	(1964,610)(1939,607)(1916,603)
	(1895,598)(1875,593)(1857,587)
	(1839,580)(1822,572)(1807,564)
	(1791,555)(1775,545)(1759,535)
	(1742,524)(1725,513)(1708,502)
	(1690,490)(1672,479)(1654,468)
	(1636,458)(1619,449)(1602,441)
	(1585,434)(1569,429)(1553,425)
	(1537,423)(1522,422)(1507,423)
	(1491,425)(1475,429)(1459,434)
	(1442,441)(1425,449)(1408,458)
	(1390,468)(1372,479)(1354,490)
	(1336,502)(1319,513)(1302,524)
	(1285,535)(1269,545)(1253,555)
	(1237,564)(1222,572)(1207,579)
	(1191,586)(1175,591)(1159,596)
	(1142,601)(1125,604)(1108,606)
	(1090,608)(1072,608)(1054,608)
	(1036,606)(1019,604)(1002,601)
	(985,596)(969,591)(953,586)
	(937,579)(922,572)(907,564)
	(891,555)(875,545)(859,535)
	(842,524)(825,513)(808,502)
	(790,490)(772,479)(754,468)
	(736,458)(719,449)(702,441)
	(685,434)(669,429)(653,425)
	(637,423)(622,422)(607,423)
	(591,425)(575,429)(559,434)
	(542,441)(525,449)(508,458)
	(490,468)(472,479)(454,490)
	(436,502)(419,513)(402,524)
	(385,535)(369,545)(353,555)
	(337,564)(322,572)(305,580)
	(287,587)(269,593)(249,598)
	(228,603)(205,607)(180,610)
	(153,613)(125,616)(97,618)
	(72,620)(51,621)(35,621)
	(26,622)(23,622)(22,622)
\end{picture}
} } 
  \end{array}
} \ar[u]
    } \\
    \calC/\calC'
  \end{array}
  & \supset & \begin{array}{c}
    \xymatrix{
      \entry{\beta} \ar[r] \ar[rd]^{\tau=d_{\star 0}\Delta^{-1}}
      & \entry{0} \\
      \entry{\alpha} \ar[u]^\Delta \ar[r]^{d_{\star 0}}
      & \entry{\tau\beta} \ar[u]
    } \\
    \calC'''
  \end{array} \\
  \parbox{1.9in}{\captionwidth=1.9in
    \caption{
      A second proof of invariance under {\rm(R2)}.
    }
  }
  & \quad & \begin{array}{c}
    \xymatrix{
      \entry{\beta} \ar[r]
        \ar[rd]^{\beta=\tau\beta}
      & \entry{0} \\
      \entry{0} \ar[u] \ar[r]
      & \entry{\gamma} \ar[u]
    } \\
    (\calC/\calC')/\calC'''
  \end{array}
\end{eqnarray*}
\end{figure}

What is $(\calC/\calC')/\calC'''$? Well, the freedom in the choice of
$\alpha$ kills the lower left corner of $\calC/\calC'$, and the freedom in
the choice of $\beta$ identifies everything in the upper left corner with
some things in the lower right corner (this is the relation
$\beta=\tau\beta$ appearing in Figure~2; in more detail it is
$(\beta,0)=(0,\tau\beta)$ in $
  \setlength{\unitlength}{0.2\standardunitlength}
  \begin{array}{c}
    {\begingroup\makeatletter\ifx\SetFigFont\undefined%
\gdef\SetFigFont#1#2#3#4#5{%
  \reset@font\fontsize{#1}{#2pt}%
  \fontfamily{#3}\fontseries{#4}\fontshape{#5}%
  \selectfont}%
\fi\endgroup%
{\renewcommand{\dashlinestretch}{30}
\begin{picture}(2144,659)(0,-10)
\thicklines
\path(22,22)(23,22)(26,22)
	(35,23)(51,23)(72,24)
	(97,26)(125,28)(153,31)
	(180,34)(205,37)(228,41)
	(249,46)(269,51)(287,57)
	(305,64)(322,72)(337,80)
	(353,89)(368,99)(384,110)
	(399,122)(415,135)(430,149)
	(444,164)(458,179)(471,195)
	(483,211)(493,227)(502,243)
	(509,260)(515,276)(519,291)
	(521,307)(522,322)(521,337)
	(519,353)(515,368)(509,384)
	(502,401)(493,417)(483,433)
	(471,449)(458,465)(444,480)
	(430,495)(415,509)(399,522)
	(384,534)(368,545)(353,555)
	(337,564)(322,572)(305,580)
	(287,587)(269,593)(249,598)
	(228,603)(205,607)(180,610)
	(153,613)(125,616)(97,618)
	(72,620)(51,621)(35,621)
	(26,622)(23,622)(22,622)
\path(722,322)(723,307)(725,291)
	(729,276)(735,260)(742,243)
	(751,227)(761,211)(773,195)
	(786,179)(800,164)(814,149)
	(829,135)(845,122)(860,110)
	(876,99)(891,89)(907,80)
	(922,72)(937,65)(953,58)
	(969,53)(985,48)(1002,43)
	(1019,40)(1036,38)(1054,36)
	(1072,36)(1090,36)(1108,38)
	(1125,40)(1142,43)(1159,48)
	(1175,53)(1191,58)(1207,65)
	(1222,72)(1237,80)(1253,89)
	(1268,99)(1284,110)(1299,122)
	(1315,135)(1330,149)(1344,164)
	(1358,179)(1371,195)(1383,211)
	(1393,227)(1402,243)(1409,260)
	(1415,276)(1419,291)(1421,307)
	(1422,322)(1421,337)(1419,353)
	(1415,368)(1409,384)(1402,401)
	(1393,417)(1383,433)(1371,449)
	(1358,465)(1344,480)(1330,495)
	(1315,509)(1299,522)(1284,534)
	(1268,545)(1253,555)(1237,564)
	(1222,572)(1207,579)(1191,586)
	(1175,591)(1159,596)(1142,601)
	(1125,604)(1108,606)(1090,608)
	(1072,608)(1054,608)(1036,606)
	(1019,604)(1002,601)(985,596)
	(969,591)(953,586)(937,579)
	(922,572)(907,564)(891,555)
	(876,545)(860,534)(845,522)
	(829,509)(814,495)(800,480)
	(786,465)(773,449)(761,433)
	(751,417)(742,401)(735,384)
	(729,368)(725,353)(723,337)(722,322)
\path(2122,622)(2121,622)(2118,622)
	(2109,621)(2093,621)(2072,620)
	(2047,618)(2019,616)(1991,613)
	(1964,610)(1939,607)(1916,603)
	(1895,598)(1875,593)(1857,587)
	(1839,580)(1822,572)(1807,564)
	(1791,555)(1776,545)(1760,534)
	(1745,522)(1729,509)(1714,495)
	(1700,480)(1686,465)(1673,449)
	(1661,433)(1651,417)(1642,401)
	(1635,384)(1629,368)(1625,353)
	(1623,337)(1622,322)(1623,307)
	(1625,291)(1629,276)(1635,260)
	(1642,243)(1651,227)(1661,211)
	(1673,195)(1686,179)(1700,164)
	(1714,149)(1729,135)(1745,122)
	(1760,110)(1776,99)(1791,89)
	(1807,80)(1822,72)(1839,64)
	(1857,57)(1875,51)(1895,46)
	(1916,41)(1939,37)(1964,34)
	(1991,31)(2019,28)(2047,26)
	(2072,24)(2093,23)(2109,23)
	(2118,22)(2121,22)(2122,22)
\end{picture}
} } 
  \end{array}
\!\!_{/v_+=0} \oplus

  \setlength{\unitlength}{0.2\standardunitlength}
  \begin{array}{c}
    {\begingroup\makeatletter\ifx\SetFigFont\undefined%
\gdef\SetFigFont#1#2#3#4#5{%
  \reset@font\fontsize{#1}{#2pt}%
  \fontfamily{#3}\fontseries{#4}\fontshape{#5}%
  \selectfont}%
\fi\endgroup%
{\renewcommand{\dashlinestretch}{30}
\begin{picture}(2144,659)(0,-10)
\thicklines
\path(22,22)(23,22)(26,22)
	(35,23)(51,23)(72,24)
	(97,26)(125,28)(153,31)
	(180,34)(205,37)(228,41)
	(249,46)(269,51)(287,57)
	(305,64)(322,72)(337,80)
	(353,89)(369,99)(385,109)
	(402,120)(419,131)(436,142)
	(454,154)(472,165)(490,176)
	(508,186)(525,195)(542,203)
	(559,210)(575,215)(591,219)
	(607,221)(622,222)(637,221)
	(653,219)(669,215)(685,210)
	(702,203)(719,195)(736,186)
	(754,176)(772,165)(790,154)
	(808,142)(825,131)(842,120)
	(859,109)(875,99)(891,89)
	(907,80)(922,72)(937,65)
	(953,58)(969,53)(985,48)
	(1002,43)(1019,40)(1036,38)
	(1054,36)(1072,36)(1090,36)
	(1108,38)(1125,40)(1142,43)
	(1159,48)(1175,53)(1191,58)
	(1207,65)(1222,72)(1237,80)
	(1253,89)(1269,99)(1285,109)
	(1302,120)(1319,131)(1336,142)
	(1354,154)(1372,165)(1390,176)
	(1408,186)(1425,195)(1442,203)
	(1459,210)(1475,215)(1491,219)
	(1507,221)(1522,222)(1537,221)
	(1553,219)(1569,215)(1585,210)
	(1602,203)(1619,195)(1636,186)
	(1654,176)(1672,165)(1690,154)
	(1708,142)(1725,131)(1742,120)
	(1759,109)(1775,99)(1791,89)
	(1807,80)(1822,72)(1839,64)
	(1857,57)(1875,51)(1895,46)
	(1916,41)(1939,37)(1964,34)
	(1991,31)(2019,28)(2047,26)
	(2072,24)(2093,23)(2109,23)
	(2118,22)(2121,22)(2122,22)
\path(2122,622)(2121,622)(2118,622)
	(2109,621)(2093,621)(2072,620)
	(2047,618)(2019,616)(1991,613)
	(1964,610)(1939,607)(1916,603)
	(1895,598)(1875,593)(1857,587)
	(1839,580)(1822,572)(1807,564)
	(1791,555)(1775,545)(1759,535)
	(1742,524)(1725,513)(1708,502)
	(1690,490)(1672,479)(1654,468)
	(1636,458)(1619,449)(1602,441)
	(1585,434)(1569,429)(1553,425)
	(1537,423)(1522,422)(1507,423)
	(1491,425)(1475,429)(1459,434)
	(1442,441)(1425,449)(1408,458)
	(1390,468)(1372,479)(1354,490)
	(1336,502)(1319,513)(1302,524)
	(1285,535)(1269,545)(1253,555)
	(1237,564)(1222,572)(1207,579)
	(1191,586)(1175,591)(1159,596)
	(1142,601)(1125,604)(1108,606)
	(1090,608)(1072,608)(1054,608)
	(1036,606)(1019,604)(1002,601)
	(985,596)(969,591)(953,586)
	(937,579)(922,572)(907,564)
	(891,555)(875,545)(859,535)
	(842,524)(825,513)(808,502)
	(790,490)(772,479)(754,468)
	(736,458)(719,449)(702,441)
	(685,434)(669,429)(653,425)
	(637,423)(622,422)(607,423)
	(591,425)(575,429)(559,434)
	(542,441)(525,449)(508,458)
	(490,468)(472,479)(454,490)
	(436,502)(419,513)(402,524)
	(385,535)(369,545)(353,555)
	(337,564)(322,572)(305,580)
	(287,587)(269,593)(249,598)
	(228,603)(205,607)(180,610)
	(153,613)(125,616)(97,618)
	(72,620)(51,621)(35,621)
	(26,622)(23,622)(22,622)
\end{picture}
} } 
  \end{array}
$). What remains is just the
arbitrary choice of $\gamma$ in the lower right corner and hence
the complex $(\calC/\calC')/\calC'''$ is isomorphic to the complex
$\calC''$ of Figure~1 and this, as there, is what we wanted to prove. \qed

\subsubsection{Invariance under (R3).} We can now turn back to the
proof of invariance under (R3). Repeat the definitions of the acyclic
subcomplexes $\calC'$ and $\calC'''$ as above but within the top layers
of each of the cubes in Equation~\eqref{eq:R3Cubes}, and then mod out each
cube by its $\calC'$ and $\calC'''$ (without changing the homology, by
Lemma~\ref{lem:Cancellation}). The resulting cubes are

\noindent{\centering\resizebox*{\textwidth}{!}{\parbox{6.2in}{\[
  \def\neg{\hspace{-2mm}}
  \def\D#1{\neg
  \setlength{\unitlength}{0.125\standardunitlength}
  \begin{array}{c}
    {\input figs/D#1.tex } 
  \end{array}
\neg}
  \def\U#1{\neg
  \setlength{\unitlength}{0.125\standardunitlength}
  \begin{array}{c}
    {\input figs/U#1.tex } 
  \end{array}
\neg}
  \begin{array}{c} \xymatrix@R=3.5mm@C=3.5mm{
    & \beta_1\in\D{101}_{/v_+=0} \ar@{->}[rr]
        \ar@{<-}'[d][dd]^(0.4){d_{1,\star 01}}
        \ar[dr]^{\beta_1=\tau_1\beta_1}
      & & 0 \ar@{<-}[dd] \\
    0 \ar@{->}[ur] \ar@{->}[rr] \ar@{<-}[dd] &
      & \gamma_1\in\D{110} \ar@{->}[ur] \ar@{<-}[dd]_(0.7){d_{1,\star 10}} \\
    & \D{001} \ar@{->}'[r][rr] & & \D{011} \\
    \D{000} \ar@{->}[rr] \ar@{->}[ur] & & \D{010} \ar@{->}[ur]
  } \end{array}
  \ \ 
  \begin{array}{c} \xymatrix@R=4mm@C=4mm{
    & \gamma_2\in\U{101} \ar@{->}[rr]
        \ar@{<-}'[d][dd]^(0.4){d_{2,\star 01}}
      & & 0 \ar@{<-}[dd] \\
    0 \ar@{->}[ur] \ar@{->}[rr] \ar@{<-}[dd] &
      & \beta_2\in\U{110}_{/v_+=0} \ar@{->}[ur]
        \ar@{<-}[dd]_(0.7){d_{2,\star 10}}
        \ar[lu]^{\tau_2\beta_2=\beta_2} \\
    & \U{001} \ar@{->}'[r][rr] & & \U{011} \\
    \U{000} \ar@{->}[rr] \ar@{->}[ur] & & \U{010} \ar@{->}[ur]
  } \end{array}
\]}}}

Now these two complexes really are isomorphic, via the map
$\Upsilon$ that keeps the bottom layers in place and ``transposes''
the top layers by mapping the pair $(\beta_1,\gamma_1)$ to the pair
$(\beta_2,\gamma_2)$. The fact that $\Upsilon$ is an isomorphism on
spaces level is obvious. To see that $\Upsilon$ is an isomorphism of
complexes we need to know that it commutes with the edge maps, and
only the vertical edges require a proof. We leave the (easy) proofs
that $\tau_1\circ d_{1,\star 01} = d_{2,\star 01}$ and $d_{1,\star
10}=\tau_2\circ d_{2,\star 10}$ as exercises for our readers. \qed

\subsection{Some phenomenological conjectures}

The following conjectures were formulated in parts by the author and by
M.~Khovanov and S.~Garoufalidis based on computations using the program
described in the next section:

\begin{conjecture} \label{conj:HRestricted} For any prime knot $L$ there
exists an even integer $s=s(L)$ and a polynomial $\Kh'(L)$ in $t^{\pm 1}$
and $q^{\pm 1}$ with only non-negative coefficients so that
\begin{eqnarray}
  \Kh_\bbQ(L) &=& q^{s-1}\left(1+q^2 + (1+tq^4)\Kh'(L)\right)
    \label{eq:QConjecture} \\
  \Kh_{\bbF_2}(L) &=& q^{s-1}(1+q^2)\left(1+(1+tq^2)\Kh'(L)\right).
    \label{F2Conjecture}
\end{eqnarray}
($\bbF_2$ denotes the field of two elements.)
\end{conjecture}

\begin{conjecture} \label{conj:HThin} For prime alternating $L$ the integer
$s(L)$ is equal to the signature of $L$ and the polynomial $\Kh'(L)$
contains only powers of $tq^2$.
\end{conjecture}

We have computed $\Kh_\bbQ(L)$ for all prime knots with up to 11 crossings
and $\Kh_{\bbF_2}(L)$ for all knots with up to 7 crossings and the results
are in complete agreement with these two conjectures\footnote{Except
that for 11 crossing prime alternating knots only the absolute values
of $\sigma$ and $s$ were compared.}.

We note that these conjectures imply that for alternating knots $\Kh'$
(and hence $\Kh_\bbQ$ and $\Kh_{\bbF_2}$) are determined by the Jones
polynomial.  As we shall see in the next section, this is not true for
non-alternating knots.

\parpic[r]{\scriptsize
  \begin{tabular}{|c|c|c|c|c|c|c|c|c|c|c|c|}
    \hline
    \backslashbox{\!$m$\!}{\!$r$\!} & \!-7\! & \!-6\! & \!-5\! & \!-4\! &
      \!-3\! & \!-2\! & \!-1\! & 0 & 1 & 2 & 3 \\ \hline
    3 & & & & & & & & & & & 1 \\ \hline
    1 & & & & & & & & & & 2 & \\ \hline
    -1 & & & & & & & & & 3 & 1 & \\ \hline
    -3 & & & & & & & & \!\!(4+1)\!\! & 2 & & \\ \hline
    -5 & & & & & & & 5 & \!\!(3+1)\!\! & & & \\ \hline
    -7 & & & & & & 6 & 4 & & & & \\ \hline
    -9 & & & & & 4 & 5 & & & & & \\ \hline
    -11 & & & & 4 & 6 & & & & & & \\ \hline
    -13 & & & 2 & 4 & & & & & & & \\ \hline
    -15 & & 1 & 4 & & & & & & & & \\ \hline
    -17 & & 2 & & & & & & & & & \\ \hline
    -19 & 1 & & & & & & & & & & \\ \hline
  \end{tabular}
}
As a graphical illustration of Conjectures \ref{conj:HRestricted} and
\ref{conj:HThin} the table on the right contains the dimensions of the
homology groups $\calH^r_m(10_{100})$ (the coefficients of $t^rq^m$ in
the invariant $\Kh(10_{100})$) for all $r$ and $m$ in the relevant
range.  Conjecture~\ref{conj:HRestricted} is the fact that if we
subtract 1 from two of the entries in the column $r=0$ (a ``pawn
move''), the remaining entries are arranged in ``knight move'' pairs of
the form {\scriptsize $\begin{array}{|c|c|}\hline&a\\\hline\ &\\\hline
a&\\\hline\end{array}$} with $a>0$. Conjecture~\ref{conj:HThin} is the
fact that furthermore all nontrivial entries in the table occur on just
two diagonals that cross the column $r=0$ at $m=\sigma\pm 1$ where
$\sigma=-4$ is the signature of $10_{100}$. Thus after the fix at the
$r=0$ column, the two nontrivial diagonals are just shifts of each
other and are thus determined by a single list of entries (1 2 4 4 6 5
4 3 2 1, in our case). This list of entries is the list of coefficients
of $\Kh'(10_{100}) = u^{-7} + 2u^{-6} + 4u^{-5} + 4u^{-4} + 6u^{-3} +
5u^{-2} + 4u^{-1} + 3 + 2u + u^2$ (with $u=tq^4$).

As an aside we note that typically $\dim\calH^r_m(L)$ is much smaller than
$\dim\calC^r_m(L)$, as illustrated in Table~\ref{table:hc}. We don't know
why this is so.

\begin{table}
\begin{center}{\scriptsize
  \def\neg{\hspace{-0.8mm}}
  \def\c#1{\neg 0/#1\neg}
  \def\hc#1#2{\neg{\textbf{#1/#2}}\neg}
  \begin{tabular}{|c|c|c|c|c|c|c|c|c|c|c|c|}
    \hline
    \backslashbox{$m$}{$r$} & -7 & -6 & -5 & -4 & -3 & -2 & -1 & 0 & 1
      & 2 & 3 \\ \hline
    3 & & & & & & & & & & & \hc{1}{1} \\ \hline
    1 & & & & & & & & \c{1} & \c{5} & \hc{2}{10} & \c{4} \\ \hline
    -1 & & & & & \c{2} & \c{13} & \c{36} & \c{59} & \hc{3}{60} & 
      \hc{1}{30} & \c{6} \\ \hline
    -3 & \c{1} & \c{10} & \c{45} & \c{120} & \c{220} & 
      \c{304} & \c{318} & \hc{5}{237} & \hc{2}{110} & \c{30} & 
      \c{4} \\ \hline
    -5 & \c{8} & \c{70} & \c{270} & \c{600} & \c{862} & 
      \c{847} & \hc{5}{564} & \hc{4}{237} & \c{60} & \c{10} & 
      \c{1} \\ \hline
    -7 & \c{28} & \c{210} & \c{675} & \c{1200} & \c{1288} & 
      \hc{6}{847} & \hc{4}{318} & \c{59} & \c{5} & & \\ \hline
    -9 & \c{56} & \c{350} & \c{900} & \c{1200} & \hc{4}{862} & 
      \hc{5}{304} & \c{36} & \c{1} & & & \\ \hline
    -11 & \c{70} & \c{350} & \c{675} & \hc{4}{600} & \hc{6}{220} & 
      \c{13} & & & & & \\ \hline
    -13 & \c{56} & \c{210} & \hc{2}{270} & \hc{4}{120} & 
      \c{2} & & & & & & \\ \hline
    -15 & \c{28} & \hc{1}{70} & \hc{4}{45} & & & & & & & & \\ \hline
    -17 & \c{8} & \hc{2}{10} & & & & & & & & & \\ \hline
    -19 & \hc{1}{1} & & & & & & & & & & \\ \hline
  \end{tabular}
}\end{center}
\caption{
  $\dim\calH^r_m(10_{100})/\dim\calC^r_m(10_{100})$ for all values of
  $r$ and $m$ for which $\calC^r_m(10_{100})\neq\emptyset$.
} \label{table:hc}
\end{table}

A further phenomenological conjecture is presented
in~\cite{Garoufalidis:ConjectureOnKhovanov}. This paper's web
page~\cite{cite:This} will follow further phenomenological developments
as they will be announced.

\section{And now in computer talk} \label{sec:ComputerTalk}

\begin{figure}
\[ \includegraphics[width=4in]{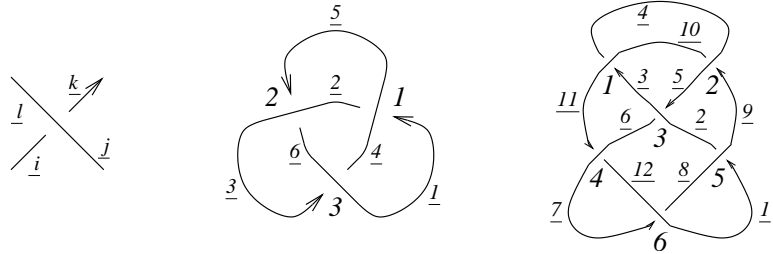} \]
\caption{
  The crossing $X_{ijkl}$, the right handed trefoil knot
  $X_{1524}X_{5362}X_{3146}$ and the Miller Institute knot (aka
  $\bar{6}_2$)
  $X_{3,10,4,11}X_{9,4,10,5}X_{5,3,6,2}X_{11,7,12,6}X_{1,9,2,8}X_{7,1,8,12}$
  (we've used a smaller font and underlining to separate the edge
  labeling from the vertex labeling).
} \label{fig:LinkNotation}
\end{figure}

In computer talk (Mathematica~\cite{Wolfram:Mathematica} dialect) we
represent every link projection by a list of edges numbered $1,\dots,n$
with increasing numbers as we go around each component and by a list
crossings presented as symbols $X_{ijkl}$ where $i,\dots,l$ are the
edges around that crossing, starting from the incoming lower thread and
proceeding counterclockwise (see Figure~\ref{fig:LinkNotation}).

\subsection{A demo run} We first start up
Mathematica~\cite{Wolfram:Mathematica} and load our categorification
package, \verb$Categorification`$ (available from~\cite{cite:This}):

{\small
\begin{verbatim}
Mathematica 4.1 for Linux
Copyright 1988-2000 Wolfram Research, Inc.
 -- Motif graphics initialized --
\end{verbatim}
}
{\In << Categorification`}
\begingroup\Print Loading Categorification`...\endgroup

\vskip 2mm
Next, we type in the trefoil knot:

{\In L = Link[X[1,5,2,4], X[5,3,6,2], X[3,1,4,6]];}

\vskip 2mm
Let us now view the edge $0\!\star\!1$ of the cube of smoothings of the
trefoil knot (as seen in Section~\ref{subsec:Notation}, this edge begins
with a single cycle labeled $1$ and ends with two cycles labeled $1$ and
$3$):

{\In \{S[L, "001"], S[L, "0*1"], S[L, "011"]\}}
\begingroup\Out {c[1], c[1] -> c[1]*c[3], c[1]*c[3]}\endgroup

\vskip 2mm
Next, here's a basis of the space $V_{011}$ (again, compare with
Section~\ref{subsec:Notation}):

{\In V[L, "011"]}
\begingroup\Out {vm[1]*vm[3], vm[3]*vp[1], vm[1]*vp[3], vp[1]*vp[3]}\endgroup

\vskip 2mm
And here's a basis of the degree $2$ elements of $V_{111}$ (remember the
shift in degrees in the definition of $V_\alpha$!):

{\In V[L, "111", 2]}
\begingroup\Out {vm[2]*vm[3]*vp[1], vm[1]*vm[3]*vp[2], vm[1]*vm[2]*vp[3]}\endgroup

\vskip 2mm
The per-edge map $d_\xi$ is a list of simple replacement rules,
sometimes replacing the tensor product of two basis vectors by a single
basis vector, as in the case of $d_{00\star}=m_{12}$, and sometimes
the opposite, as in the case of $d_{0\star 1}=\Delta^{13}$:

{\In d[L, "00*"]}
\begingroup\Out {vp[1]*vp[2] -> vp[1], vm[2]*vp[1] -> vm[1], vm[1]*vp[2] -> vm[1], 

   vm[1]*vm[2] -> 0}\endgroup
{\In d[L, "0*1"]}
\begingroup\Out {vp[1] -> vm[3]*vp[1] + vm[1]*vp[3], vm[1] -> vm[1]*vm[3]}\endgroup

\vskip 2mm
Here's a simple example. Let us compute $d_{1\star 1}$ applied to
$V_{101}$:

{\In V[L, "101"] /. d[L, "1*1"]}
\begingroup\Out {vm[1]*vm[2]*vm[3], vm[2]*vm[3]*vp[1], 

   vm[1]*(vm[3]*vp[2] + vm[2]*vp[3]), vp[1]*(vm[3]*vp[2] + \
vm[2]*vp[3])}\endgroup

\vskip 2mm
And now a more complicated example. First, we compute the degree $0$
part of $\llbracket\righttrefoil\rrbracket^1$. Then we apply $d^1$
to it, and then $d^2$ to the result. The end result better be a
list of zeros, or else we are in trouble! Notice that each basis
vector in $\llbracket\righttrefoil\rrbracket^{1,2}$ is tagged with
a symbol of the form \verb$v[...]$ that indicates the component of
$\llbracket\righttrefoil\rrbracket^{1,2}$ to which it belongs.

{\In chains = KhBracket[L, 1, 0]}
\begingroup\Out {v[0, 0, 1]*vm[1], v[0, 1, 0]*vm[1], v[1, 0, 0]*vm[1]}\endgroup
{\In boundaries = d[L][chains]}
\begingroup\Out {v[1, 0, 1]*vm[1]*vm[2] + v[0, 1, 1]*vm[1]*vm[3], 

   v[1, 1, 0]*vm[1]*vm[2] - v[0, 1, 1]*vm[1]*vm[3], 

   -(v[1, 0, 1]*vm[1]*vm[2]) - v[1, 1, 0]*vm[1]*vm[2]}\endgroup
{\In d[L][boundaries]}
\begingroup\Out {0, 0, 0}\endgroup

\vskip 2mm
Because of degree shifts, the degree $3$ part of $\calC^1(\righttrefoil)$
is equal to the degree $0$ part of $\llbracket\righttrefoil\rrbracket^1$:

{\In CC[L, 1, 3] == KhBracket[L, 1, 0]}
\begingroup\Out True\endgroup

\vskip 2mm
It seems that $\calH^2(\righttrefoil)$ is one dimensional, and that the
non trivial class in $\calH^2(\righttrefoil)$ lies in degree $5$ (our
program defaults to computations over the rational numbers if no other
modulus is specified):

{\In qBetti[L, 2]}
\begingroup\Out q^5\endgroup

\vskip 2mm
Here's Khovanov's invariant of the right handed trefoil along if its
evaluation at $t=-1$, the unnormalized Jones polynomial
$\hatJ(\righttrefoil)$:

{\In kh1 = Kh[L]}
\begingroup\Out q + q^3 + q^5*t^2 + q^9*t^3\endgroup
{\In kh1 /. t -> -1}
\begingroup\Out q + q^3 + q^5 - q^9\endgroup

\vskip 2mm
We can also compute $\Kh_{\bbF_2}(\righttrefoil)$ and use it to compute
$\hatJ(\righttrefoil)$ again (we leave it to the reader to verify
Conjecture~\ref{conj:HRestricted} in the case of $L=\righttrefoil$):

{\In kh2 = Kh[L, Modulus -> 2]}
\begingroup\Out q + q^3 + q^5*t^2 + q^7*t^2 + q^7*t^3 + q^9*t^3\endgroup
{\In kh2 /. t -> -1}
\begingroup\Out q + q^3 + q^5 - q^9\endgroup

\vskip 2mm
\parpic[r]{\begin{tabular}{c}
  \includegraphics[width=1.3in]{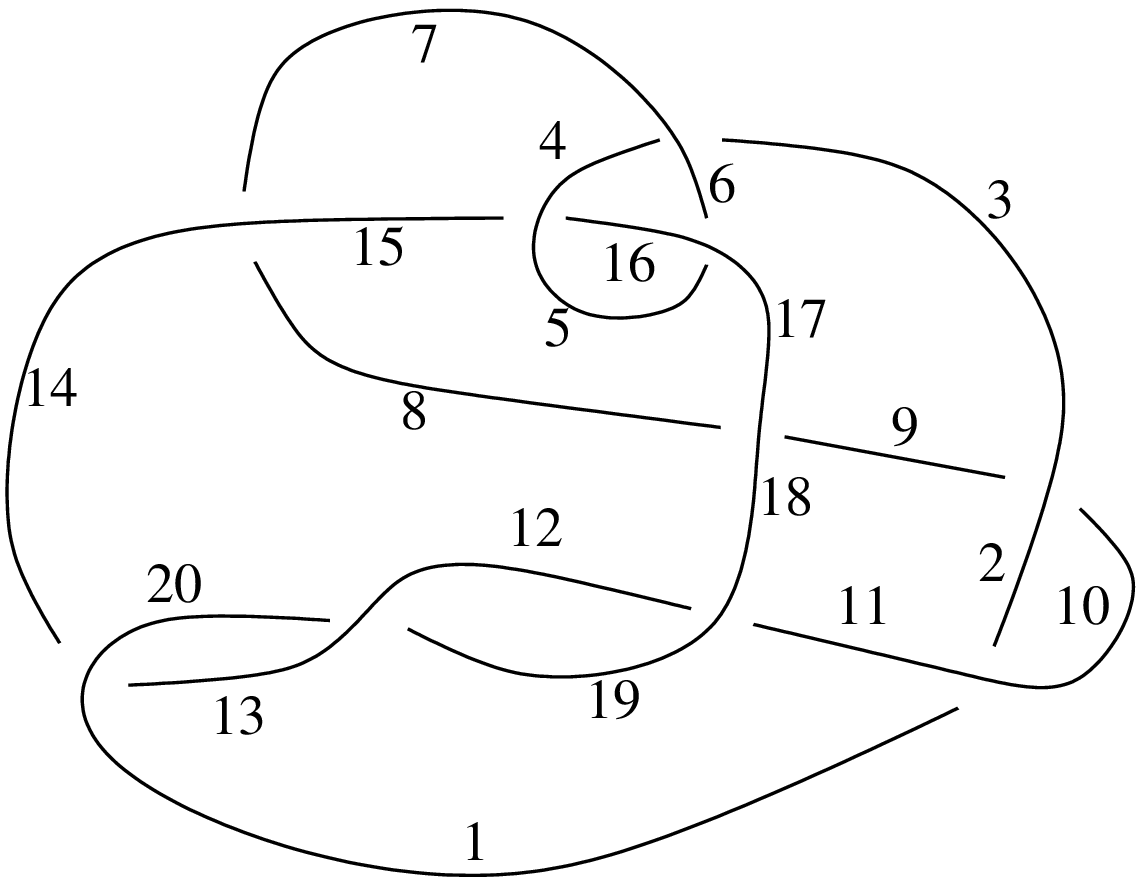}
\end{tabular}}
The package \verb$Links`$ (available from~\cite{cite:This}) contains
the definitions of many interesting knot and link projections including
Millett's 10 crossing hard-to-simplify unknot (shown on the right) and the knots $5_1$ and
$10_{132}$ (knot numbering as in Rolfsen's~\cite{Rolfsen:KnotsLinks}):

{\In << Links`} 
\begingroup\Print Loading Links`...\endgroup
{\In MillettUnknot}
\begingroup\Out Link[X[1, 10, 2, 11], X[9, 2, 10, 3], X[3, 7, 4, 6], X[15, 5, 16, 4], \

 X[5, 17, 6, 16], X[7, 14, 8, 15], X[8, 18, 9, 17], 

   X[11, 18, 12, 19], X[19, 12, 20, 13], X[13, 20, 14, 1]]
\endgroup
{\In Kh[MillettUnknot]}
\begingroup\Out q^(-1) + q\endgroup
{\In kh3 = Kh[Knot[5, 1]]}
\begingroup\Out q^(-5) + q^(-3) + 1/(q^15*t^5) + 1/(q^11*t^4) + 1/(q^11*t^3) + 

   1/(q^7*t^2)\endgroup
{\In kh4 = Kh[Knot[10, 132]]}
\begingroup\Out q^(-3) + q^(-1) + 1/(q^15*t^7) + 1/(q^11*t^6) + 1/(q^11*t^5) + 

   1/(q^9*t^4) + 1/(q^7*t^4) + 1/(q^9*t^3) + 1/(q^5*t^3) + 

   2/(q^5*t^2) + 1/(q*t)\endgroup
{\In (kh3 /. t -> -1) == (kh4 /. t -> -1)}
\begingroup\Out True\endgroup

\vskip 2mm
These are excellent news! We have just learned that our program is not
confused by complicated mess, and even better, we have just learned that
Khovanov's invariant is strictly stronger than the Jones polynomial, for
$J(5_1)=J(10_{132})$ whereas $\Kh(5_1)\neq\Kh(10_{132})$.

Here are two further pieces of good news:
{\In diff1 = Together[Kh[Knot[9, 42]] - Kh[Mirror[Knot[9, 42]]]]}
\begingroup\Out (1 + q^4*t - t^2 + q^4*t^2 - q^4*t^3 + q^6*t^3 + q^8*t^3 - q^4*t^4 + 

    q^10*t^4 - q^6*t^5 - q^8*t^5 + q^10*t^5 - q^10*t^6 + q^14*t^6 - 

    q^10*t^7 - q^14*t^8)/(q^7*t^4)\endgroup
{\In diff2 = Expand[q^9*t^5*(Kh[Knot[10, 125]]-Kh[Mirror[Knot[10, 125]]])]}
\begingroup\Out 1 + q^4*t - t^2 + q^4*t^2 - q^4*t^3 + q^6*t^3 + q^8*t^3 - q^4*t^4 + 

   q^10*t^4 - q^6*t^5 - 2*q^8*t^5 + 2*q^10*t^5 + q^12*t^5 - q^8*t^6 + \

 q^14*t^6 - q^10*t^7 - q^12*t^7 + q^14*t^7 - q^14*t^8 + q^18*t^8 - 

   q^14*t^9 - q^18*t^10\endgroup
{\In \{diff1, diff2\} /. t -> -1}
\begingroup\Out {0, 0}\endgroup

\vskip 2mm Thus we see that $\Kh$ detects the facts that
$9_{42}\neq\overline{9_{42}}$ and $10_{125}\neq\overline{10_{125}}$
whereas the Jones polynomial doesn't detect that. See also
Section~\ref{subsec:NewSeparations}.

\subsection{The program}

The program {\tt Categorification.m} and the data files {\tt Data.m}
and {\tt Links.m} demonstrated in this article are available at

{\tt http://www.maths.warwick.ac.uk/agt/ftp/aux/agt-2-16/}

(with a link from the home page of this paper) and also from the arXiv
at at~\cite{cite:This}.  A complete listing of the package {\tt
Categorification.m} takes up less than 70 lines of code, demonstrating
that categorification must be quite simple.

\subsection{$\Kh'(L)$ for prime knots with up to 10 crossings}


Conjecture~\ref{conj:HRestricted} on page~\pageref{conj:HRestricted}
introduces an integer $s=s(L)$ and a polynomial $\Kh'(L)$. By direct
computation using our program we verified that these quantities are
determined by $\Kh_\bbQ(L)$ for all knots with up to 11 crossings. These
quantities easily determine $\Kh_\bbQ(L)$ (and also $\Kh_{\bbF_2}(L)$,
at least up to knots with 7 crossings), as in the statement of
Conjecture~\ref{conj:HRestricted}.

There are many fewer terms in $\Kh'(L)$ as there are in $\Kh_\bbQ(L)$
or in $\Kh_{\bbF_2}(L)$ and thus with the rain forests in our minds,
we've tabulated $s$ and $\Kh'(L)$ rather than $\Kh_\bbQ(L)$ and/or
$\Kh_{\bbF_2}(L)$. To save further space, we've underlined negative
numbers ($\underline{1}:=-1$), used the notation $a^r_m$ to denote the
monomial $at^rq^m$ and suppressed all ``$+$'' signs. Thus $\Kh'(7_7) =
\frac{1}{q^6t^3} + \frac{2}{q^4t^2} + \frac{1}{q^2t} + 2 + 2q^2t
+ q^4t^2 + q^6t^3$ is printed as $1^{\underline{3}}_{\underline{6}}2^{\underline{2}}_{\underline{4}}1^{\underline{1}}_{\underline{2}}2^{0}_{0}2^{1}_{2}1^{2}_{4}1^{3}_{6}$.

Staring at the tables below it is difficult not to notice that $s(L)$
is often equal to the signature $\sigma=\sigma(L)$ of $L$, and that
most monomials in most $\Kh'(L)'s$ are of the form $t^rq^{2r}$ for some
$r$. We've marked the exceptions to the first observation by the flag
$\clubsuit$ and the knots where exceptions to the second observation occur
by the flag $\spadesuit$. All exceptions occur at non-alternating knots.
(And for your convenience, these are marked by the flag $\diamondsuit$).

\par\noindent{\bf Acknowledgement and Warning.} The combinatorial data
on which I based the computations was provided to me by A.~Stoimenow
(see~\cite{Stoimenow:Polynomials}), who himself borrowed it from J.~Hoste
and M.~Thistlethwaite~\cite{HosteThistlethwaite:Knotscape}, and was
translated to our format by a program written by D.~Thurston. The
knot pictures below were generated using R.~Scharein's program
KnotPlot~\cite{Scharein:KnotPlot}. The assembly of all this information
involved some further programming and manual work. I hope that no errors
crept through, but until everything is independently verified, I cannot
be sure of that. I feel that perhaps other than orientation issues (some
knots may have been swapped with the mirrors) the data below is reliable.
Finally, note that we number knots as in
Rolfsen's~\cite{Rolfsen:KnotsLinks}, except that we have removed $10_{162}$
which is equal to $10_{161}$ (this is the famed ``Perko pair''). Hence
Rolfsen's $10_{163,\ldots,166}$ are ours $10_{162,\ldots,165}$.

All data shown here is available in computer readable format at~\cite[the
file {\tt Data.m}]{cite:This}.

\begin{center} \small
\def\neg{\hspace{-3mm}}

\begin{longtable}{|cl|cl|}
  \hline
    $L$
    & \neg $\begin{array}{l}
      \parbox{13mm}{$n_k$;} \sigma; \quad s; \quad \text{flags} \\[1.5pt]
      \Kh'(L)
    \end{array}$
    & $L$
    & \neg $\begin{array}{l}
      \parbox{13mm}{$n_k$;} \sigma; \quad s;\quad \text{flags} \\[1.5pt]
      \Kh'(L)
    \end{array}$
    \\
  \hline \endhead
  \hline \endfoot


\parbox[c]{12mm}{     \includegraphics[height=12mm,width=12mm,keepaspectratio]{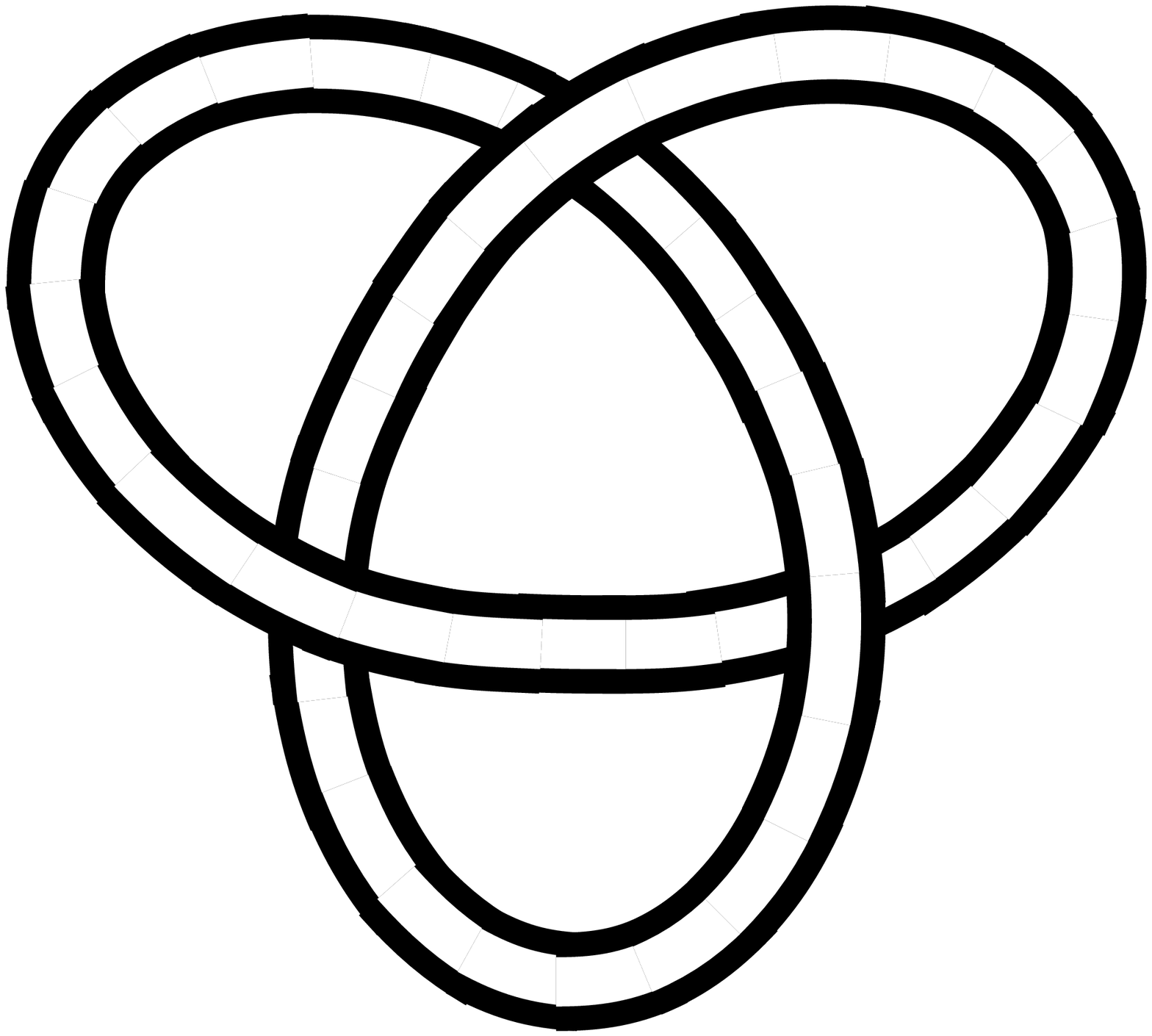}}& \neg $\begin{array}{l}\parbox{13mm}{$3_{1}$;} \underline{2}; \quad \underline{2}\\[1.5pt] 1^{\underline{3}}_{\underline{6}}\neg \end{array}$ \neg& 
& \\ \hline 

 \hline
\parbox[c]{12mm}{     \includegraphics[height=12mm,width=12mm,keepaspectratio]{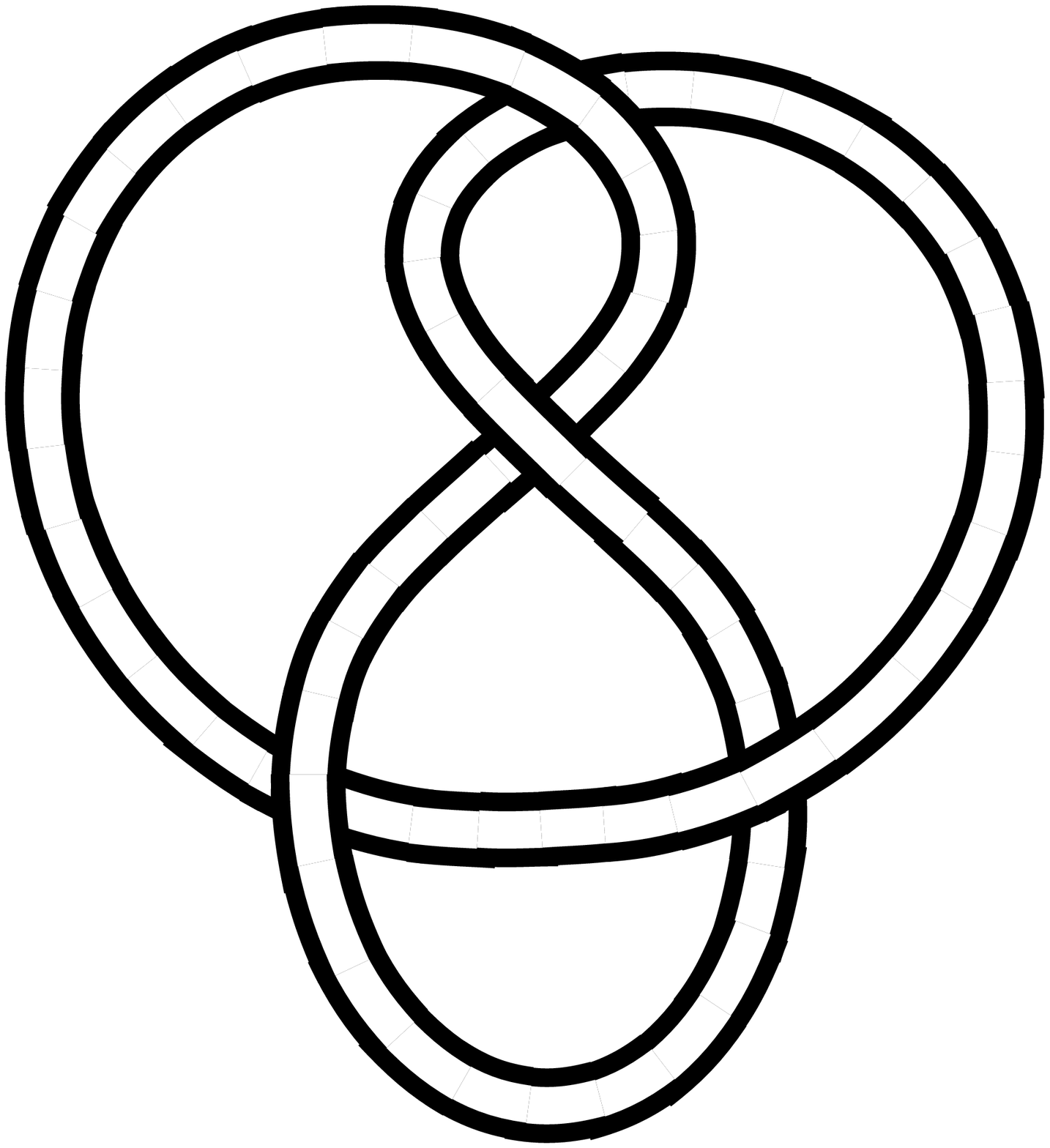}}& \neg $\begin{array}{l}\parbox{13mm}{$4_{1}$;} 0; \quad 0\\[1.5pt] 1^{\underline{2}}_{\underline{4}}1^{1}_{2}\neg \end{array}$ \neg& 
& \\ \hline 

 \hline
\parbox[c]{12mm}{     \includegraphics[height=12mm,width=12mm,keepaspectratio]{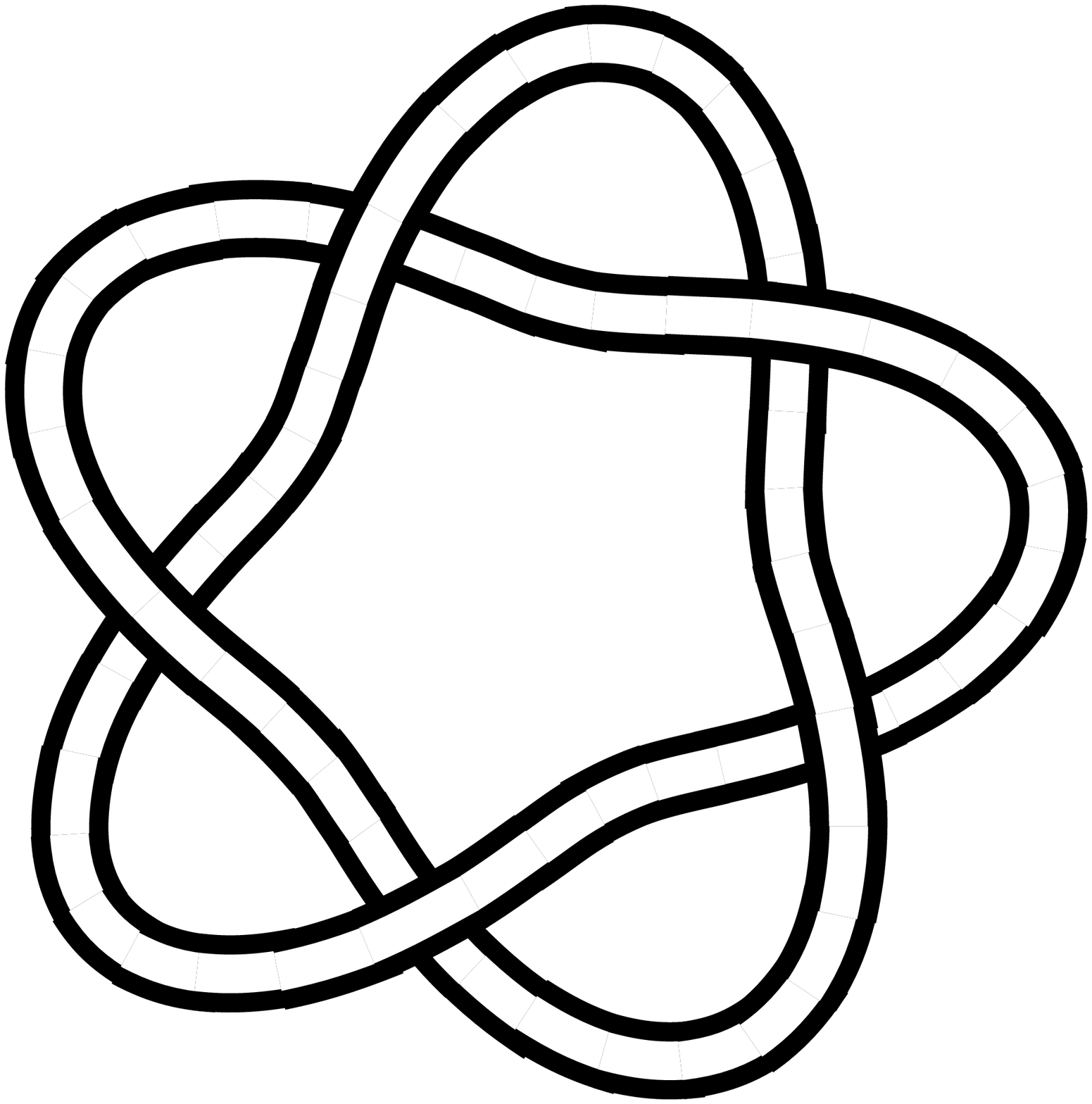}}& \neg $\begin{array}{l}\parbox{13mm}{$5_{1}$;} \underline{4}; \quad \underline{4}\\[1.5pt] 1^{\underline{5}}_{\underline{10}}1^{\underline{3}}_{\underline{6}}\neg \end{array}$ \neg& 
\parbox[c]{12mm}{     \includegraphics[height=12mm,width=12mm,keepaspectratio]{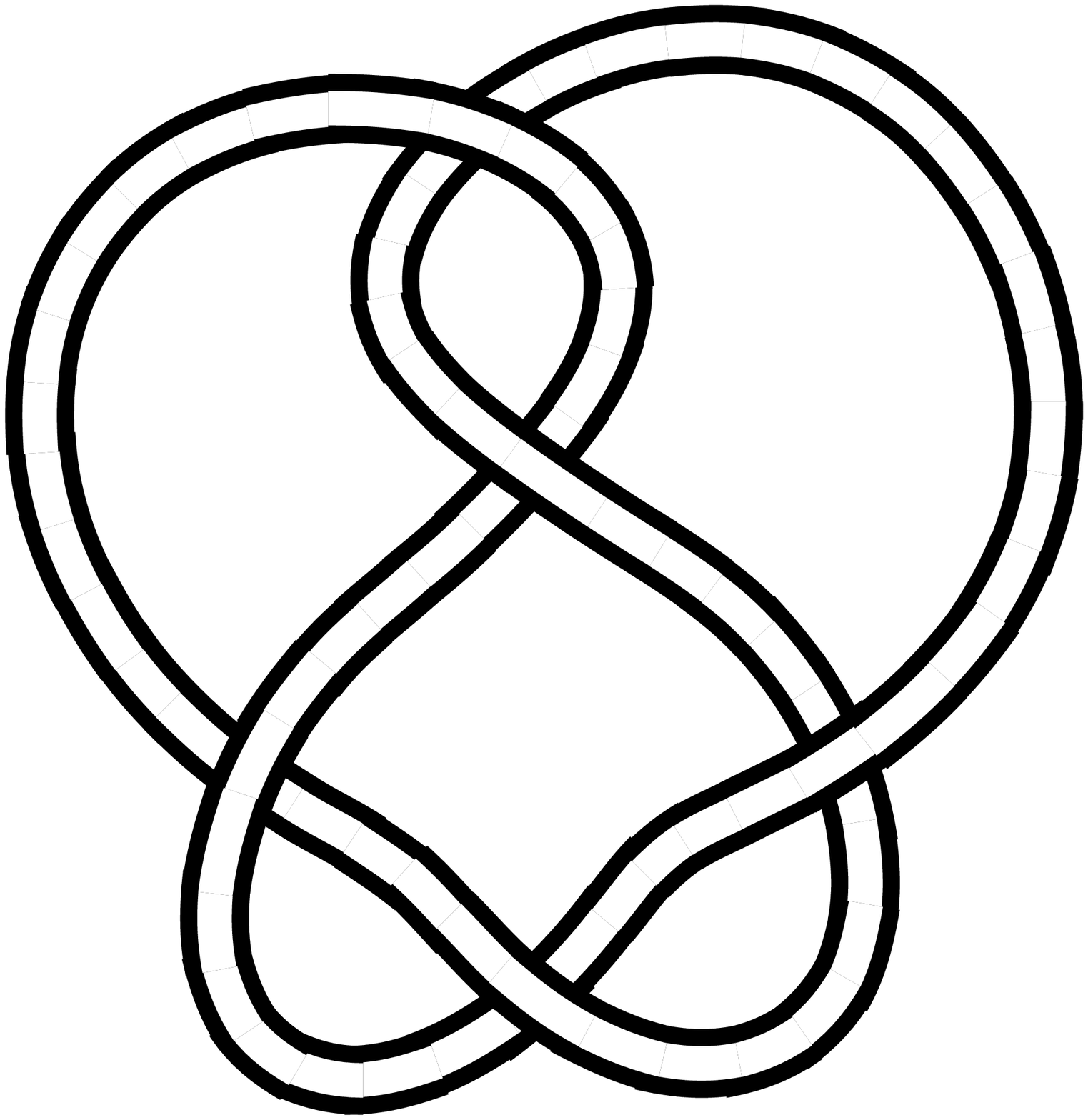}}& \neg $\begin{array}{l}\parbox{13mm}{$5_{2}$;} \underline{2}; \quad \underline{2}\\[1.5pt] 1^{\underline{5}}_{\underline{10}}1^{\underline{3}}_{\underline{6}}1^{\underline{2}}_{\underline{4}}\neg \end{array}$ \neg\\ \hline 

 \hline
\parbox[c]{12mm}{     \includegraphics[height=12mm,width=12mm,keepaspectratio]{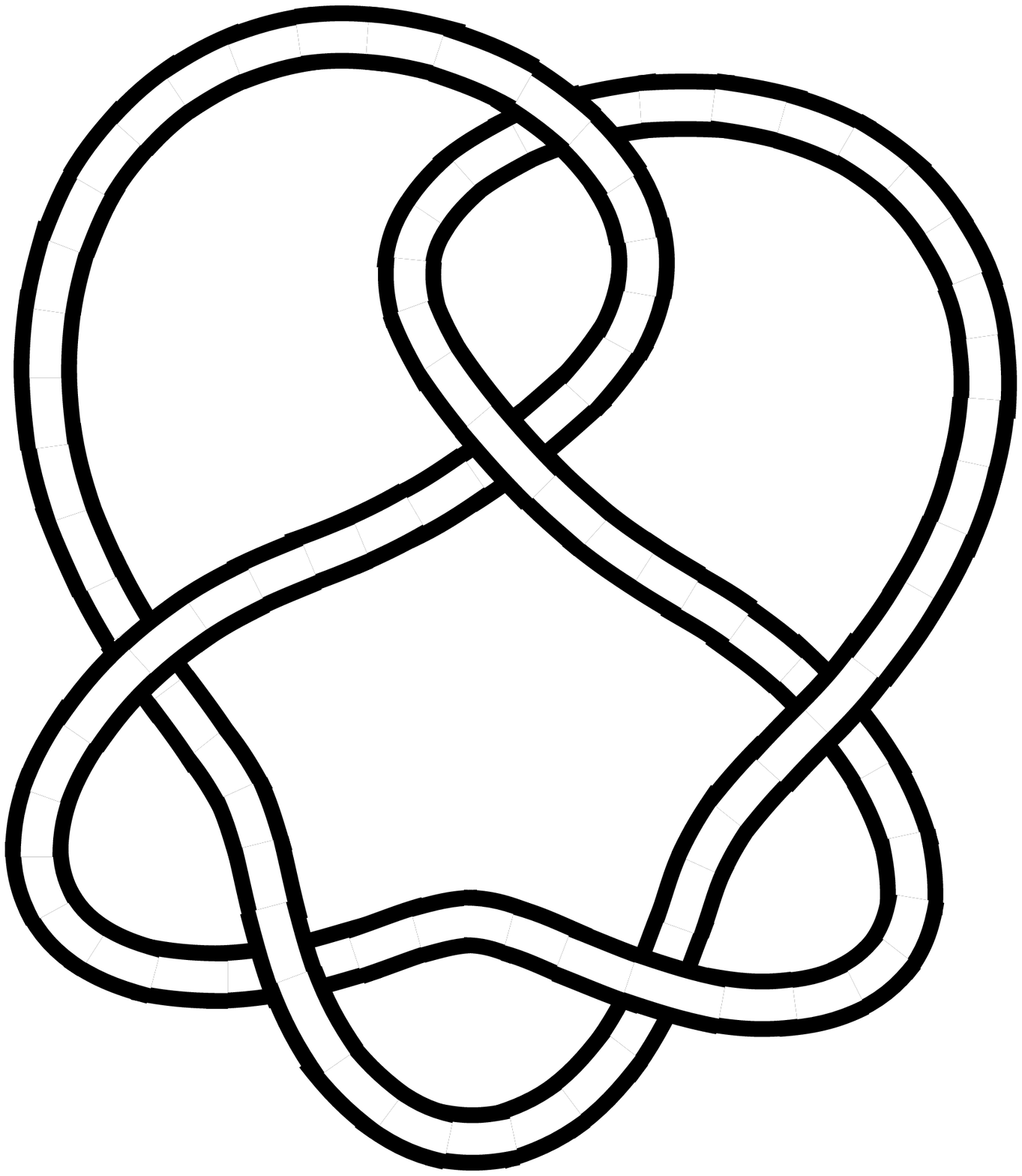}}& \neg $\begin{array}{l}\parbox{13mm}{$6_{1}$;} 0; \quad 0\\[1.5pt] 1^{\underline{4}}_{\underline{8}}1^{\underline{2}}_{\underline{4}}1^{\underline{1}}_{\underline{2}}1^{1}_{2}\neg \end{array}$ \neg& 
\parbox[c]{12mm}{     \includegraphics[height=12mm,width=12mm,keepaspectratio]{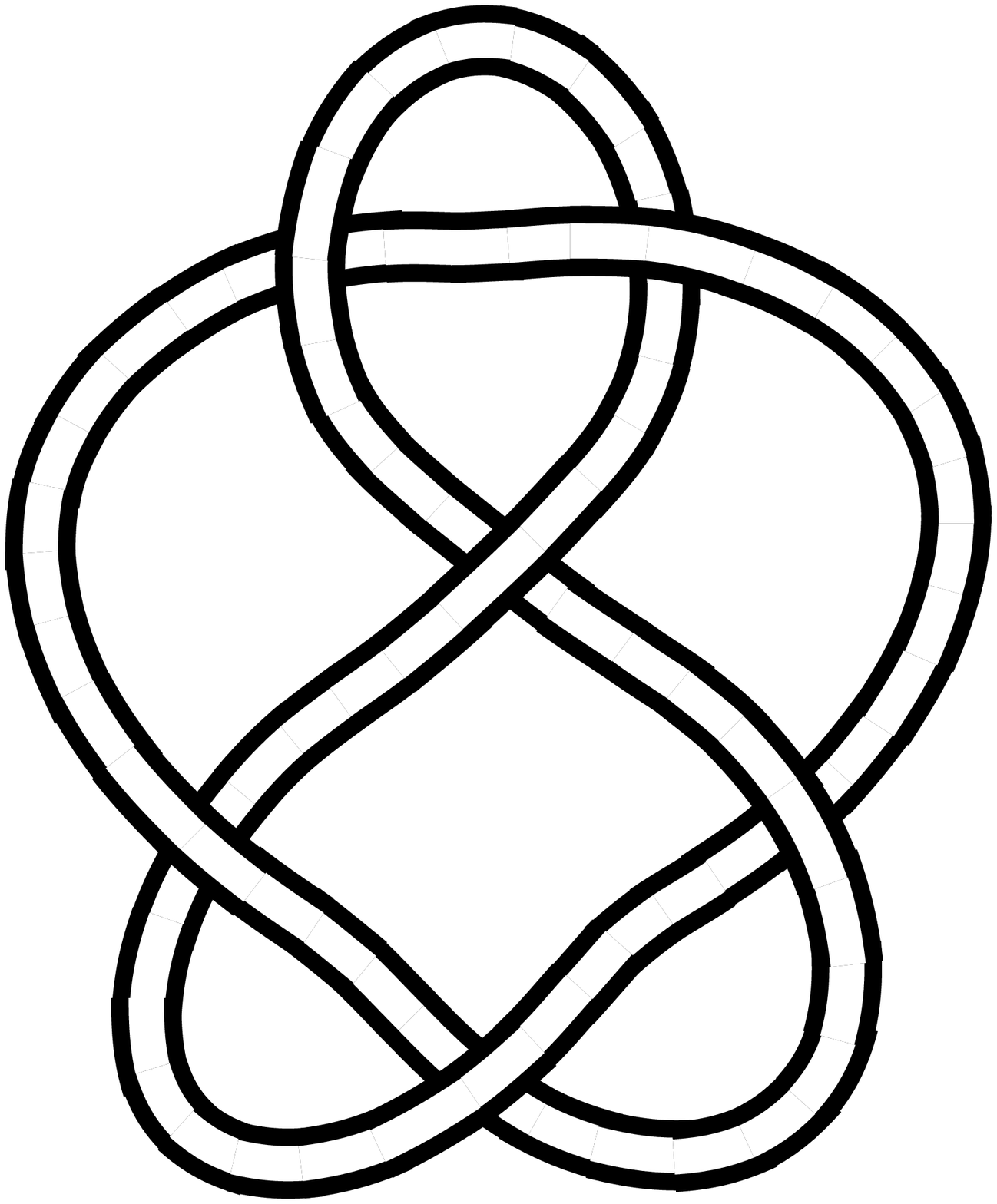}}& \neg $\begin{array}{l}\parbox{13mm}{$6_{2}$;} \underline{2}; \quad \underline{2}\\[1.5pt] 1^{\underline{4}}_{\underline{8}}1^{\underline{3}}_{\underline{6}}1^{\underline{2}}_{\underline{4}}1^{\underline{1}}_{\underline{2}}1^{1}_{2}\neg \end{array}$ \neg\\
\hline

\end{longtable}

\end{center}

\newpage

\subsection{$\Kh'(L)$ for prime knots with 11 crossings} This data is
available as a 20-page appendix to this paper (titled ``Khovanov's
invariant for 11 crossing prime knots'') and in computer readable
format from~\cite{cite:This}.

\subsection{New separation results} \label{subsec:NewSeparations}
Following is the complete list of pairs of prime knots with up to 11
crossings whose Jones polynomials are equal but whose rational Khovanov
invariants are different:
$(4_1, 11^n_{19})$,
$(5_1, 10_{132})$,
$(5_2, \overline{11^n_{57}})$,
$(7_2, \overline{11^n_{88}})$,
$(8_1, \overline{11^n_{70}})$,
$(9_2, \overline{11^n_{13}})$,
$(9_{42}, \overline{9_{42}})$,
$(9_{43}, 11^n_{12})$,
$(10_{125}, \overline{10_{125}})$,
$(10_{130}, \overline{11^n_{61}})$,
$(10_{133}, \overline{11^n_{27}})$,
$(10_{136}, 11^n_{92})$,
$(11^n_{24}, \overline{11^n_{24}})$,
$(11^n_{28}, 11^n_{64})$,
$(11^n_{50}, \overline{11^n_{133}})$,
$(11^n_{79}, \overline{11^n_{138}})$,
$(11^n_{82}, \overline{11^n_{82}})$,
$(11^n_{132}, \overline{11^n_{133}})$.

\subsection{$\Kh(L)$ for links with up to 11 crossings} For links with
more than one components, we have computed $\Kh(L)$ (not $\Kh'(L)$,
which does not make sense) for $L$ with up to 11 crossings. The results
are available as a 16 page appendix to this paper (up to 10 crossings)
and as a 26 page appendix (11 crossings) and in computer readable
format from~\cite{cite:This}. Below we only display the results for
links with up to 6 crossings. The same acknowledgement and warning of
the previous section still applies:


\vspace{5mm}

\begin{center} \scriptsize
\begin{longtable}{|ccl|ccl|}
  \hline
    $n^c_k$ & $L$ & \Kh(L) & $n^c_k$ & $L$ & \Kh(L) \\
  \hline \endhead
  \hline \endfoot

$2^2_1$ &     \hspace{-1.5mm}\parbox[c]{12mm}{     \includegraphics[height=12mm,width=12mm,keepaspectratio]{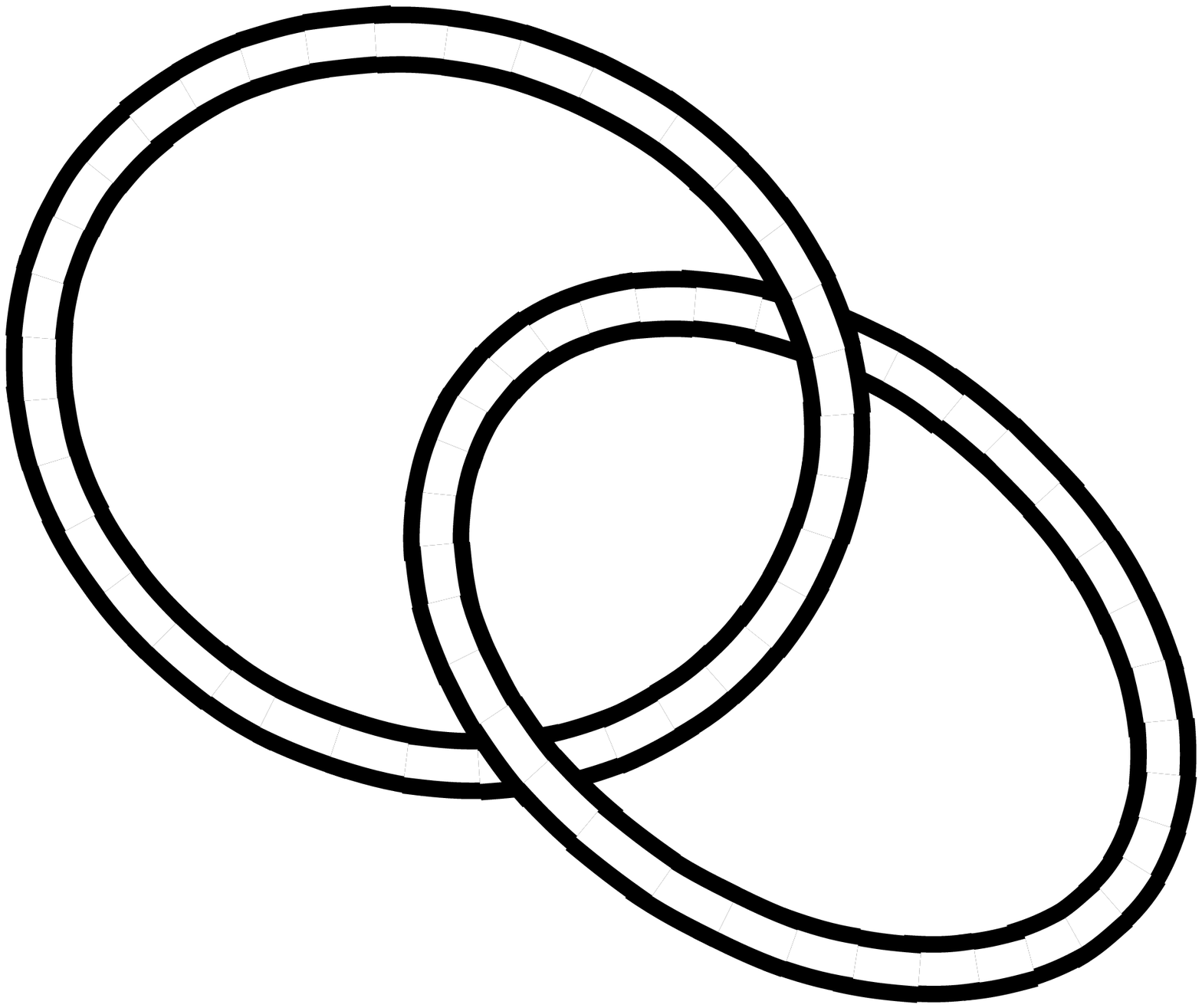}     } \hspace{-4mm} & $1^{0}_{0}1^{0}_{2}1^{2}_{4}1^{2}_{6}$ & $4^2_1$ &     \hspace{-1.5mm}\parbox[c]{12mm}{     \includegraphics[height=12mm,width=12mm,keepaspectratio]{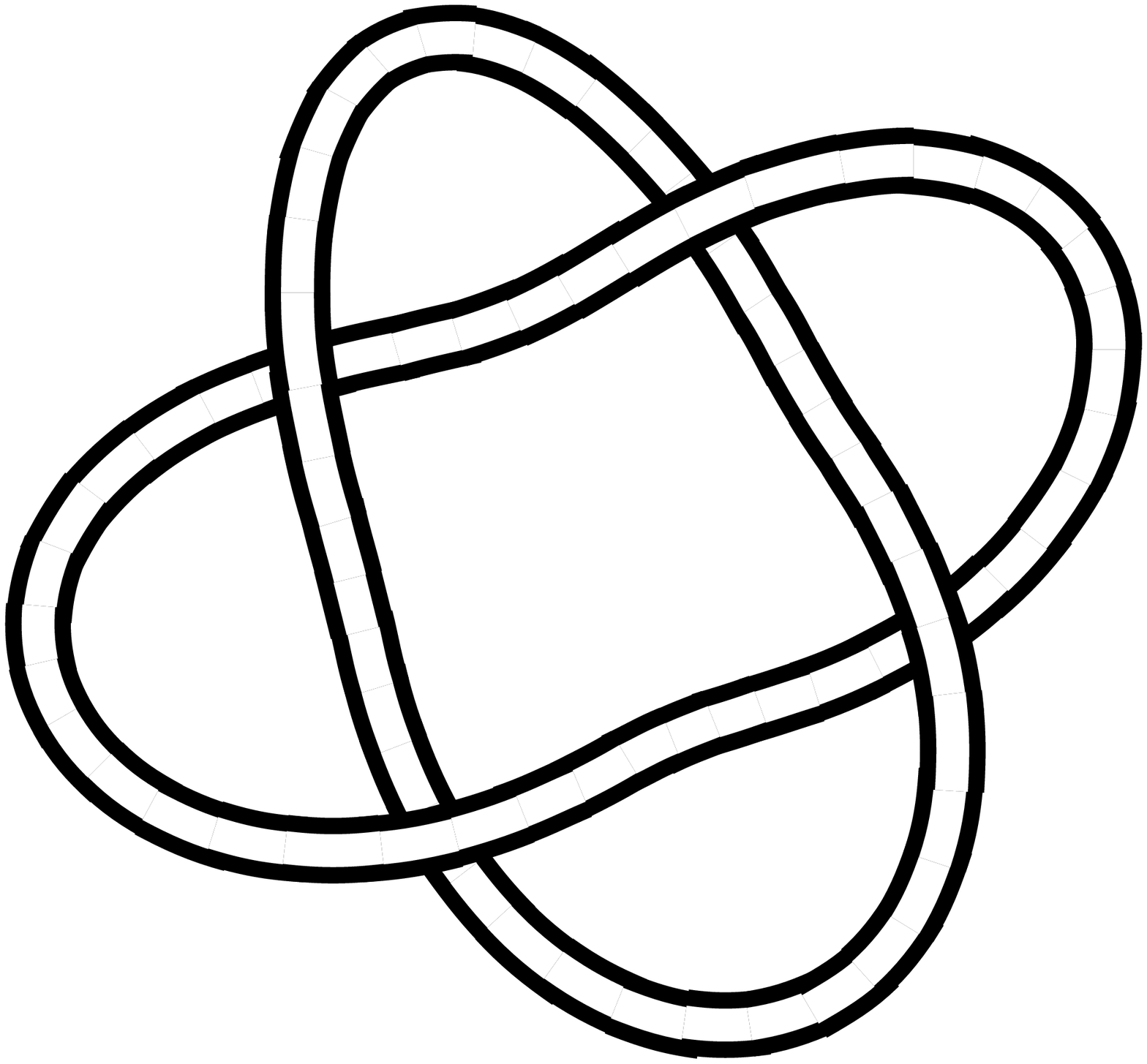}     } \hspace{-4mm} & $1^{\underline{4}}_{\underline{12}}1^{\underline{4}}_{\underline{10}}1^{\underline{3}}_{\underline{10}}1^{\underline{2}}_{\underline{6}}1^{0}_{\underline{4}}1^{0}_{\underline{2}}$ \\ \hline
$5^2_1$ &     \hspace{-1.5mm}\parbox[c]{12mm}{     \includegraphics[height=12mm,width=12mm,keepaspectratio]{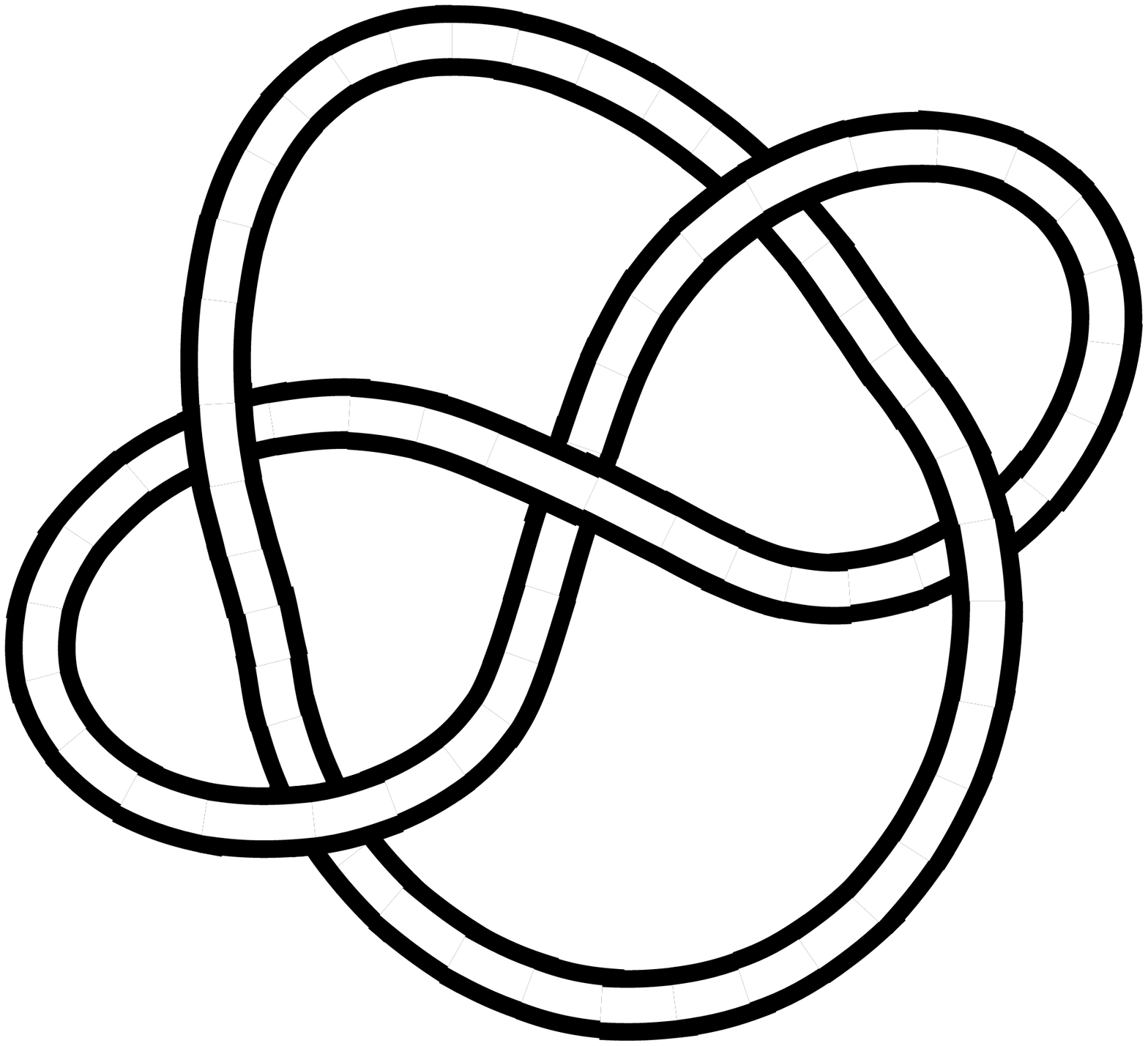}     } \hspace{-4mm} & $1^{\underline{3}}_{\underline{8}}1^{\underline{2}}_{\underline{6}}1^{\underline{2}}_{\underline{4}}1^{\underline{1}}_{\underline{2}}2^{0}_{\underline{2}}2^{0}_{0}1^{1}_{0}1^{2}_{4}$ & $6^2_1$ &     \hspace{-1.5mm}\parbox[c]{12mm}{     \includegraphics[height=12mm,width=12mm,keepaspectratio]{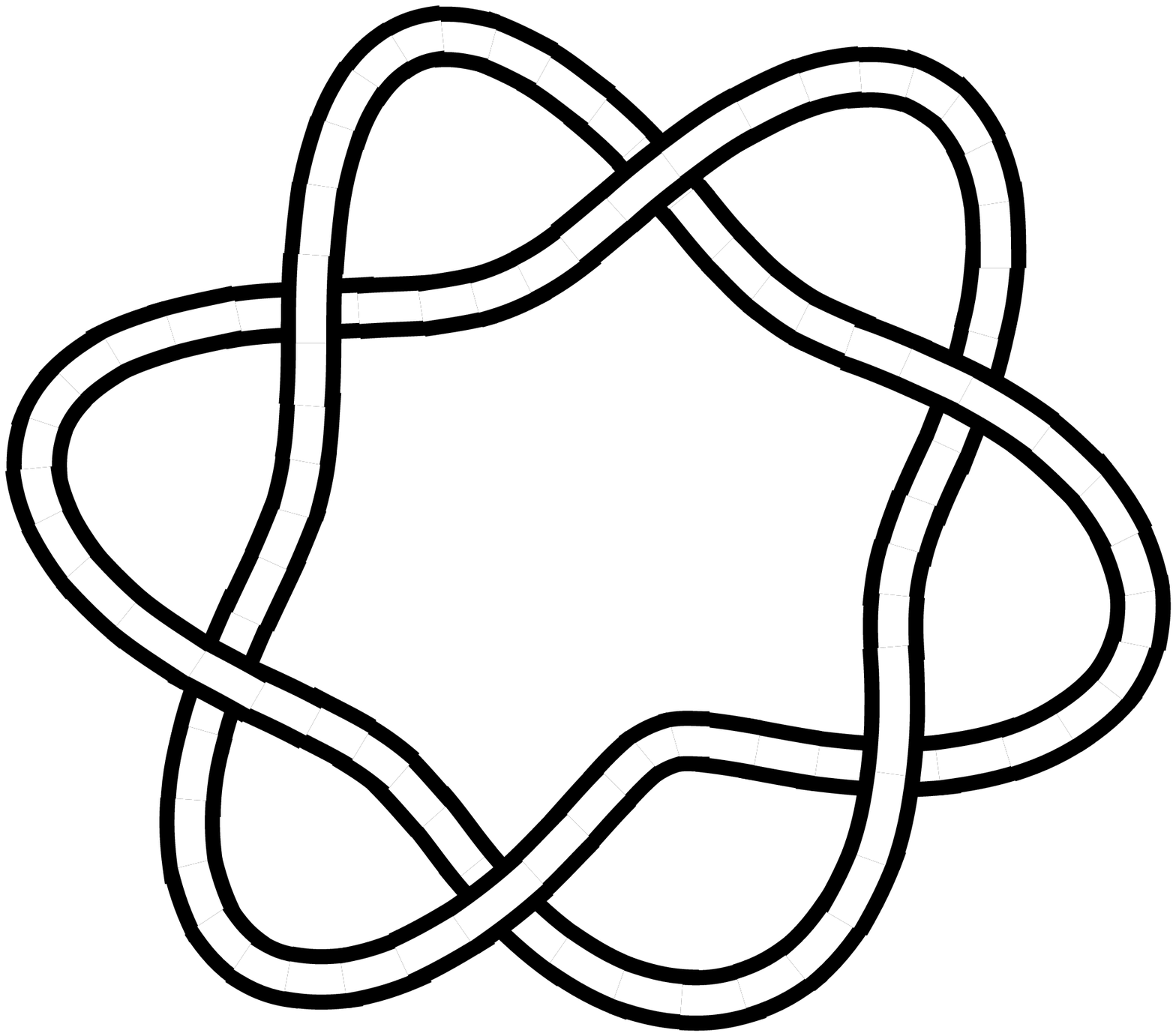}     } \hspace{-4mm} & $1^{\underline{6}}_{\underline{18}}1^{\underline{6}}_{\underline{16}}1^{\underline{5}}_{\underline{16}}1^{\underline{4}}_{\underline{12}}1^{\underline{3}}_{\underline{12}}1^{\underline{2}}_{\underline{8}}1^{0}_{\underline{6}}1^{0}_{\underline{4}}$ \\ \hline
$6^2_2$ &     \hspace{-1.5mm}\parbox[c]{12mm}{     \includegraphics[height=12mm,width=12mm,keepaspectratio]{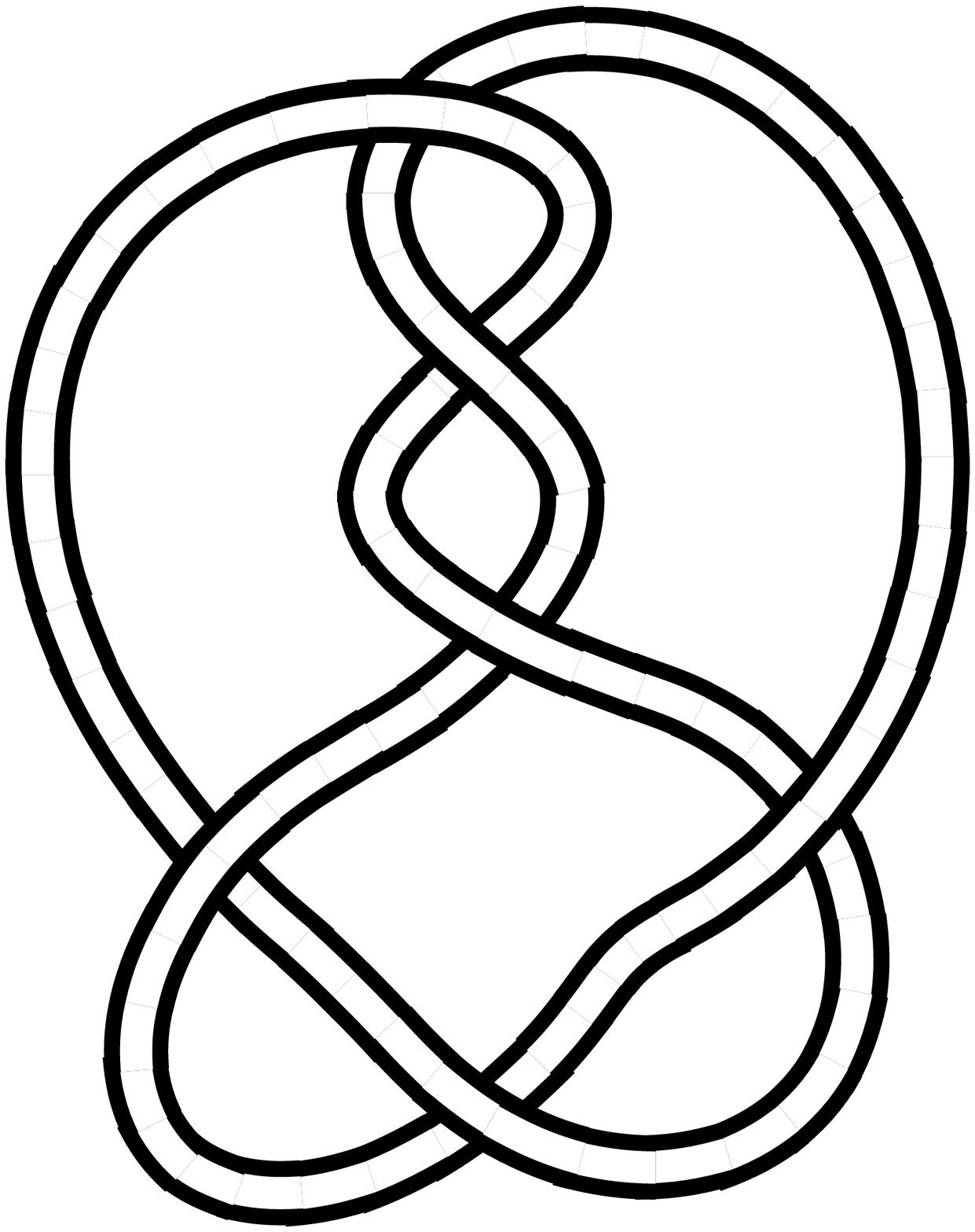}     } \hspace{-4mm} & $1^{0}_{2}1^{0}_{4}1^{1}_{4}1^{2}_{6}1^{2}_{8}1^{3}_{8}1^{3}_{10}1^{4}_{10}1^{4}_{12}1^{5}_{14}1^{6}_{14}1^{6}_{16}$ & $6^2_3$ &     \hspace{-1.5mm}\parbox[c]{12mm}{     \includegraphics[height=12mm,width=12mm,keepaspectratio]{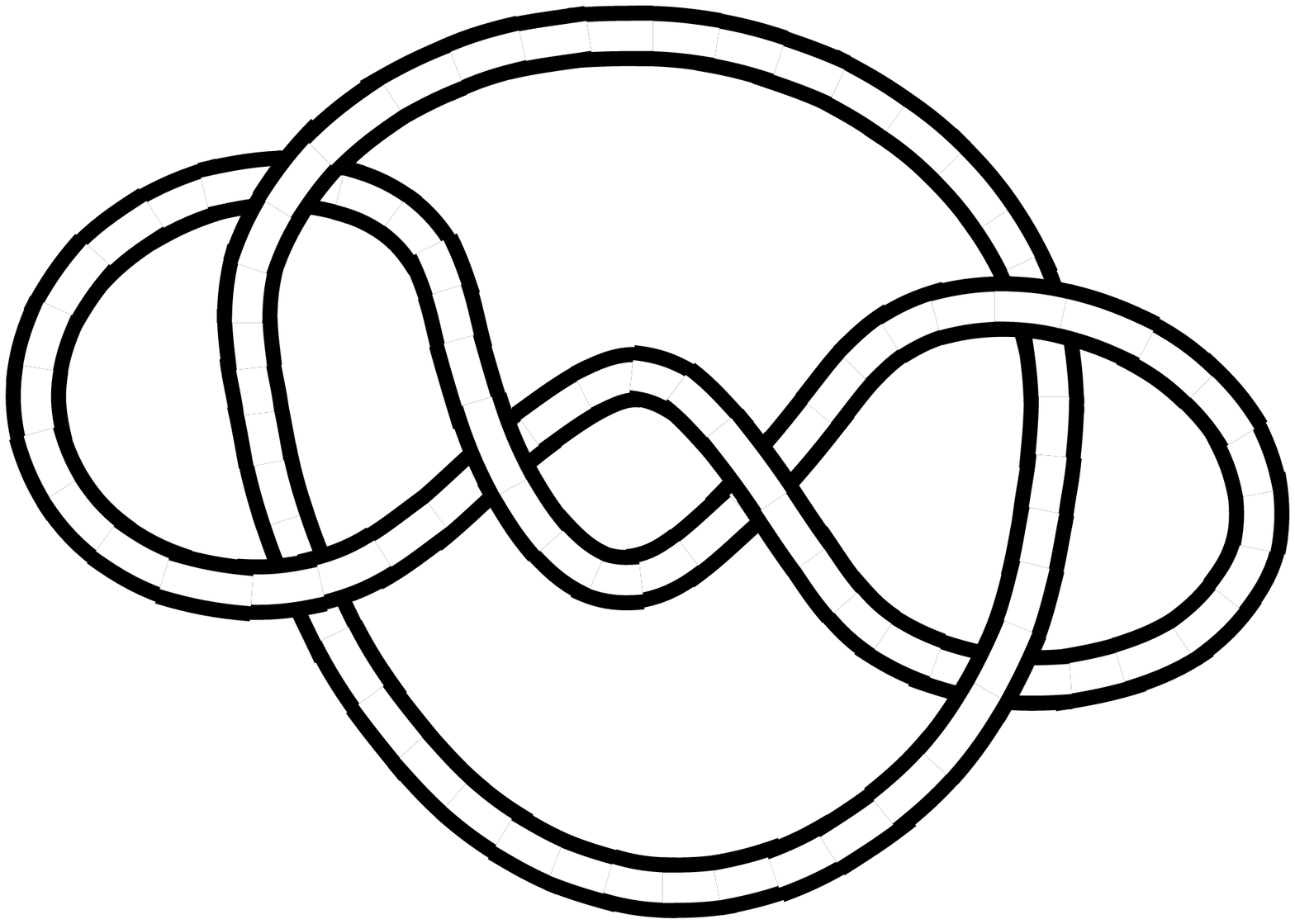}     } \hspace{-4mm} & $1^{\underline{6}}_{\underline{16}}1^{\underline{5}}_{\underline{14}}1^{\underline{5}}_{\underline{12}}1^{\underline{4}}_{\underline{12}}2^{\underline{4}}_{\underline{10}}2^{\underline{3}}_{\underline{10}}1^{\underline{2}}_{\underline{8}}2^{\underline{2}}_{\underline{6}}1^{\underline{1}}_{\underline{4}}1^{0}_{\underline{4}}1^{0}_{\underline{2}}$ \\ \hline
$6^3_1$ &     \hspace{-1.5mm}\parbox[c]{12mm}{     \includegraphics[height=12mm,width=12mm,keepaspectratio]{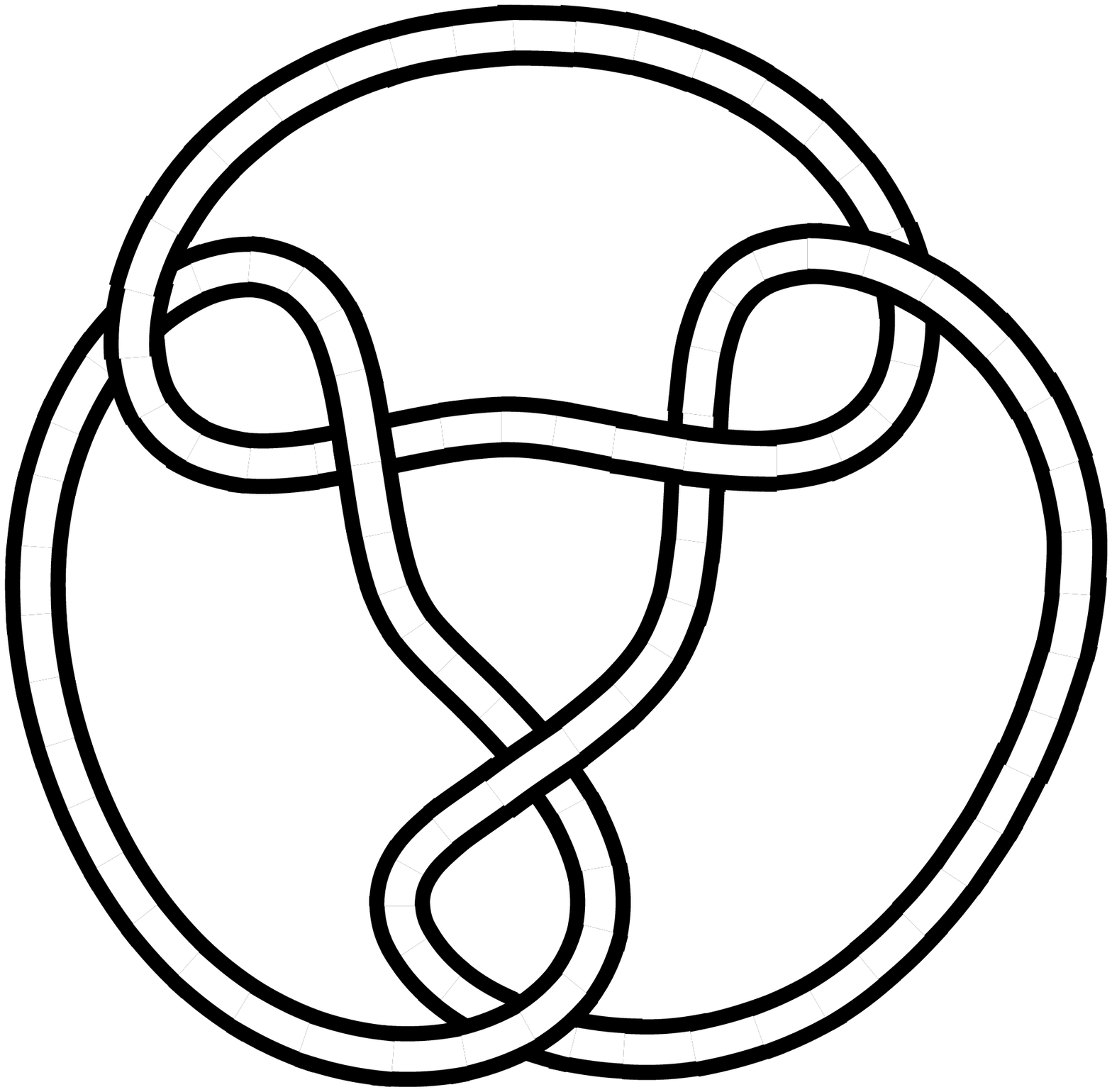}     } \hspace{-4mm} & $1^{\underline{6}}_{\underline{15}}1^{\underline{5}}_{\underline{11}}3^{\underline{4}}_{\underline{11}}3^{\underline{4}}_{\underline{9}}1^{\underline{3}}_{\underline{9}}2^{\underline{2}}_{\underline{7}}1^{\underline{2}}_{\underline{5}}2^{\underline{1}}_{\underline{3}}1^{0}_{\underline{3}}1^{0}_{\underline{1}}$ & $6^3_2$ &     \hspace{-1.5mm}\parbox[c]{12mm}{     \includegraphics[height=12mm,width=12mm,keepaspectratio]{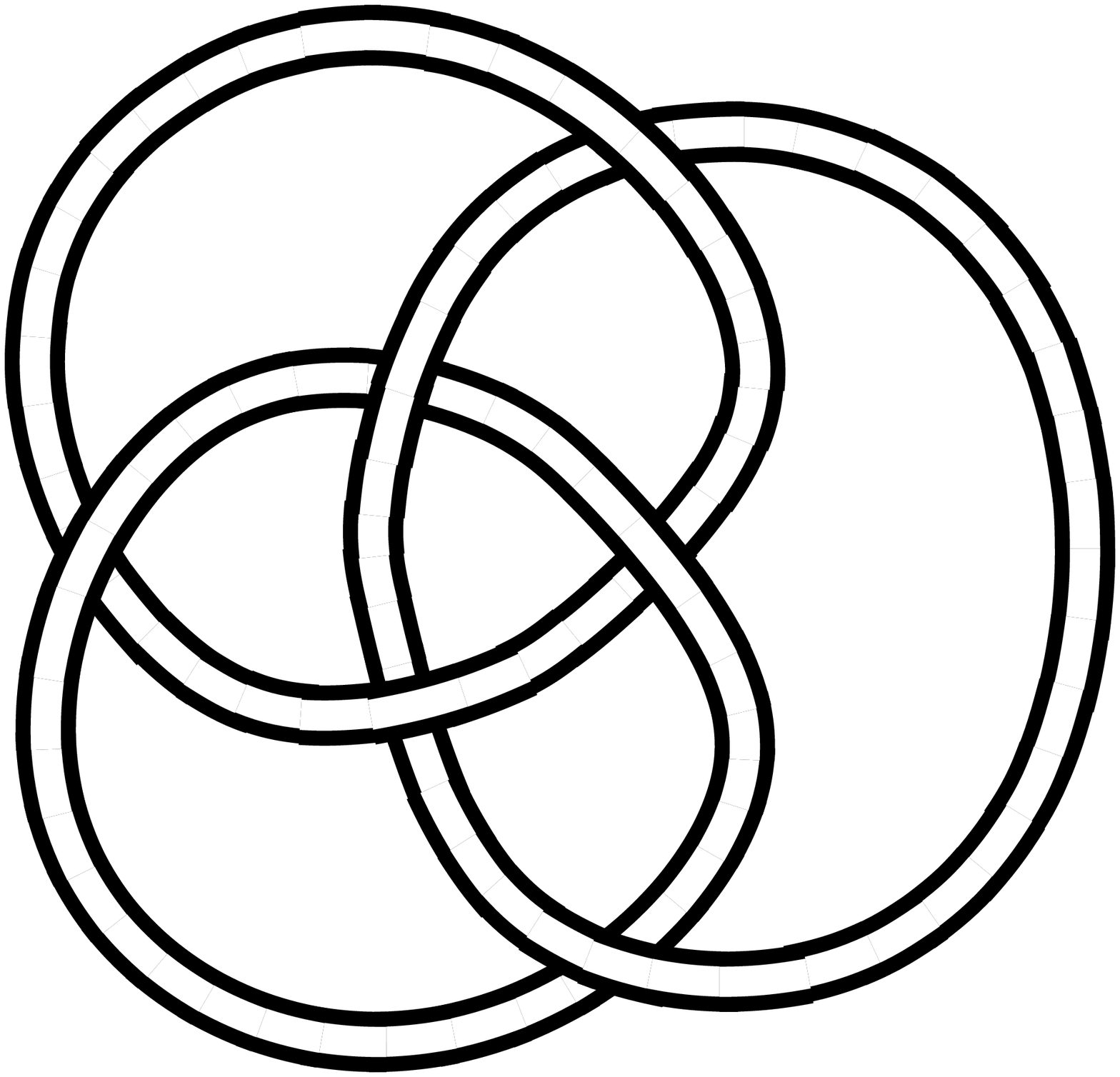}     } \hspace{-4mm} & $1^{\underline{3}}_{\underline{7}}2^{\underline{2}}_{\underline{5}}1^{\underline{2}}_{\underline{3}}2^{\underline{1}}_{\underline{1}}4^{0}_{\underline{1}}4^{0}_{1}2^{1}_{1}1^{2}_{3}2^{2}_{5}1^{3}_{7}$ \\ \hline
$6^3_3$ &     \hspace{-1.5mm}\parbox[c]{12mm}{     \includegraphics[height=12mm,width=12mm,keepaspectratio]{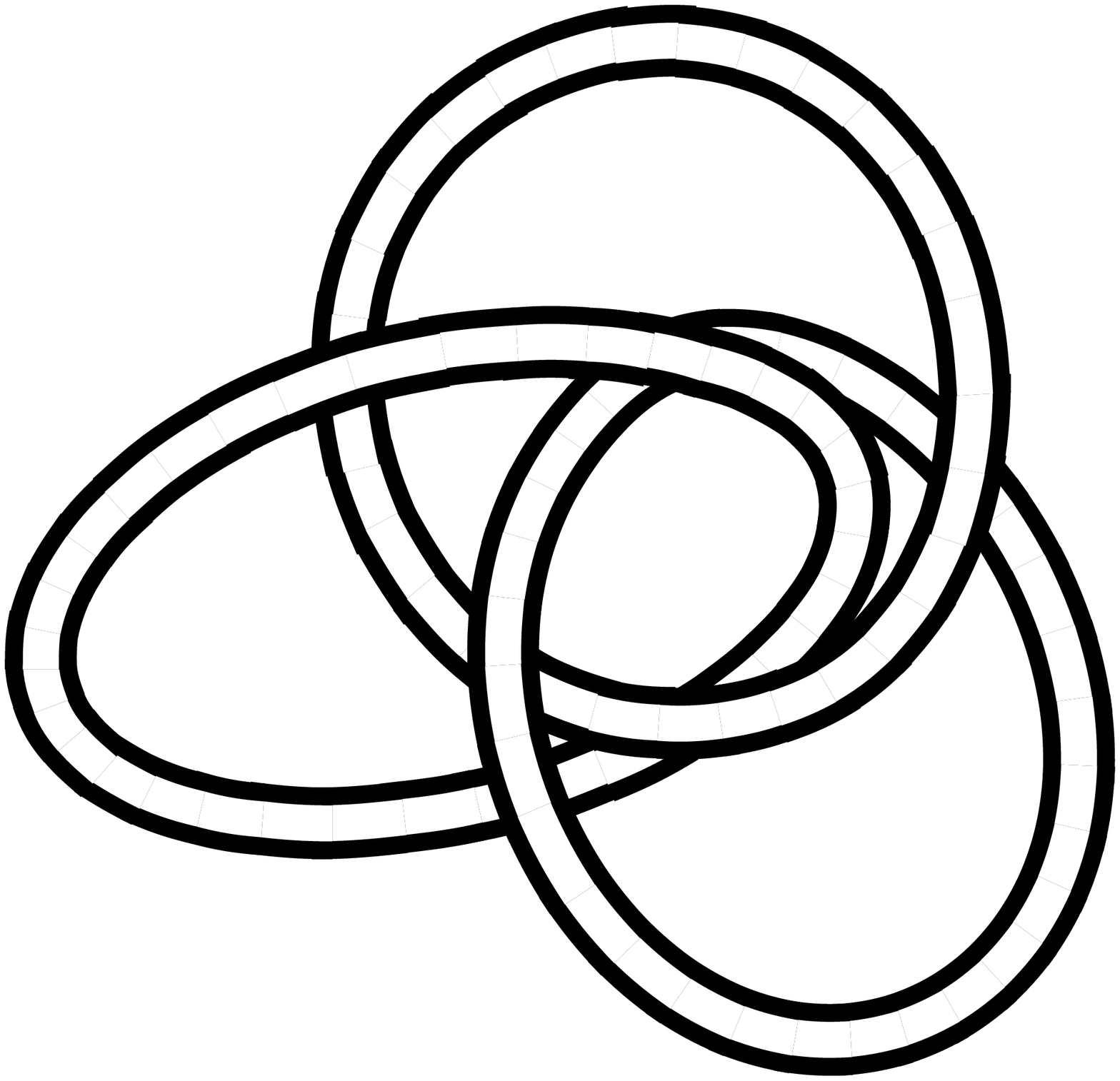}     } \hspace{-4mm} & $1^{0}_{3}1^{0}_{5}1^{2}_{7}1^{3}_{11}1^{4}_{9}3^{4}_{11}2^{4}_{13}$ & & & \\ \hline

\end{longtable}
\end{center}

\Addresses\recd

\begin{figure}[p]
\centering{\fbox{\resizebox*{!}{7.2in}{\parbox{7.2in}{
  \input qrg.tex
}}}}
\caption{
  A quick reference guide -- cut, fold neatly and place in your wallet.
}
\end{figure}

\end{document}
\endinput